\documentclass[11pt]{article} 
\usepackage{amssymb} 
\usepackage{cite}
\usepackage{epsfig} 
\usepackage{rotate}

\newcommand{\beq}{\begin{equation}} 
\newcommand{\eeq}{\end{equation}}
\newcommand{\Leq}[1]{\label{#1}\end{equation}}
\newcommand{\bdm}{\begin{displaymath}} 
\newcommand{\edm}{\end{displaymath}}

\newtheorem{theorem}{Theorem} 
\newtheorem{corollary}[theorem]{Corollary}
\newtheorem{lemma}[theorem]{Lemma}

\def\Z{{\bf Z}}
\def\R{\mathbb R}
\def\Z{\mathbb Z}
\def\M{{\cal M}}
\def\mod#1{\,({\rm mod\ }#1) }
\def\proof{\medskip\noindent P{\sc roof.}\hskip 10pt}
\def\endproof{\hfill\vbox{\hrule \hbox{\vrule\kern4pt\vbox{\kern4pt
\kern4pt}\kern4pt\vrule}\hrule}\bigskip}
\def\iff{\Leftrightarrow}

\title{Approach to a rational rotation number\\ 
in a piecewise isometric system} \author{J. H. Lowenstein and F. Vivaldi\dag } 
\date{\it\small 
Dept.~of Physics, New York University, 2 Washington Place, New York, NY 10003, USA
\\ 
\dag School of Mathematical Sciences, Queen Mary, University of London, London E1 4NS, UK 
}

\begin{document} 
\renewcommand{\labelenumi}{$(\roman{enumi})$} 
\maketitle
\begin{abstract} We study a parametric family of piecewise rotations of 
the torus, in the limit in which the rotation number approaches the rational
value 1/4. There is a region of positive measure where the discontinuity 
set becomes dense in the limit; we prove that in this region the area
occupied by stable periodic orbits remains positive. The main device is 
the construction of an induced map on a domain with vanishing measure;
this map is the product of two involutions, and each involution
preserves all its atoms. Dynamically, the composition of these involutions 
represents linking together two {\it sector maps;} this dynamical system
features an orderly array of stable periodic orbits having a smooth
parameter dependence, plus irregular contributions which become 
negligible in the limit.
\end{abstract}
\vspace*{50pt}
\centerline{\today} 

\section{Introduction}

\begin{figure}[t]
\hfil\epsfig{file=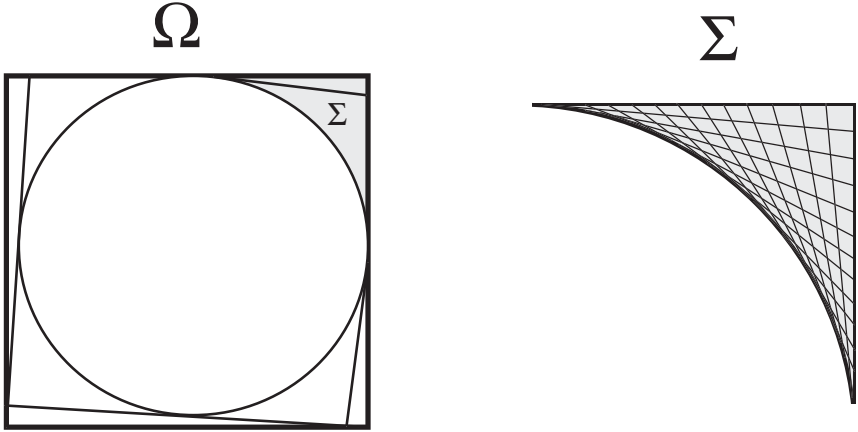,width=12cm}\hfil
\caption{\label{fig:BasicDiscontinuitySet}\rm\small
Left: Mismatch between the action of $F$ on $\Omega$ and a rotation by $\pi/4$,
for a parameter value close to zero. Shown are the first few images of the 
boundary of $\Omega$, which are tangent to the central island.
Right: Dominant features of the discontinuity set in the North-East 
sector $\Sigma$.
}
\end{figure}

This paper is devoted to the study of the action of the matrix
\begin{equation}\label{eq:Matrix} 
C=\left(\matrix{\lambda & -1\cr 1 & 0}\right)\hskip 30pt 
\lambda=2\cos(2\pi\rho) \end{equation} 
on the unit square, in the limit $\lambda\to 0$. 
This deceptively simple dynamical system has a surprisingly complicated 
behaviour. In figure \ref{fig:BasicDiscontinuitySet} we display
the phase portrait in the North-East corner of the unit square, for 
$\lambda=2^{-6}$. The picture suggests the existence of an infinite 
hierarchy of stable islands, immersed in a `pseudo-chaotic' sea.
 
More precisely, we let $\Omega=[0,1)^2$ and we define the map
\begin{equation}\label{eq:Map}
F:\Omega\to\Omega
\hskip 40pt
(x,y)\mapsto (\lambda x -y+\iota(x,y),x)
\hskip 30pt
\iota(x,y)=-\lfloor \lambda x-y \rfloor,
\end{equation}
where $\lfloor\cdot\rfloor$ denotes the floor function.
This map is area-preserving and has time-reversal symmetry \cite{LambRoberts}
\begin{equation}\label{eq:Symmetry}
F^{-1}=G\circ F\circ G^{-1}\hskip 40pt G:\Omega\to\Omega\qquad (x,y)\mapsto (y,x).
\end{equation}
Occasionally, we will regard $F$ as a map of the torus, rather than of the square. 
It must be understood however, that the function $\iota$ in (\ref{eq:Map}) is 
defined over $\Omega$, not over $\R^2/\Z^2$.

In the specified parameter range, the map $F$ is non-ergodic, and has zero topological 
entropy \cite{Buzzi}. It is linearly conjugate to a piecewise rotation on a rhombus 
with rotation number $\rho$, and indeed it's a {\it piecewise isometry\/} with respect 
to the metric induced by the inner product
\begin{equation}\label{eq:Metric}
{\cal Q}(U,V)=U_xV_x+U_yV_y-\frac{\lambda}{2}\left(UxVy+U_yV_x\right).
\end{equation}
The quantity ${\cal Q}(U,U)$ is the invariant quadratic form of the matrix $C$ in (\ref{eq:Matrix}).
For $\lambda=0$, this metric reduces to the ordinary Euclidean metric.
Throughout this paper, the term isometry will always refer to the metric (\ref{eq:Metric}).
Piecewise isometries are dynamical systems that generalize to higher dimensions 
the construct of interval exchange maps.
Recently, these systems have attracted a great deal of attention, both theoretically 
\cite{Goetz:96,Goetz:98,Goetz:00,GoetzPoggiaspalla,AdlerEtAl,
BruinEtAl,Kahng,KouptsovLowensteinVivaldi,AshwinGoetz,VowdenVowden},
and in applications 
\cite{ChuaLin:88,ChuaLin:90b,WuChua,Davies,AshwinEtAl,
LowensteinHatjispyrosVivaldi,LowensteinVivaldi:00}.  

\begin{figure}
\hfil\hfil\epsfig{file=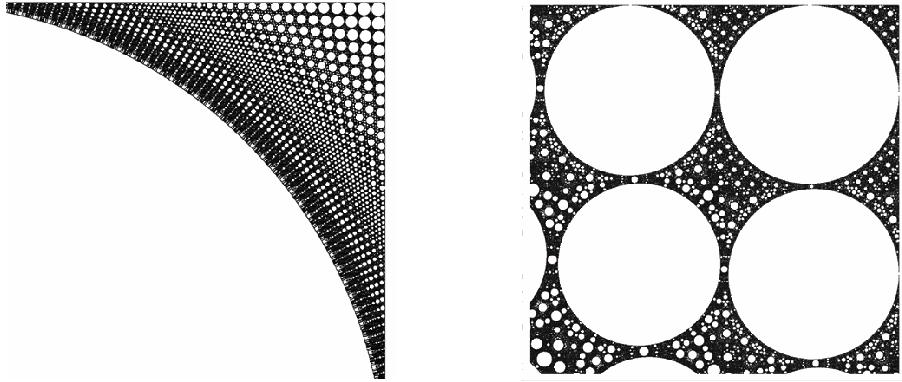,width=12cm}\hfil
\caption{\label{fig:NorthEastSector}\rm\small
Left: Detail of the phase portrait of the map $F$ for $\lambda=2^{-6}$,
corresponding to an irrational value of the rotation number $\rho$, 
close to $1/4$. A portion of the large central island appears in 
the SW corner of the picture. Right: Magnified view of the region near the NE corner.
}
\end{figure}

In the analysis of these systems, a decisive simplifying factor is the presence 
of self-similarity, which allows a satisfactory ---occasionally complete--- 
characterization of the dynamics \cite{Poggiaspalla}.
Self-similarity has been invariably found in maps with rational rotation 
number (the parameter $\rho$ in (\ref{eq:Matrix})), although the
occurrence of scaling in this case has never been justified by a general theory.
For rational rotations, the stable regions in phase space ---the ellipses of 
figure \ref{fig:NorthEastSector}--- become convex polygons, and the system 
parameters (such as the quantity $\lambda$ in (\ref{eq:Matrix})) are algebraic numbers.
Much of recent research has been devoted to quadratic parameter values, plus some 
scattered results for cubic parameters 
\cite{GoetzPoggiaspalla,LowensteinKouptsovVivaldi,Kahng:02,Poggiaspalla,Lowenstein}. 
Rational rotation numbers with prime denominator were considered in \cite{Goetz:05} 
from a ring-theoretic angle, in a rather general setting.
In all cases in which computations have been performed, the complement of the 
cells, namely the closure of the so-called {\it discontinuity set\/}
(the dark area in figure \ref{fig:NorthEastSector}), has been found to have 
zero Lebesgue measure.
However, no general result has been established in this direction.

The case of {\it irrational rotations\/} ---the generic one--- has
stubbornly resisted attack, largely due to an apparent lack of self-similarity. 
In 1997, P.~Ashwin conjectured that for irrational rotations, the 
complement of the set of elliptic islands has positive Lebesgue measure,
and this measure depends continuously on the parameter \cite{Ashwin}.
Subsequently, the semi-continuity of the measure at irrational rotation
numbers was established rigorously \cite{Goetz:01}; however, the important 
question of positivity of measure remains unresolved, and constitutes the 
motivation for the present work. Ashwin's conjecture is analogous to
the well-known conjecture of the positivity of measure of chaotic 
motions in Hamiltonian systems with divided phase space.

To identify elusive instances of scaling behaviour, in this paper
we examine irrational rotations in the vicinity of a rational one.
For values of the parameter $\lambda$ approaching zero, the map $F$
approaches a rotation of the unit square about its centre, by the angle $\pi/2$ 
(figure \ref{fig:BasicDiscontinuitySet}). To a first approximation, the 
mismatch between the actual and the limit dynamics causes the images 
of the boundary of the square (the discontinuity set) slowly to build 
a regular envelope of the central disk, by line segments. 
These segment intersect transversally, forming an orderly array of 
{\it pseudo-hyperbolic points} \cite{VivaldiLowenstein}.
It is not difficult to see that in the limit the pseudo-hyperbolic points 
become dense in the four corner sectors, and any stable island present in 
these sectors must be confined within the meshes of the grid, 
which are of vanishing size. 
On the other hand, at $\lambda=0$ the discontinuity set is just the 
boundary of the unit square, and all this complicated structure disappears. 

This paper is devoted to the study of the limit $\lambda\to 0$.
We have two main results. Firstly, we construct an induced mapping that 
describes in this limit the dynamics of $F$ outside the main island.
This map has interesting properties, summarized in the following theorem.

\bigskip\noindent
{\bf Theorem A.}{\sl\hskip 5pt
For $\lambda$ positive and sufficiently small, there is a mapping $L$, induced 
by $F$ on a domain $\Lambda$ of area $O(\lambda)$, which serves as a return map 
for all orbits outside the main island, apart from a set of measure $O(\lambda^2)$. 
The mapping $L$ is the composition of two involutions; each involution preserves 
all its atoms, apart from a set of irregular atoms of total measure $O(\lambda^3)$.
}

Theorem A is the summary of theorems \ref{theorem:Lin} and \ref{theorem:Lout},
which describe the structure of the involutions comprising $L$, 
and theorem \ref{theorem:Covering}, which estimates the measure of the 
orbits that pass through the so-called {\it regular atoms\/} of the map $L$.
This result leads us to consider the symmetric fixed points of $L$, which
form a distinguished family of stable periodic orbits of the map $F$
---the {\it regular sequence}. Our second result concerns the area of 
the corresponding stable regions.

\bigskip\noindent
{\bf Theorem B.}{\sl\hskip 5pt
As $\lambda\to 0^+$, the orbits of $F$ that correspond to the fixed points 
of $L$ form an infinite family of stable symmetric cycles, with the property 
that the total area of the associated islands approaches a positive limit.
}

The above theorem will appear in the text as theorems 
\ref{theorem:atomintersections} and \ref{theorem:Area},
which also give detailed information about these periodic orbits.
The `perturbative' approach, whereby irregular contributions of vanishing 
measure are neglected, is an essential element of our analysis, since the 
complexity of the neglected dynamics is considerable.

This paper is organised as follows.
In section \ref{section:Preliminaries} we provide the basic definitions and constructs. 
In section \ref{section:ReturnMap} we first define the induced map $L$,
and then provide a geometrical proof of the main part of theorem A 
(theorems \ref{theorem:Lin} and \ref{theorem:Lout}).
This proof is centred around the idea of {\it sector maps\/}, continuous 
transformations of a neighbourhood of an island, which are close to the 
identity. 
We will show that the dynamics of $F$ is the result of linking together
two sector maps (see figures \ref{fig:Sector} and \ref{fig:TwoSectors}).

In section \ref{section:QuantitativeResults}, we provide an algebraic proof 
of theorem \ref{theorem:Lout}, and also derive explicit formulae for the atoms 
of the return maps. These will be essential for the subsequent analysis.
In section \ref{section:PeriodicPoints} we prove theorem B.
We begin by constructing the regular sequence of periodic points, which shape the 
phase space outside the main island (figure \ref{fig:NorthEastSector}).
These orbits are parametrized by pairs of atoms of the two involutions 
comprising the mapping $L$ of theorem A. 
Then we use the formulae developed earlier to establish some asymptotic 
properties of the atoms of the aforementioned involutions.
With these formulae, we prove that, as $\lambda\to 0$, 
the regular sequence of stable periodic points of $L$ becomes infinite,
and the limiting area of the ellipses associated with the regular sequence 
approaches a positive limit (theorem \ref{theorem:Area}).
Our proof is constructive, and we derive an analytic formula for the total area.

In section \ref{section:Covering} we complete the proof of theorem A; we prove 
that the domain of definition of the induced map $L$ ---indeed the subset of 
it constituted by the regular atoms of $L$--- is a surface of section for all 
orbits, except for a set of measure $O(\lambda^2)$ (theorem \ref{theorem:Covering}).
Finally, in section \ref{section:Extensions} we briefly consider 
the extension of our results to the case of negative $\lambda$ and to
the dynamics of irregular atoms.
Several formulae and proofs are collected in an appendix.

\section{Preliminaries}
\label{section:Preliminaries}

The map $F$ consists of the action of the matrix (\ref{eq:Matrix}),
followed by a translation by $\iota$, given in equation (\ref{eq:Map}).
The quantity $\iota$ takes a finite set of integer values. In particular
\begin{equation} \label{eq:Alphabet}
\iota(x,y)\in\{0,1\} \qquad \mbox{if}\qquad 0 \leq \lambda < 1.
\end{equation}
We shall now assume that $\lambda$ belongs to this range, approaching 
zero from above; later on, we shall restrict the parameter further, 
to a smaller right neighbourhood of zero.
The case of negative parameters will be considered briefly in 
section \ref{section:Extensions}.

The unit square $\Omega$ is partitioned into two {\bf atoms\/} $\Omega_i$,
which are convex polygons given by the level sets of the function $\iota$
$$
\Omega_i=\iota^{-1}(j)\qquad j=0,1.
$$ 
Given a code\footnote{We use the symbol $\iota$ for both the coding function 
and the code.}
$\iota=(\iota_0,\iota_1,\ldots)$, with integer symbols $\iota_t$
taken from the alphabet (\ref{eq:Alphabet}), we consider the
set ${\cal C}(\iota)$ of the points $z\in\Omega$
for which $\iota(F^t(x,y))=\iota_t$, for $t=0,1,\ldots$.
These are the points whose images visit the atoms in the order
specified by the code.
If ${\cal C}(\iota)$ is non-empty, then it is called a {\bf cell.}
Generically, a cell is an open ellipse together with a subset of its 
boundary \cite[proposition 2]{Goetz:01}. 
In coordinates relative to their centre, these ellipses are similar
to the ellipse ${\cal Q}(x,y)=1$, cf.~equation (\ref{eq:Metric}). 
A cell can also be a convex
polygon (when the rotation number $\rho$ is rational), or a point.

The set of all images and pre-images of the boundary
of the atoms constitutes the {\bf discontinuity set\/} $\Gamma$
$$
\Gamma=\bigcup_{t=-\infty}^\infty F^t(\partial\Omega)
\hskip 40pt \partial\Omega=\bigcup_i\partial\Omega_i
$$
which consists of a countable set of segments.
For the purpose of generating $\Gamma$, the set $\partial\Omega$ may 
be replaced with any set that covers $\partial\Omega$ under iteration.
One verifies that the discontinuity set of the map (\ref{eq:Map}) is 
generated by the oriented segment $\gamma$, with first end-point 
$(0,0)$ (included) and second end-point $(0,1)$ (excluded). 
The map (\ref{eq:Map}), viewed as a map of the torus, 
is discontinuous only on $\gamma$, from the left. Indeed
$$
F(1-x,y)-F(x,y)=(\lambda +\iota(1-x,y)-\iota(x,y),\,1-2x)
 \equiv (\lambda,-2x)\mod{\Z^2}.
$$
Letting $x\to 0^+$, we see that $F$ has a {\it sliding singularity\/}
at $\gamma$ \cite{Kahng:09}, which causes a shift by $\lambda$ in the 
horizontal direction. By time-reversal symmetry, the map $F^{-1}$ 
has a siding singularity on $G(\gamma)$.

The generator $\gamma$ is oriented in such a way that the map $F$ be 
continuous on the segment's right side. The images and pre-images 
of $\gamma$ inherit an orientation with the same property.
In view of this, we'll say that a segment in the discontinuity set 
is {\bf glued\/} to any domain tangent to it on the right.
Furthermore, if the first and second end-points of a segment belong 
to two sets $A$ and $B$, respectively, we say that such a segment 
{\bf connects\/} $A$ to $B$.

Necessarily, any periodic orbit of $F$ will have a point whose 
cell is tangent to $\gamma$, from either the left or the right. 
For $t\in\Z$, the set $F^t(\gamma)$ is a collection of segments. 
Given a cell ${\cal E}$, the {\bf regular component\/}
of $F^t(\gamma)$ with respect to ${\cal E}$ is the set of segments 
in $F^t(\gamma)$ that are tangent to ${\cal E}$. 
A non-empty regular component consists of a single segment; indeed 
all components of the image of $\gamma$ must have the same orientation,
so that two tangent segments would have to lie at opposite
points of ${\cal E}$. However, for $\lambda\not=0$, there is no
cell touching $\gamma$ on opposite sides, and the image of a tangent point
is a single point. The same property holds also for images of the segment
$G(\gamma)$, and indeed of any other segment.
We will also speak of the regular component with respect to a cycle, 
because, given a cycle and an integer $t$, there is at most one ellipse in 
the cycle such that a component of $F^t(\gamma)$ is tangent to it. 
Regular components will be important in section 
\ref{section:ReturnMap}, in the construction of the induced map.

At $\lambda=0$ the dynamics is trivial. The atom $\Omega_0$ (a triangle)
collapses to the segment $G(\gamma)$, and the discontinuity set is
just $\Gamma=\gamma\cup G(\gamma)$.
All points in $\Gamma$ have period 2, except for a fixed point at the origin.
The rest of the space ---the interior of $\Omega$--- consists of a single 
square cell; its points have period four, except for a fixed point at 
the centre.

We are interested in periodic orbits for parameter values near zero.
A periodic orbit $(\overline{x_0,x_1,\ldots,x_{t-1}})$ corresponds to 
a periodic code $\iota=(\overline{\iota_0,\iota_1,\ldots,\iota_{t-1}})$. 
As the parameter varies, the points of the orbit move in phase space;
in \cite[theorem 2]{VivaldiLowenstein} it was shown that the coordinates
$x_j$ are rational functions of $\lambda$. 
Specifically, one has
\begin{equation}\label{eq:RationalFunctions}
x_j(\lambda,\iota)= \frac{{\cal X}_t(\lambda,\sigma^j(\iota))}{{\cal M}_t(\lambda)}
\qquad j=0,\ldots,t-1,
\end{equation}
where ${\cal X}_t$ and $\M_t$ are polynomials in $\Z[\lambda]$,
and $\sigma$ is the left shift map.
The denominator $\M_t$ is monic of degree $\lfloor (t+2)/2\rfloor$ 
and it depends only on the period; the numerator ${\cal X}_t$ has degree
$\lfloor (t+1)/2\rfloor$, and it depends on the orbit via the code $\iota$. 
This algebraic structure is common for periodic points of piecewise affine 
systems, see, e.g., \cite{BirdVivaldi,BosioVivaldi}.
In section \ref{section:PeriodicPoints} we will derive {\it ad hoc\/} 
expressions for these functions, tailored for the $\lambda\to 0$ regime.

The symbolic dynamics is far from being complete, so not all periodic 
codes $\iota$ correspond to an actual periodic orbit.
For this to be the case, the value of all functions (\ref{eq:RationalFunctions})
must belong to the half-open interval $[0,1)$. The function that is closest 
to the boundary of the unit interval determines the common size of all cells 
of the orbit corresponding to that code.
Accordingly, we define the {\bf radius\/} $r$ 
of a periodic orbit as
\begin{equation}\label{eq:Radius}
r(\lambda)=\min_{0\leq j\leq t-1} \{x_j(\lambda),\,1-x_j(\lambda)\}.
\end{equation}
The function $r$ is piecewise rational; besides the singularities
inherited from the functions $x_j$, the radius will typically have
points with discontinuous first derivative, due to a change of the 
index $j$ for which the minimum in (\ref{eq:Radius}) is attained.
An orbit exists if it has positive radius, and the radius becomes zero
at the bifurcation parameter values. An orbit of zero radius exists 
if and only if none of its rational functions assume the value 1,
which is the excluded point in the unit interval. On the torus,
this distinction becomes irrelevant.

Formulae of the type (\ref{eq:RationalFunctions}) also describe the
parameter dependence of pseudo-hyperbolic points \cite[theorem 6]{VivaldiLowenstein}.
These are the points that recur to the boundary of the atoms. More precisely,
a point in $\Omega$ is pseudo-hyperbolic if it maps to $\gamma$ in the forward
time direction, and to $G(\gamma)$, in the backward time direction. For an irrational 
rotation number, these points correspond to transversal intersections of two segments 
of the discontinuity set, which act like `pseudo-separatrices'. 
It can be shown that, generically, an orbit on the discontinuity set can contain 
only finitely many pseudo-hyperbolic points, the first of which is in $G(\gamma)$,
and the last in $\gamma$. 
A pseudo-hyperbolic point can be periodic only for parameter values which are 
algebraic of degree greater than one \cite[theorem 9]{VivaldiLowenstein}.
Thus, in general, a pseudo-hyperbolic sequence belongs to an infinite orbit;
such a sequence is determined by a finite symbolic dynamics, which describes the
itinerary between the endpoints on $G(\gamma)$ and $\gamma$.

Pseudo-hyperbolic points will occur as vertices of the atoms of the return map $L$,
to be defined in the next section. As with periodic orbits, we shall derive
explicit formulae for these points, appropriate for the $\lambda\to0$ limit.

\section{Return map}
\label{section:ReturnMap}

In this section we prove the main part of theorem A.
We will construct an induced mapping $L$ that describes the 
dynamics of $F$ in the limit $\lambda\to 0^+$.
This map is the composition of two involutions on a 
$\lambda$-dependent domain whose area vanishes in the limit
(theorems \ref{theorem:Lin} and \ref{theorem:Lout}).

From the time-reversal symmetry (\ref{eq:Symmetry}), we find that $F$
can be written as the composition of two orientation-reversing involutions
\begin{equation}\label{eq:HG}
F=H\circ G \qquad \mbox{with}\qquad H=F\circ G,\quad G^2=H^2={Id}
\end{equation}
with $H$ given by $(x,y)\mapsto (\{\lambda y -x\},\,y)$, where 
$\{\cdot\}$ denotes the fractional part.
Let ${\rm Fix}\, G$ and ${\rm Fix}\, H$ be the sets of
fixed points of the involutions $G$ and $H$, respectively.
The set ${\rm Fix}\, G$ is the segment $x=y$, independent of $\lambda$;
for $0\leq \lambda <1$, the set ${\rm Fix}\, H$ consists 
of two segments, with end-points $(0,0),\, (\lambda/2,1)$, and 
$(1/2,0),\, ((\lambda+1)/2,1)$, respectively. 

We begin by examining the rational functions of the main periodic orbits of $F$ 
(cf.~equation (\ref{eq:RationalFunctions}))
\begin{equation}\label{eq:BasicCycles}
\vcenter{\baselineskip 18pt\halign{
\quad 
 $#$ \hfil&\qquad
 $#$ \hfil&\qquad
 $\displaystyle #$ \hfil &\hskip 20pt 
 $\displaystyle #$ \hfil &\hskip 20pt 
 $\displaystyle #$ \hfil \quad
\cr
t &\mbox{code}  & \mbox{denominator} & \mbox{numerators}&\mbox{radius} \cr
\noalign{\vskip 4pt\hrule\vskip 5pt}
1&(\overline{1})  & 2-\lambda&\quad (1)&1-x_0(\lambda)\cr
2&(\overline{01}) & 4-\lambda^2&\quad (2,\lambda)&x_1(\lambda)\cr
}}
\end{equation}
These cycles are symmetric ($G$-invariant), and
the 2-cycle belongs to ${\rm Fix}\, H$.
At $\lambda=0$ the fixed point has radius $1/2$,
while the 2-cycle has zero radius, but it still exists 
on the discontinuity set:\/ $\{(1/2,0),(0,1/2)\}$.
It turns out that, in the $\lambda\to 0^+$ limit, the dynamics 
is dominated by these two periodic orbits,
which determine the coarse features of the phase space. 
The fine structure also depends on other periodic orbits, 
some of which will be analysed in section \ref{section:Extensions}.

The ellipse ${\cal E}$ bounding the large cell of the 1-cycle is tangent 
to the segment $\gamma$ at the point $(1,\tau)$, where $\tau=(1+\lambda)/2$.
By $G$-symmetry, ${\cal E}$ is also tangent to $\gamma_0=G(\gamma)$ at the 
point $T_0=(\tau,1)$. 
The open region delimited by the segments $\gamma$, $\gamma_0$, and 
the elliptical arc connecting $T_0$ to $G(T_0)$ clockwise will be 
called the {\bf 1-sector,} 
denoted by $\Sigma$. 
By construction, the set $\Sigma$ is invariant under $G$.
The 1-sector ---indeed a vanishingly small subset of it--- will 
serve as a surface of section for orbits outside the main island. 
Not all such orbits intersect $\Sigma$, for example, the 2-cycle
displayed in (\ref{eq:BasicCycles}), and the cycles mentioned in 
section \ref{section:Extensions}.
However, the missing orbits will be shown to have vanishing measure.
More important, the distinction between the orbits that enter $\Sigma$ 
and those that don't will be essential to our analysis.

For $\lambda$ near zero, the map $F^4$ is close to the identity,
and hence the set $\Sigma$ is `close to' being invariant under $F^4$.
There are however orbits that enter and leave $\Sigma$; to characterize
them, we define the sets
\begin{equation}\label{eq:InOut}
\Lambda=\Sigma\setminus F^4(\Sigma)
\hskip 40pt
\Xi=F^4(\Sigma)\setminus \Sigma
\end{equation}
where, by time-reversal symmetry, 
$$ 
G(\Lambda)=\Sigma\setminus F^{-4}(\Sigma)
\hskip 40pt
\Xi=F^4\circ G(\Lambda).
$$
Thus $G(\Lambda)$ is the set of points that leave $\Sigma$ under $F^4$.
Below, we will define two induced maps over the set $\Lambda$.

\begin{figure}[t]
\begin{center}
\epsfig{file=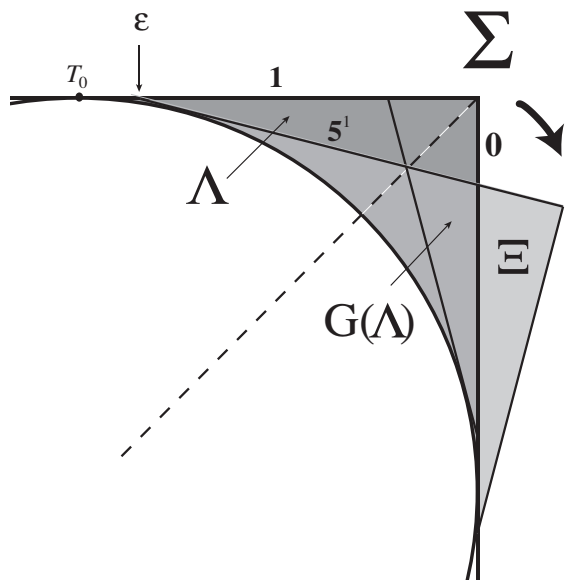,width=8cm}
\end{center}
\caption{\label{fig:Sector}
\small The sector $\Sigma$, and its image under $F^4$, corresponding to a 
small clockwise rotation. 
The darker region is the entry sub-domain $\Lambda$; the exit sub-domain 
$G(\Lambda)$ is the reflection  of $\Lambda$ with respect to the symmetry 
axis ${\rm Fix}\, G$ (the dashed line).
All points in $\Sigma$ will eventually map to the region
$\Xi=F^4\circ G (\Lambda)$, located on the opposite side of the
discontinuity line ${\bf 0}$.
}
\end{figure}

To determine the set $\Lambda$, we construct the first 
few images of the segment $\gamma$, which are listed in table 
\ref{table:gammaImages} below as pairs of end-points. 
We adopt the following conventions: \/
$i)$ an integer {\bf n} in boldface denotes the $n$-th iterate of the 
generator $\gamma$ (e.g., ${\bf 3}=F^3(\gamma)$); \/
$ii)$ segments crossing $\gamma$ have negative $x$-coordinate 
in the end-point on the left of $\gamma$; \/
$iii)$ if an end-point of ${\bf n}$ belongs to more than one 
set ${\bf m}$, we list only the value of $m$ corresponding 
to the lowest order iterate.

$$
\vcenter{\baselineskip 20pt\halign{
\quad $\displaystyle  #$ \hfil&\qquad
 $\displaystyle #$ \hfil&\quad
 $\displaystyle #$ \hfil \cr
                & \mbox{first} & \mbox{second} \cr
\noalign{\vskip -10pt}
\mbox{segment}  & \mbox{end-point} &\mbox{end-point} \cr
\noalign{\vskip 4pt\hrule\vskip 5pt}
\quad {\bf 0}  & (0,0) & (0,1)\cr
\quad {\bf 1}  & (1,0)\in {\bf 0} & (0,0)\in{\bf 0}\cr
\quad {\bf 2}  & (\lambda,1)\in {\bf 1} & (0,0)\in{\bf 0}\cr
\quad {\bf 3}  & (\lambda^2,\lambda)\in {\bf 2} & (1,0)\in {\bf 0}\cr
\quad {\bf 4}  & (-\lambda+\lambda^3,\lambda^2)\in {\bf 3}& (\lambda,1)\in {\bf 1}\cr
\quad {\bf 5}^1  & (\lambda-2\lambda^2+\lambda^4,1-\lambda+\lambda^3)\in {\bf 4}
               & (g-1,1)\in {\bf 1}\cr
\quad {\bf 5}^2  & (g',0)\in {\bf 1} & (\lambda^2,\lambda)\in {\bf 2}\cr
}}
$$

\begin{equation}\label{eq:ggprime}
g=\frac{1+2\lambda-\lambda^2-\lambda^3}{2-\lambda^2}
\hskip 40pt
g'=\frac{1-\lambda^2}{2-\lambda^2}.
\end{equation}
\begin{table}[h]
\caption{\label{table:gammaImages}
The first five iterates of the generator $\gamma$ of the discontinuity set.}
\end{table}

The segments ${\bf 0},\ldots,{\bf 4}$ are regular components with respect to
both the 1-cycle and the 2-cycle. 
Since $0\leq\lambda<1$, the segment ${\bf 4}$ crosses the discontinuity line 
$\gamma$, and so its image consists of two segments ${\bf 5}^1$ and ${\bf 5}^2$.
These are the images of the portions of ${\bf 4}$ lying to the left and to the 
right of $\gamma$, respectively, and they are regular components of ${\bf 5}$
with respect to the 1-cycle and the 2-cycle, respectively. 

We define approximate orthogonality relations.
If the angle between two segments ${\bf m}$ and ${\bf n}$ is equal to 
$0,\pi+O(\lambda)$, we say that ${\bf m}$ and ${\bf n}$ are 
{\bf quasi-parallel,} and write ${\bf m}\parallel {\bf n}$.
Likewise, two segments are {\bf quasi-perpendicular\/} (${\bf m}\perp{\bf n}$)
if the angle between them is $\pm\pi/2+O(\lambda)$. Plainly,
${\bf m}$ and ${\bf n}$ are quasi-parallel (perpendicular) if $m-n$
is even (odd). 

There are also incidence relations. We use the notation 
${\bf m}\dashv{\bf n}$ to indicate that one endpoint of ${\bf m}$ lies 
in the interior of ${\bf n}$ (and similarly for ${\bf m}\vdash{\bf n}$). 

The images of $\gamma$ have the following behaviour under symmetry
\begin{eqnarray}
G({\bf m})&=&G\circ F^m(\gamma)=F^{-m}\circ G(\gamma)=F^{-m}
 \circ F(\gamma)={\bf -m+1}\label{eq:Gm}\\
H({\bf m})&=&H\circ F^m(\gamma)=H\circ H\circ G\circ F^{m-1}(\gamma)
= G\circ F^{m-1}(\gamma)={\bf -m+2}.\label{eq:Hm}
\end{eqnarray}
We specify polygons/sectors by listing their boundary components between
angle brackets, arranged clockwise, e.g.,
$$
\Omega_0=\langle\,{\bf 0},{\bf 1},{\bf -1}\,\rangle
\hskip 40pt
\Sigma=\langle\,{\cal E},{\bf 1},{\bf 0}\,\rangle.
$$
Since $G$ and $H$ are orientation-reversing, mapping a polygon entails
mapping its sides, and reversing their order, e.g.,
$$
G(\Omega_0)=\langle\,G({\bf 0}),G({\bf -1}),G({\bf 1})\,\rangle
=\langle\,{\bf 1},{\bf 2},{\bf 0}\,\rangle.
$$

The segment ${\bf 5}^1$ connects ${\bf 4}$ to ${\bf 1}$,
and ${\bf 5}^1\cap \Sigma$ is the segment with end-points 
\begin{equation}\label{eq:gdoubleprime}
(1,g''),\,(g,1)\hskip 40pt 
g''=\frac{1}{1+\lambda-\lambda^2}\qquad
\end{equation}
where $g$ was given above.
As $\lambda$ approaches zero, $g\to 1/2$ and $g''\to 1$.
For our purpose it is sufficient to require that $g''<g$; 
so we replace the bounds in (\ref{eq:Alphabet}) by the more 
restrictive conditions
\begin{equation}\label{eq:lambdaRange}
0<\lambda <\lambda_+<1 
\qquad \mbox{where}\qquad
\lambda_+=2\cos(2\pi/9).
\end{equation}
(The number $\lambda_+$ is algebraic, being a root of the polynomial $x^3-3x+1$.)

The segment ${\bf 5}^1$ is tangent to the ellipse ${\cal E}$ at the point $T_1$;
hence ${\bf 5}^1$ decomposes $\Sigma$ into three regions (figure \ref{fig:Sector});
a large triangular sector $F^4(\Sigma)\cap\Sigma$, of area $O(1)$, 
a right triangle of area $O(\lambda)$, and a small triangular sector 
$\varepsilon$ of area $O(\lambda^3)$. 
The set $\Lambda$ is the union of the last two regions. 
We find
\begin{equation}\label{eq:DomainBoundaries}
\Lambda=\langle\,{\bf 0},{\bf 5}^1,{\bf 1}\rangle \cup\,\varepsilon
\hskip 40pt
\Xi=\langle\,{\bf 0},{\bf 5}^1,{\bf 4}\rangle.
\end{equation}
Our analysis will be perturbative, and in what follows we shall neglect
domains of area $O(\lambda^2)$. In particular, we shall omit $\varepsilon$ 
from all considerations, and represent $\Lambda$ as a triangle. 

We intend to study the dynamics outside the main island via the first return 
map $L$ induced by $F$ on $\Lambda$. This map is constructed as the 
composition of two transit maps
\begin{equation}\label{eq:Lbar}
L=\bar L^{\rm out}\circ \bar L^{\rm in}\hskip 30pt
\bar L^{\rm in}: \Lambda\to G(\Lambda)\hskip 30pt
\bar L^{\rm out}: G(\Lambda)\to \Lambda
\end{equation}
where $\bar L^{\rm in}$ and $\bar L^{\rm out}$ are the 
first-return maps induced by $F$ on the respective sets. 
The map $\bar L^{\rm in}$ is built up from 
iterations of $F^4$ {inside} $\Sigma$, the map $\bar L^{\rm out}$ 
is built up (after the initial step out of $\Sigma$) by the iteration 
of $F$ {outside} $\Sigma$.

The definition of the maps $\bar L^{\rm in,out}$ requires some care. 
As noted above, some orbits outside the main island do not enter $\Sigma$ 
at all, and it is conceivable that some orbits exiting $\Sigma$ 
will never return there. From Poincar\'e recurrence theorem,
we know that these points have zero measure; so we'll ignore their 
contribution and accept that the map $\bar L^{\rm out}$ may be undefined
on a zero measure subset of $G(\Lambda)$. In particular, the 
transit time of the map $\bar L^{\rm out}$ may be unbounded. The same is 
true for $\bar L^{\rm in}$, as a result of having dropped the small 
set $\varepsilon$ from the domain $\Lambda$.

The following result has important consequences for the dynamics of the sector $\Sigma$.

\begin{lemma}\label{lemma:Continuity} 
The restriction to the sector $\Sigma$ of the maps $F^{\pm 4}$ is continuous.
Continuity extends to the boundary of $\Sigma$, apart from a sliding 
singularity of $F^4$ on $\gamma$, and of $F^{-4}$ on $G(\gamma)$.
The symbolic dynamics is $(1,1,1,1)$ for points in $\Sigma\setminus G(\Lambda)$
and $(1,1,1,0)$ for points in $G(\Lambda)$.
\end{lemma}

\proof Since $F$ is continuous away from $\gamma$, and $\gamma={\bf 0}$ is 
a sliding singularity (the segment ${\bf 0}$ does not detach itself
from $\Sigma$ under iteration) we have 
\begin{equation} \label{eq:SigmaImages}
F^t(\Sigma)=\langle {\cal E},{\bf t+1},{\bf t}\,\rangle,\quad t=0,\ldots,3 
\hskip 40pt
F^4(\Sigma)=\langle {\cal E},{\bf 5}^1,{\bf 4}\,\rangle.
\end{equation}
From table \ref{table:gammaImages}, we find that, in the parameter 
range (\ref{eq:lambdaRange}), of the first four images of 
${\bf 0}$, only ${\bf 4}$ intersect $\gamma$. However, the point
at which ${\bf 4}$ is tangent to ${\cal E}$ lies to the left of
$\gamma$. Hence none of the first three images of $\Sigma$ intersects 
$\gamma$, and therefore $F^4$ is continuous on $\Sigma$. 
The continuity of the inverse follows from time-reversal symmetry.
From (\ref{eq:SigmaImages}), we find that $F(\Sigma)\cap \Sigma=\emptyset$,
and hence all points in $\Sigma \setminus G(\Lambda)$ return to $\Sigma$ 
in four iterations of $F$ and no fewer, with return symbolic dynamics $(1,1,1,1)$. 
Furthermore, from (\ref{eq:Hm}) and (\ref{eq:DomainBoundaries}), we
have $G(\Lambda)=\langle{\bf -4},{\bf 1},{\bf 0}\rangle$,
and since $G(\Lambda)\subset \Sigma$ and $F^3$ is continuous there, we have
$F^3(G(\Lambda))=\langle {\bf -1},{\bf 4},{\bf 3}\,\rangle\subset \Omega_0$.
So the symbolic dynamics of a point in $G(\Lambda)$ is (1,1,1,0).
\endproof

The transit maps $\bar L^{\rm in, out}$ are piecewise isometries, and our 
next step is to partition their respective domains $\Lambda$ and $G(\Lambda)$ 
into atoms. 
Since the two domains are mirror image of each other under the involution
$G$ (see equation (\ref{eq:InOut})), we replace (\ref{eq:Lbar}) by the
more convenient decomposition 
\begin{equation}\label{eq:L}
L=L^{\rm out}\circ L^{\rm in};\qquad
L^{\rm in, out}: \Lambda\to \Lambda,\qquad
L^{\rm in}= G\circ \bar L^{\rm in},\qquad 
L^{\rm out}= \bar L^{\rm out}\circ G
\end{equation}
where $\Lambda$ is the common domain of the two maps.
We will often represent the right triangle $\Lambda$ using scaled coordinates, 
with the shorter leg being four times the length of the longest one 
---see figure \ref{fig:Atoms}. 

The distinction between {\it regular\/} and {\it irregular\/} atoms will be 
important in what follows. All regular atoms share a simple geometrical structure 
(apart from two exceptions), determined by regular components with respect 
to a periodic island. Asymptotically, these atoms occupy the entire domain 
$\Lambda^{\rm in}$ of definition of the maps. 
Irregular atoms are more varied, but their detailed description is not
required in asymptotic calculations, because in the limit, their total 
area becomes negligible.

\subsection{Preliminary lemmas}\label{section:PreliminaryLemmas}

The proof of theorem A will be quite laborious. We begin with some
preliminary lemmas on involutivity and symmetry. We call an atom
{\bf maximal\/} if it is not properly contained in another atom.

\begin{lemma}\label{lemma:MaximalAtoms}
Let $L$ be a piecewise isometric involution. Then the image of
a maximal atom is an atom. In particular, a maximal atom 
containing a fixed point is invariant.
\end{lemma}

\proof Since $L$ is an involution, its atoms are the same as those 
of $L^{-1}$. If $\Omega$ is an atom of $L$, we claim that $L$ is 
continuous on $L(\Omega)$. Indeed if it were not, we could find 
arbitrarily close points $z,w\in L(\Omega)$ mapping far apart in 
$\Omega$, thereby contradicting the fact that $L$ is an isometry 
on $\Omega$. Thus $L(\Omega)\subset \Omega'$, where $\Omega'$ is
another atom. If this inclusion were proper, then $L(\Omega')$ 
would be a connected set properly containing $\Omega$, because 
$L$ is continuous on $\Omega'$ and volume-preserving. 
But since $L$ is continuous on $L(\Omega')$, then $\Omega$ would 
not be maximal. So, if $\Omega$ is maximal, then $L(\Omega)=\Omega'$, 
as desired. Furthermore, if $\Omega$ contains a fixed point of $L$, 
then $\Omega$ intersects $L(\Omega)$, and is therefore invariant.
\endproof

\begin{lemma} \label{lemma:Involutivity}
The maps $L^{\rm in, out}$ are involutions. 
\end{lemma}

\proof 
Let $z\in\Lambda$. Since $\bar L^{\rm in}$ is an induced map,
we have $\bar L^{\rm in}(z)=F^t(z)$, for some $t=t(z)$. From
(\ref{eq:Symmetry}), a straightforward induction gives 
$F^{-t}=G\circ F^t\circ G$, for all $t\in \Z$, 
and hence
\begin{eqnarray*}
L^{\rm in}\circ L^{\rm in}(z)&=&G\circ F^t\circ G\circ F^t(z)
=F^{-t}\circ F^t(z)=z.
\end{eqnarray*}
Since $z$ was arbitrary, $L^{\rm in}$ is an involution.
The involutive nature of $L^{\rm out}$ is proved with
an analogous argument.
\endproof

Next we look at the fixed sets of these involutions.

\begin{lemma} \label{lemma:FixedSets}
Let $\Lambda_k^{\rm in,out}$ be an atom of $L^{\rm in,out}$, and let
$t=t(k)$ be its transit time (number of iterations of $F$).
Then
$$
\vcenter{\baselineskip 20pt\halign{
 $\displaystyle #$\hfil& \hskip 1pt
 $\displaystyle #$\hfil&\qquad
 $\displaystyle #$\hfil \cr
{\rm Fix}\, L^{\rm in}_k&=F^{-s}({\rm Fix}\, G) \cap \Lambda_k^{\rm in}
    & \mbox{\rm if}\quad t=2s\cr
{\rm Fix}\, L^{\rm in}_k&=F^{-s}\circ G({\rm Fix}\, H)\cap \Lambda_k^{\rm in} 
    & \mbox{\rm if}\quad t=2s+1\cr
{\rm Fix}\, L^{\rm out}_k&=F^{s}({\rm Fix}\, G)\cap \Lambda_k^{\rm out} 
    & \mbox{\rm if}\quad t=2s\cr
{\rm Fix}\, L^{\rm out}_k&=F^{s}({\rm Fix}\, H)\cap \Lambda_k^{\rm out} 
    & \mbox{\rm if}\quad t=2s+1\cr
}}
$$
where $L_k^{\rm in, out}$ is the restriction of $L^{\rm in, out}$ 
to $\Lambda_k^{\rm in, out}$.
\end{lemma}

\proof
From (\ref{eq:L}), the fixed points equations become
\begin{equation}\label{eq:FixedPointsEquations}
L^{\rm in}(z)=z \iff G\circ F^{t}(z)=z\qquad
L^{\rm out}(z)=z \iff F^{t}\circ G(z)=z
\end{equation}
and we are interested in solutions (if any) within $\Lambda_k^{\rm in,out}$.
First, we consider the inner map $L^{\rm in}$, and $z\in\Lambda_k^{\rm in}$.
If $t=2s$, equation (\ref{eq:FixedPointsEquations}) becomes
\begin{eqnarray*}
G\circ F^{2s}(z)=z&\iff& F^{-2s}\circ G(z)=z\\
 &\iff& F^{-s}\circ G(z)=F^s(z) \\
 &\iff& G\circ F^s(z)=F^s(z).
\end{eqnarray*}
This shows that $F^s(z)\in {\rm Fix}\, G$, namely that $z\in F^{-s}({\rm Fix}\, G)$.
If $t=2s+1$ is odd, we obtain, using (\ref{eq:HG})
\begin{eqnarray*}
G\circ F^{2s+1}(z)=z&\iff& F^{-(s+1)}\circ G(z)=F^s(z)\\
 &\iff& F\circ F^s(z)=G\circ F^s(z)\\
 &\iff& H\circ G\circ F^s(z)= G\circ F^s(z).
\end{eqnarray*}
We find that $G\circ F^s(z)\in {\rm Fix}\, H$, or
$z\in F^{-s} \circ G({\rm Fix}\, H)$.

Similarly, for $z\in\Lambda_k^{\rm out}$, and $t$ even, the outer map gives
\begin{eqnarray*}
F^{2s}\circ G(z)=z&\iff& F^{s}\circ G(z)=F^{-s}(z)\\
 &\iff& G\circ F^{-s}(z)=F^{-s}(z),
\end{eqnarray*}
so $F^{-s}(z)\in {\rm Fix}\,G$, that is, $z\in F^{s}({\rm Fix}\,G)$.
Finally, if $t$ is odd
\begin{eqnarray*}
F^{2s+1}\circ G(z)=z &\iff& F^{s+1}\circ G(z)=F^{-s}(z)\\
 &\iff& G\circ F^{-1}\circ F^{-s}(z)= F^{-s}(z)\\
 &\iff& H\circ F^{-s}(z)=F^{-s}(z)
\end{eqnarray*}
which gives $F^{-s}(z)\in {\rm Fix}\,H$, that is, 
$z\in F^{s}({\rm Fix}\,H)$.
\endproof

We will show that for a significant class of atoms, the transit times
of $L^{\rm in}$ and $L^{\rm out}$ are even and odd, respectively, and 
the fixed lines of the involutions intersect these atoms, which are
themselves invariant. This will be done in the next two sections, where 
we look in detail at the transit maps.

\subsection{Structure of $L^{\rm in}$}\label{section:Lin}
 
The map $L^{\rm in}$ results from restricting to $\Sigma$
the action of the matrix $C^4$ (see equation (\ref{eq:Matrix})) with 
respect to the fixed point of $F$.
We consider again the point $T_0$ at which the ellipse ${\cal E}$ 
is tangent (on the torus) to the segment ${\bf 1}$, see figure \ref{fig:Sector}.
Let $T_1,T_2,\ldots$ be the images of $T_0$ under $F^4$, and let 
\begin{equation}\label{eq:N}
N=\left\lfloor \frac{1}{4(1-4\rho)}-\frac{1}{4}\right\rfloor
\hskip 40pt
\rho=\frac{1}{2\pi}\cos^{-1}(\lambda/2),
\end{equation}
where $\rho$ is the rotation number, see equation (\ref{eq:Matrix}).
Then the first $N$ points of this sequence lie on the ellipse 
$\cal E$, within the boundary of $\Sigma$. 
 
\begin{figure}[t]
\begin{center}
\epsfig{file=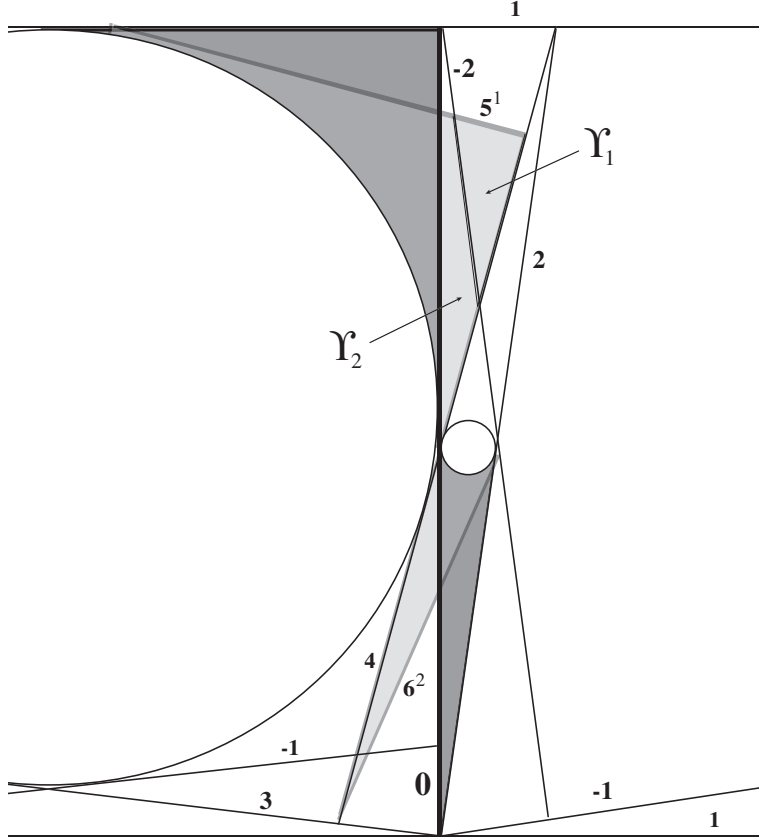,width=10cm}
\end{center}
\caption{\label{fig:TwoSectors}
\small The linked sector maps of the 1-cycle (the large circle) and the 2-cycle 
(the small circle), which lie on opposite sides of the discontinuity 
line ${\bf 0}$. Under the action of the map $F^4$, part of the sectors 
$\Sigma$ and $\Sigma'$ (the darker regions) escape to the triangles 
$\Xi$ and $\Xi'$, respectively, located on the opposite side of ${\bf 0}$. 
The points in $\Xi$ that do not end up in $\Sigma'$ (under $F^2$)
comprise the triangle $\Upsilon_1$; those that do comprise the
quadrilateral $\Upsilon_2$.
}
\end{figure}

\begin{figure}[t]
\begin{center}
\epsfig{file=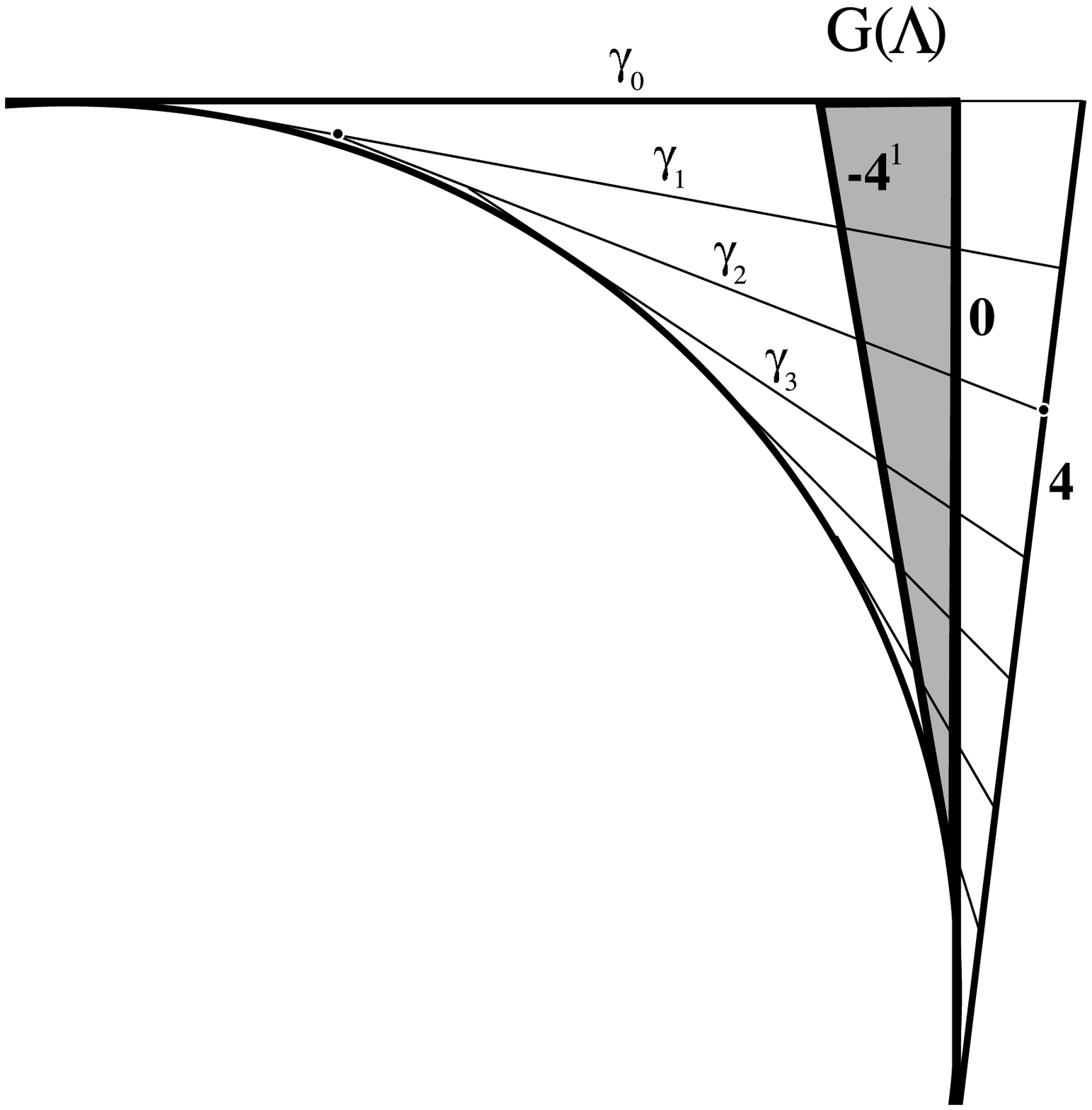,width=7cm}
\end{center}
\caption{\label{fig:RegularComponents}
\small The geometry of lemma \ref{lemma:RegularComponents}: the regular components 
$\gamma_n$ intersect transversally the boundaries ${\bf -4}^1$ and ${\bf 0}$ of 
the exit domain $G(\Lambda)$, creating pseudo-hyperbolic points.
}
\end{figure}

For $n=0,\ldots, N-1$, we let $\gamma_n=({\bf 4n+1})^1$ be the regular 
component of $F^{4n}({\bf 1})$ with respect to ${\cal E}$. Then $\gamma_n$ 
is tangent to ${\cal E}$ at the point $T_n$.
For any $k$, $\gamma_n$ and $\gamma_{n+k}$ are quasi-parallel,
and their orientations agree.
The segments $\gamma_n$ will be used to construct the atoms of $L^{\rm in}$.
We begin with a lemma, whose geometric content is illustrated in figure
\ref{fig:RegularComponents}.

\begin{lemma}\label{lemma:RegularComponents}
For $n=1,\ldots,N-1$, the regular component $\gamma_n$ connects ${\bf 4}$ 
to $\gamma_{n-1}$, and it intersects transversally the segments ${\bf 0}$ 
and ${\bf -4}^1$ on the boundary of $G(\Lambda)$, apart from a possible 
non-generic degeneracy of $\gamma_{N-1}\cap {\bf -4}^1$ into a segment.
\end{lemma}

\proof 
From (\ref{eq:DomainBoundaries}) we verify that ${\bf 0}$ and
${\bf -4}^1$ do indeed belong to the boundary of $G(\Lambda)$.
Let $n=1$. The segment $\gamma_1={\bf 5}^1$ connects {\bf 4} to 
$\gamma_0={\bf 1}$ --- see table \ref{table:gammaImages}. 
For the parameter $\lambda$ in the range (\ref{eq:lambdaRange}), 
the point $T_1$ lies outside $G(\Lambda)$, and therefore $\gamma_1$ crosses 
${\bf 0}$ and ${\bf -4}^1$ transversally.
We proceed by induction on $n$, and assume that for some $n$ in the range 
$1\leq n < N-1$, the segment $\gamma_n$ connects ${\bf 4}$ 
to $\gamma_{n-1}$, and it is tangent to ${\cal E}$ at a point $T_n$
lying outside $G(\Lambda)$. Then $\gamma_n$ intersects $\gamma$ transversally, 
so we consider the segment $\gamma_n^-=\gamma_n\cap\Sigma$.
By construction, the segment $\gamma_n^-$ connects $\gamma$ to $\gamma_{n-1}$.
From lemma \ref{lemma:Continuity}, the set $\gamma_{n+1}=F^4(\gamma_n^-)$
is a segment tangent to ${\cal E}$ at a point $T_{n+1}$, which connects 
${\bf 4}$ to $\gamma_n$ (the statement concerning the first end-point 
is justified by the sliding nature of the singularity on $\gamma$). 
Since $F^4$ is close to the identity, the orientations of $\gamma_n$
and $\gamma_{n+1}$ agree. Because $n<N-1$, the point $T_{n+1}$ 
lies outside $G(\Lambda)$, which completes the induction.
An exception may occur for $n=N-1$, in the non-generic
case in which $T_{N-1}$ lies on the boundary segment ${\bf -4}^2$ 
of $G(\Lambda)$. In this case $\gamma_{N-1}$ is tangent to that boundary.
\endproof

From the above lemma, we see that the segments $\gamma_n$ decompose
$G(\Lambda)$ into $N-1$ quadrilaterals, plus a small 
residual region near the tip. We are then led to consider the
half-open quadrilaterals $\Lambda_n^{\rm in}\subset \Lambda$, given by
\begin{eqnarray}
\Lambda_n^{\rm in}&=&
\langle\, \gamma_0,\, G(\gamma_{n-1}),\, \gamma_1, \, G(\gamma_{n})\,\rangle \nonumber\\
&=&
 \langle\,{\bf 1},({\bf -4n+4})^1,{\bf 5}^1,({\bf -4n})^1\,\rangle 
\qquad n=1,\ldots,N-1. \label{eq:AtomsLin}
\end{eqnarray}
The sides are oriented in such a way that
\begin{equation}\label{eq:AtomsBoundaries}
\gamma_0,\, G(\gamma_{n-1})\cap\Lambda_n^{\rm in}=\emptyset\qquad
\gamma_1, \, G(\gamma_{n}) \cap\Lambda_n^{\rm in}\not=\emptyset.
\end{equation}
We also represent polygons as lists of vertices
\begin{equation}\label{eq:AtomsVertices}
\Lambda_n^{\rm in}=[\,Q_0(n),Q_0(n-1),Q_1(n-1),Q_1(n)\,]\qquad
Q_{0,1}(n)=\gamma_{0,1}\cap G(\gamma_n).
\end{equation}
Below, we will show that the sets $\Lambda_n^{\rm in}$ are atoms of 
$L^{\rm in}$. Since these domains are bounded by regular components, 
we shall call them {\bf regular atoms};
all other atoms will be termed {\bf irregular}.

From (\ref{eq:AtomsVertices}), we see that the vertices of $\Lambda_n^{\rm in}$
are pseudo-hyperbolic points of the map $F$. Indeed, 
for $n=1,\ldots,N-1$, the point $Q_0(n)$ is the intersection of $\gamma_0=G(\gamma)$ 
and a segment in the discontinuity set, and hence it is the left 
end-point of a pseudo-hyperbolic sequence (see remarks at the end of
section \ref{section:Preliminaries}). 
The point $Q_1(n-1)=F^4(Q_0(n)$ belongs to the same orbit, and since it 
is the transversal intersection of two segments (by virtue of lemma 
\ref{lemma:RegularComponents}), it also is pseudo-hyperbolic.
Under the action of $F^4$, these points propagate across the sector 
$\Sigma$, forming a regular pattern of intersecting pseudo-separatrices 
---see figure \ref{fig:BasicDiscontinuitySet}. The point $G(Q_0(n))$,
which belongs to $\gamma$, is the right endpoint of the pseudo-hyperbolic
sequence that starts at $Q_0(n)$.

\begin{theorem}\label{theorem:Lin} The regular atoms of the map $L^{\rm in}$
are the half-open quadrilaterals (\ref{eq:AtomsLin}), with the boundary
specified in (\ref{eq:AtomsBoundaries}).
For each $n$, the atom $\Lambda_n^{\rm in}$ is invariant under $L^{\rm in}$, 
with return symbolic dynamics $(1^{4(n-1)})$ (the exponent denotes repetition 
of symbols). The fixed set of $L_n^{\rm in}$ is the segment connecting 
the vertex $Q_1(n)$ to the opposite vertex $Q_0(n-1)$, where $Q_{0,1}$
is defined in (\ref{eq:AtomsVertices}).
There are at most five other atoms, of total area $O(\lambda^3)$.
\end{theorem}

\proof Consider the quadrilateral 
$$
G(\Lambda_n^{\rm in})=
\langle \gamma_{n-1},\,G(\gamma_0),\,\gamma_n,\,G(\gamma_1)\rangle.
$$
From the proof of lemma \ref{lemma:RegularComponents}, we know 
that $\gamma_n$ and $G(\gamma_1)$ are glued to 
$G(\Lambda_n^{\rm in})$, but the other two sides are not.
Applying $F^{-4(n-1)}$ to each segment in $G(\Lambda_n^{\rm in})$ we obtain 
the quadrilateral 
$\langle\,\gamma_0,\, G(\gamma_{n-1}),\, \gamma_1, \, G(\gamma_{n})\,\rangle$, 
and since the iterates of these segments remain within $\Sigma$, we
deduce from (\ref{eq:AtomsLin}) that $\Lambda_{n}^{\rm in}$ 
contains an atom of $L^{\rm in}$, 
glued to $\gamma_1$ and $G(\gamma_n)$, with transit time $4(n-1)$.

It remains to show that $\Lambda_n^{\rm in}$ is actually an atom.
The case $n=1$ is trivial, since $\Lambda_1^{\rm in}=\Lambda\cap G(\Lambda)$,
with zero transit time. 
For $n>1$, we must verify that, in the recursive construction
of $\Lambda_{n}^{\rm in}$, no forward image of the segment 
$\gamma_{n-1}^+=\gamma_{n-1}\cap \Xi$ (the piece of $\gamma_{n-1}$ that gets cut off
by the discontinuity) reaches $G(\Lambda_n^{\rm in})$ before $\gamma_n$,
namely that
\begin{equation}\label{eq:NonInterference}
\bigcup_{t=1}^4F^t(\gamma_{n-1}^+)\cap G(\Lambda_n^{\rm in})=\emptyset
\qquad n>1.
\end{equation}
Let $z\in\gamma_{n-1}^+$. From lemma \ref{lemma:RegularComponents} 
and table \ref{table:gammaImages}, we have $z\in G(\Omega_0)$. 
We have two cases.
If $F(z)\in\Omega_0$, then $F^2(z)\in G(\Omega_0)$, and $F^3(z)$ is either
again in $\Omega_0$, or in $\langle\,{\bf -1},{\bf 2},{\bf 3}\,\rangle$;
in either case $F^4(z)$ is outside $\Sigma$, hence outside $G(\Lambda_n^{\rm in})$. 
If $F(z)\in\Omega_1$, then $F(z)\in\langle\,{\bf -1},{\bf 2},{\bf 3}\,\rangle$,
hence $F^2(z)\in\langle\,{\bf 0},{\bf 3},{\bf 4}\,\rangle$. The points of the
latter triangle map either to $G(\Omega_0)$ or to $\Lambda$. From neither set
it is possible to reach $G(\Lambda_n^{\rm in})$ in one iteration.
We have established equation (\ref{eq:NonInterference}), and 
$\Lambda_n^{\rm in}$ is an atom of $L^{\rm in}$.

Generically, the point $T_{N}$ lies between $G(T_0)$ and $G(T_1)$,
and the residual region consists of one quadrilateral, two triangles and 
two triangular sectors. At isolated parameter values (determined by the 
condition that $1/4(1-4\rho)$ be an integer), we have $T_{N}=G(T_0)$, 
and the residual region reduces to one triangular sector,
while the regular atom $\Lambda_{N}$ degenerates to a triangle.
By construction, the total area of the residual region is $O(\lambda^3)$
in both cases.

From (\ref{eq:AtomsBoundaries}), we see that $\gamma_1$ and $G(\gamma_n)$ are 
part of the included boundary of $\Lambda_n^{\rm in}$, and since
$\gamma_n=F^{4(n-1)}(\gamma_1)$, the transit time of this atom is equal 
to $4(n-1)$, which is even.
With the notation of lemma \ref{lemma:FixedSets}, we have $s=2(n-1)$, 
and hence ${\rm Fix}\, L^{\rm in}_n=F^{2(1-n)}({\rm Fix}\, G)\cap \Lambda_n^{\rm in}$. 
Now, by symmetry, 
\begin{equation}\label{eq:PointsInFixG}
\gamma_n\cap G(\gamma_n)\,\in\, {\rm Fix}\, G\cap \Sigma \qquad 
   n=0,\ldots,\lfloor (N-1)/2\rfloor.
\end{equation}
We also have the symmetric points
\begin{equation}\label{eq:ExtraPointsInFixG}
F^{4n-2}(\gamma_0)\cap G\circ F^{4n-2}(\gamma_0)\,\in\,{\rm Fix}\,G \qquad
   n=1,\ldots,\lfloor (N-1)/2\rfloor
\end{equation}
which lie in the South-West sector. (With a slight abuse of notation,
we have denoted by $F^{4n-2}(\gamma_0)$ the regular component of that set,
with respect to ${\cal E}$.)
Applying the map $F^{-2}$ to the points (\ref{eq:ExtraPointsInFixG}), 
we obtain a second sequence of points in $\Sigma$
$$
\gamma_{n-1}\cap G(\gamma_{n}) \in F^{-2}({\rm Fix}\, G) \cap \Sigma\qquad
   n=1,\ldots,\lfloor (N-1)/2\rfloor.
$$
Then, by applying repeatedly the identity 
$$
F^{-4}(\gamma_j \cap G(\gamma_k))=\gamma_{j-1}\cap G(\gamma_{k+1})
$$
to pairs of adjacent points in the above two sequences, we 
translate images of ${\rm Fix}\, G$ within $\Lambda$, and we
see that ${\rm Fix}\, L^{\rm in}_n$ is the segment connecting
$\gamma_0\cap G(\gamma_{n-1})$ to $\gamma_1\cap G(\gamma_{n})$.
Comparison with (\ref{eq:AtomsLin}) shows that this is the segment 
connecting the North-East and South-West corners of the $n$th atom
(with the obvious modification if the last regular atom
degenerates into a triangle).

Now, $L^{\rm in}$ is an involution (lemma \ref{lemma:Involutivity}),
and its regular atom $\Lambda_n^{\rm in}$ is maximal, being bound by 
images of the discontinuity line. Then every regular atom is invariant, 
because it contains an invariant segment (lemma \ref{lemma:MaximalAtoms}).
\endproof

\subsection{The sector map of the 2-cycle}

The dynamics of the map $L^{\rm out}$ is dominated by the
2-cycle of $F$ given in table (\ref{eq:BasicCycles}), which 
plays a role analogous to the 1-cycle for the map $L^{\rm in}$.
The 2-cycle generates the regular atoms of the map.
There are also irregular atoms, but they have negligible measure.

We noted that the 2-cycle is symmetrical, and it belongs to ${\rm Fix}\, H$ 
(see beginning of section \ref{section:ReturnMap}). We consider the element 
of the cycle that lies within the atom $\Omega_1$ of $F$, denoting by 
${\cal E}'$ the ellipse bounding its cell (figure \ref{fig:TwoSectors}). 
From the data (\ref{eq:BasicCycles}), we find that the
segment $\gamma={\bf 0}$ is tangent to ${\cal E}'$ at the point $(0,1/2)$.
By $H$-symmetry, the segment $H(\gamma)=H({\bf 0})={\bf 2}$ ---see (\ref{eq:Hm})---
is also tangent to ${\cal E}'$ at the point $T_0'=H(0,1/2)=(\lambda/2,1/2)$.
The open region delimited by the segments ${\bf 0}$ and ${\bf 2}$,
and the elliptical arc connecting $T_0'$ to $H(T_0')$ clockwise
will be called the {\bf 2-sector\/} $\Sigma'$ 
\begin{equation}\label{eq:SigmaII}
\Sigma'=\langle\,{\cal E}',{\bf 2},{\bf 0}\,\rangle.
\end{equation}
This region is $H$-invariant, and both ${\bf 0}$ and ${\bf 2}$ are are glued to it.

The analysis of the dynamics of the 2-sector now proceed as for the 1-sector;
we'll use the same notation, with primed symbols.
First we define entry and exit domains
---the analogue of (\ref{eq:InOut})
$$
\Lambda'=\Sigma'\setminus F^4(\Sigma')
\hskip 40pt
\Xi'=F^4(\Sigma')\setminus \Sigma' =F^4\circ H(\Lambda').
$$
Let 
\begin{equation}\label{eq:M}
M=\left\lfloor \frac{1}{2(1-4\rho)}-\frac{1}{2}\right\rfloor\,\sim\,
\left\lfloor \frac{\pi}{2\lambda}\right\rfloor\qquad \mbox{\rm as}\quad \lambda\to 0
\end{equation}
where $\rho$ is the rotation number; this is the analogue of equation (\ref{eq:N}).
For $m=0,\ldots, M-1$, we let $\gamma_m'=({\bf 4m+2})^2$ be the regular component 
of $F^{4m}({\bf 2})$ with respect to ${\cal E}'$,
denoting by $T_m'$ the corresponding point of tangency. 
The segments $\gamma_m'$ are pairwise quasi-parallel, and their 
orientations agree.
The following result is the twin of lemma \ref{lemma:Continuity}.

\begin{lemma}\label{lemma:ContinuityII}
The restriction to the closure of the sector $\Sigma'$ of the maps $F^{\pm 4}$ 
is continuous. The symbolic dynamics is (1,0,1,0) for points in 
$\Sigma'\setminus H(\Lambda')$ and (1,0,1,1) for points in $H(\Lambda')$.
\end{lemma}

\proof 
From (\ref{eq:SigmaII}) we find that
$F(\Sigma')=\langle\,{\cal E}',{\bf 3},{\bf 1}\,\rangle$.
Considering that $ F(T_0')=((1+\lambda^2)/2,\, \lambda/2)\in\Omega_0 $,
from table \ref{table:gammaImages} we deduce that $F(\Sigma')\subset\Omega_0$. 
Then we have $F^2(\Sigma')\subset G(\Omega_0)$, hence $F^3(\Sigma')\cap \gamma=\emptyset$,
and so $F^4$ is continuous on $\Sigma'$. By time-reversal symmetry, the inverse
$F^{-4}$ is also continuous on $\Sigma'$.
Since the boundaries of $\Sigma'$ are glued to it, and remain glued under iteration,
continuity extends to the boundary.

Every point in $\Sigma'\setminus H(\Lambda')$ follows the 2-cycle,
and so the return dynamics is (1,0,1,0). We determine the itinerary of
$H(\Lambda')$. 
We have $\Lambda'=\langle{\bf 0},{\bf 6}^2,{\bf 2}\rangle$, and hence
$H(\Lambda')=\langle{\bf 2},{\bf 0},{\bf -4}^2\rangle\in\Omega_1$.
This gives $F\circ H(\Lambda')=\langle{\bf 3},{\bf 1},{\bf -3}\rangle \subset \Omega_0$,
whence $F^2\circ H(\Lambda')\subset \Omega_1$, and finally 
$F^3\circ H(\Lambda')=\langle{\bf 5}^2,{\bf 3},{\bf -1}\rangle \subset \Omega_1$.
All inclusion relations may be verified with the aid of table 
\ref{table:gammaImages}, together with (\ref{eq:Gm}) and (\ref{eq:Hm}).
\endproof

The continuity of $F^4$ on $\Sigma'$ implies the existence of an 
orderly array of regular components.

\begin{lemma}\label{lemma:RegularComponentsII}
For $m=1,\ldots,M-1$, the regular component $\gamma_m'$ intersects 
$\gamma$ transversally, apart from a possible non-generic degeneracy of 
$\gamma_{M-1}'\cap H(\gamma_1')$ into a segment.
For $m>1$, $\gamma_m'$ connects $\gamma_{m-1}'$ to ${\bf 4}$.
\end{lemma}

\proof
The proof of this statement is analogous to that of lemma 
\ref{lemma:RegularComponents}, with the added simplification that
continuity now extends to the boundary. The same applies to the argument
concerning the degeneracy for $\gamma_{M-1}'$.
We omit the details for the sake of brevity. 
\endproof

The following result is the analogue of theorem \ref{theorem:Lin}
for the 2-sector $\Sigma'$.

\begin{lemma}\label{lemma:Lin2}
For $m=1,\ldots,M-1$, the regular atoms $\Lambda_m'$ of the transit 
map $\Lambda'\to H(\Lambda')$ induced by $F$ are the half-open 
quadrilaterals
\begin{equation}\label{eq:Lin2}
\Lambda_m'=\langle\gamma_0',\,H(\gamma_{m-1}'),\,\gamma_1',\,H(\gamma_{m}')\rangle
=
  \langle\, {\bf 2},({\bf -4m+4})^2,{\bf 6}^2,({\bf -4m})^2 \rangle.
\end{equation}
The boundary sides  $\gamma_0'$ and $H(\gamma_m')$ belong to $\Lambda_m'$, the other don't.
The return symbolic dynamics of the $m$-th atom is $((1,0)^{2m-2})$.
There are at most five other atoms, of total area $O(\lambda^3)$.
\end{lemma}

\proof The proof is analogous to that of theorem \ref{theorem:Lin}, with $H$
in place of $G$, and primed variables.
Identity (\ref{eq:NonInterference}) now reads
$$
\bigcup_{t=1}^4F^t(\gamma_{m-1} \cap Xi')\cap \Lambda_m'=\emptyset \qquad m>1.
$$
To establish this identity, we note that $\gamma_{m-1}'$ intersects $\gamma$.
The right component $\gamma_{m-1}'\cap \Sigma'$ will iterate to produce 
$\gamma_m'$, and we must ensure that the left component 
$\gamma_{m-1}'\cap \Xi'$ does not interfere with this process. 
We find 
$$
\Xi_1'\stackrel{\rm def}{=}\Xi'\cap \Omega_1\subset F^{-1}(\Lambda)=
\langle {\bf -1},{\bf 4},{\bf 0}\rangle
$$
and the map $F^4$ is continuous over the closure of that domain, from 
lemma \ref{lemma:Continuity}, apart from a sliding singularity on 
${\bf -1}$ and ${\bf 0}$.
Thus 
$$
\vcenter{\halign{$\displaystyle #$\hfil&\qquad$\displaystyle #$\hfil\cr
F^1(\Xi_1')\subset \langle{\bf 0},{\bf 5}^1,{\bf 1}\rangle   &
F^2(\Xi_1')\subset \langle{\bf 1},{\bf 6}^1,{\bf 4}\rangle \cr
F^3(\Xi_1')\subset \langle{\bf 2},{\bf 7}^1,{\bf 3}\rangle &
F^4(\Xi_1')\subset \langle{\bf 3},{\bf 8}^1,{\bf 4}\rangle \cr
}}
$$
all of which lie outside $\Xi'$, from table \ref{table:gammaImages}.
Similarly, we find 
$$
\Xi_0'\stackrel{\rm def}{=}\Xi'\cap \Omega_0\subset \langle {\bf -1},{\bf 0},{\bf 3} \rangle
$$
and one verifies that $F^3$ is continuous on the closure of that set. 
Hence 
$$
F^1(\Xi_0)\subset\langle {\bf 0},{\bf 1},{\bf 4}\rangle\qquad
F^2(\Xi_0)\subset\langle {\bf 1},{\bf 2},{\bf 5}^2\rangle\qquad
F^3(\Xi_0)\subset\langle {\bf 2},{\bf 3},{\bf 6}^2\rangle.
$$
None of these sets intersects $\Xi'$; the last one is adjacent to $\Xi'$ an
intersects $\gamma$. Thus $F^4(\Xi_0)$ is contained in the union of two triangles, 
which clearly lie outside $\Xi'$. 
\endproof

As for the sector $\Sigma$, the vertices of the atoms $\Lambda_m'$ of the 
transit map over $\Sigma'$ are pseudo-hyperbolic points of the map $F$. 
Indeed such vertices result from transversal intersections of regular 
components ---see equation (\ref{eq:Lin2}) and figure \ref{fig:TwoSectors}. 
Two vertices of the atom $\Lambda_1'$ belong to ${\bf 0}=\gamma$, and 
hence they are right end-points of pseudo-hyperbolic sequences.

\subsection{Structure of $L^{\rm out}$}

To construct the atoms of the map $L^{\rm out}$, we need to connect dynamically 
the exit domain of the sector $\Sigma$ to the entry domain of the sector 
$\Sigma'$, and vice-versa.
These connections will not match perfectly, but the mismatch will involve sets 
of negligible measure. In addition, there will be one atom of $L^{\rm out}$ 
by-passing $\Sigma'$ altogether (see figure \ref{fig:TwoSectors}).

All points leaving $\Sigma$ via $\Lambda^{\rm out}$ end up in $\Xi$, and since 
we want to identify the points that will enter $\Sigma'$, we examine 
inverse images of the entry domain $\Lambda'$. The second pre-image of $\Lambda'$
intersects $\Xi$, and so in order to connect the 1-sector to the 2-sector, 
we consider the intersection and symmetric difference 
of the sets $\Xi$ and $F^{-2}(\Lambda')$.
Likewise, connecting back the 2-sector to the 1-sector will involve the sets
$\Xi'$ and $F^{-1}(\Lambda)$.

Our first lemma quantifies these intersections and symmetric differences. 
Its proof will also provide information about the atoms of $L^{\rm out}$.

\begin{lemma}\label{lemma:Mismatch}
The points of $\Lambda^{\rm out}$ return to $\Sigma$ either via the
domain $\Upsilon_1=\langle {\bf -2},{\bf 5}^1,{\bf 4}\rangle$,
or via $\Lambda'$, apart from a set of 
points of measure $O(\lambda^3)$.
\end{lemma}

\proof
First, we connect $\Xi$ to $\Lambda'$.
By the continuity of $F^{-2}$ on $\Sigma'$
(lemma \ref{lemma:ContinuityII}), we have
$ F^{-2}(\Lambda')=\langle {\bf -2},{\bf 0},{\bf 4}\rangle.$
Defining
\begin{equation}\label{eq:Upsilon12}
\Upsilon_1=\langle {\bf -2},{\bf 5}^1,{\bf 4}\rangle
\hskip 40 pt
\Upsilon_2=\langle {\bf -2},{\bf 4},{\bf 0},{\bf 5}^1 \rangle
\end{equation}
one verifies that
\begin{equation}\label{eq:IntersectionSymmetricDifference}
\Xi\,\,\Delta\,\,F^{-2}(\Lambda')\,=\,\Upsilon_1\cup\langle{\bf -2},{\bf 5}^1,{\bf 0}\rangle
\hskip 40pt
\Xi\,\cap\,F^{-2}(\Lambda')\,=\,\Upsilon_2.
\end{equation}
From the orthogonality relations
${\bf -2}\,\parallel\, {\bf 0}\, \perp\, {\bf 5}^1\,\parallel\, {\bf 1}$,
and the fact that the length of the side ${\bf 0}$ of 
$\langle{\bf -2},{\bf 5}^1,{\bf 0}\rangle$ is $O(\lambda)$, we conclude 
that $|\langle{\bf -2},{\bf 5}^1,{\bf 0}\rangle|=O(\lambda^3)$, 
where $|\cdot|$ denotes the two-dimensional area determined by the 
metric ${\cal Q}$. Thus
\begin{equation}\label{eq:EstimateI}
|\Xi\,\,\Delta\,\,F^{-2}(\Lambda')|\,=\,|\Upsilon_1| + O(\lambda^3).
\end{equation}

\begin{figure}[t]
\hfil
\epsfig{file=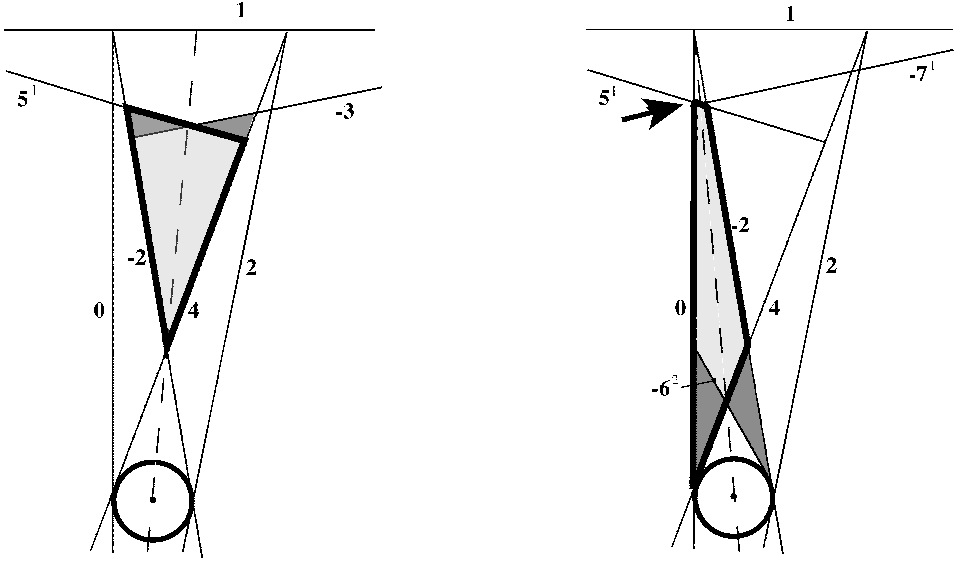,height=8cm}
\hfil
\caption{\label{fig:Turnstiles}
\small Left: Symmetrisation of the set $\Upsilon_1$ (the thick-sided triangle),
with respect to the symmetry line ${\rm Fix}\,H$ (the dashed line).
The two darker triangles represent the turnstile.
Right: Symmetrisation of the quadrilateral $\Upsilon_2$, with the appearance 
of two turnstiles; the arrow denotes the location of the smaller turnstile.
The dashed line is ${\rm Fix}\,F^{-4}\circ H$.
The symmetrised polygons (light grey) are isometric images of
the first and second regular atom of $L^{\rm out}$, respectively.
}
\end{figure}

We must now consider the dynamics of $\Upsilon_1$ and of $\Upsilon_2$, 
from $\Xi$ to $F^{-1}(\Lambda)$. 
We begin with the former. From (\ref{eq:DomainBoundaries}) and 
(\ref{eq:IntersectionSymmetricDifference}), we find, using continuity
$$
F^2(\Upsilon_1)=\langle {\bf 0},{\bf 7}^{+},{\bf 6}^2\rangle
\hskip 40pt
F^{-1}(\Lambda)=\langle {\bf -1},{\bf 0},{\bf 4}\rangle
$$
where ${\bf 7}^{+}=F^2({\bf 5}^1\cap \Xi)$.
Hence $F^2(\Upsilon_1)\cap F^{-1}(\Lambda)\not=\emptyset$, and so 
some points in $\Upsilon_1$ reach $\Lambda$ in 3 iterations.
However, not all of them do. To see this, we note that
${\bf 5}^1$ intersects ${\rm Fix}\, H$ transversally, and this intersection
belongs to the boundary of $\Upsilon_1$. Then, by symmetry,
\begin{equation}\label{eq:Turnstile}
\Upsilon_1 \,\Delta\, H(\Upsilon_1)=\Upsilon_1'\cup H(\Upsilon_1')
\qquad\mbox{\rm where}\qquad
\Upsilon_1'= \langle {\bf -2},{\bf 5}^1,{\bf -3}\rangle.
\end{equation}
The $H$-symmetric set (\ref{eq:Turnstile}) will be called a 
{\bf turnstile} (figure \ref{fig:Turnstiles}); by construction, 
$|\Upsilon_1'|=|H(\Upsilon_1')|$. 
Then the set $F^2(\Upsilon_1'\cup H(\Upsilon_1'))$ is
also a turnstile, which lies on the boundary of $F^{-1}(\Lambda)$. Finally,
$$
\Upsilon_1' \subset \Upsilon_1
\hskip 40pt
F^2(\Upsilon_1')\cap F^{-1}(\Lambda)=\emptyset
$$
which identifies $\Upsilon_1'$ as the set of points of $\Lambda$ that map 
into $\Upsilon_1$, but do not return to $\Lambda$ with it.
From the relations ${\bf -3}\parallel{\bf 5}^1$ and the fact that
the length of $\overline{\Upsilon}_1\cap {\bf 5}^1$ is $O(\lambda)$,
we have $|\Upsilon_1'|=O(\lambda^3)$. Thus
\begin{equation}\label{eq:EstimateII}
|F^3(\Upsilon_1)\cap \Lambda)|\,=\,|\Upsilon_1| + O(\lambda^3).
\end{equation}

There is an analogous phenomenon in the dynamics of $\Upsilon_2$,
defined in equation (\ref{eq:Upsilon12}). 
Since $H({\bf 0})={\bf 2}$, we have $F^{-4}\circ H({\bf 0})={\bf -2}$,
and $F^{-4}\circ H$ is an involution.
The boundary of $\Upsilon_2$ intersects ${\rm Fix}\,(F^{-4}\circ H)$ in two 
points, generating two turnstiles (figure \ref{fig:Turnstiles}, right). 
We are interested in the intersection 
${\bf 5}^1\cap {\rm Fix}\,(F^{-4}\circ H)$, which corresponds to the turnstile
\begin{equation}\label{eq:Turnstile2}
\Upsilon_2'\cup F^{-4}\circ H(\Upsilon_2')
\qquad\mbox{\rm where}\qquad
\Upsilon_2'= \langle {\bf 0},{\bf 5}^1,{\bf -7}\rangle.
\end{equation}
Then the set $F^6(\Upsilon_2'\cup F^{-4}\circ H(\Upsilon_2'))$ is
also a turnstile, which lies on the boundary of $F^{-1}(\Lambda)$. 
Finally,
$$
\Upsilon_2' \subset \Upsilon_2
\hskip 40pt
F^6(\Upsilon_2')\cap F^{-1}(\Lambda)=\emptyset
$$
which identifies $\Upsilon_2'$ as the set of points of $\Lambda$ that map 
into $\Upsilon_2$, hence to $\Lambda'$, but do not return to $\Lambda$ with it.
From the relation ${\bf 5}^1\parallel {\bf -7}$, and the fact that the 
length of $\overline{\Upsilon}_2\cap {\bf 5}^1$ is $O(\lambda^2)$,
we have $|\Upsilon_2'|=O(\lambda^4)$. This estimate, together with
equations (\ref{eq:IntersectionSymmetricDifference}), and the estimates 
(\ref{eq:EstimateI}) and (\ref{eq:EstimateII}), give the desired result.
\endproof

We now identify the regular atoms of $L^{\rm out}$. 
Let $\Upsilon_{1,2}$ be as in (\ref{eq:Upsilon12}), and let
$\Lambda_m'$ be as in lemma \ref{lemma:Lin2}.
We define the sets
\begin{eqnarray}
\Lambda_1^{\rm out}&=& G\circ F^{-4}(\Upsilon_1\cap H(\Upsilon_1))
\,=\, \langle {\bf 0},{\bf 8},{\bf 7},{\bf 1}\rangle \cap \Lambda
  \nonumber\\
\Lambda_2^{\rm out}&=& G\circ F^{-4}(\Upsilon_2\cap F^4\circ H(\Upsilon_2))
   \,=\, \langle {\bf 0},{\bf 12},{\bf 5},{\bf 11},{\bf 1},{\bf 7}\rangle \cap \Lambda
\label{eq:AtomsSidesII}\\
\Lambda_m^{\rm out}&=& G\circ F^{-6}(\Lambda_{m-1}')
    \,=\, \langle {\bf 1},{\bf 4m-1},{\bf 5},{\bf 4m+3}\rangle \cap \Lambda
\qquad 3\leq m \leq M.\nonumber
\end{eqnarray}
The rightmost equalities mean that the polygon on the left can be assembled with 
an appropriate choice of the components listed on the right, in the given order. 
The issue of uniqueness of a polygon with such a labelling is unimportant, 
as these polygons are continuous images of well-defined objects.
The boundaries of the sets $\Lambda_m^{\rm out}$ are specified by the
following non-empty intersections
\begin{equation} \label{eq:AtomsBoundariesII}
{\bf 5},{\bf 4m-1} \cap \Lambda_m^{\rm out}\not=\emptyset \qquad m>1
\end{equation}
and no other. In particular, $\Lambda_1^{\rm out}$ is open.
The sets (\ref{eq:AtomsSidesII}) are indeed quadrilaterals, 
apart from $\Lambda_2^{\rm out}$, which is a hexagon.

All atoms are tangent to ${\bf 1}$, and we let $P_0(m)$ be the left 
end-point of the intersection of ${\Lambda}_m^{\rm out}$ with ${\bf 1}$.
Then we let $P_1(m-1)$ be the vertex of $\Lambda_{m}^{\rm out}$ which is 
opposite to $P_0(m)$ (the term `opposite' refers to the fact that all
atoms have an even number of sides). We find that $P_1(0)$ and $P_1(1)$ 
are, respectively, the lower end-point of the intersection of $\Lambda_1^{\rm out}$ 
and $\Lambda_2^{\rm out}$ with ${\bf 0}$; likewise, for $m\geq 2$, the point 
$P_1(m)$ is the left end-point of the intersection of ${\Lambda}_m^{\rm out}$ 
with ${\bf 5}^1$. With these prescription we can characterize the 
sets $\Lambda_m^{\rm out}$ in terms of the vertices $P_{0,1}$, as follows
\begin{equation}\label{eq:AtomsVerticesII}
\Lambda_m^{\rm out}= [\,P_0(m),P_0(m-1),P_1(m-1),P_1(m)\,]
\qquad m >2.
\end{equation}
A complete description for the cases $m=1,2$ it not needed at this stage,
and will be given in section \ref{subsection:atomvertices}.

We can now state the main result of this section.

\begin{theorem}\label{theorem:Lout} 
The regular atoms $\Lambda_m^{\rm out}$ of the map $L^{\rm out}$ are the 
polygons (\ref{eq:AtomsSidesII}), with the boundary
specified in (\ref{eq:AtomsBoundariesII}).
For each $m$, the atom $\Lambda_m^{\rm out}$ is invariant under $L^{\rm out}$, 
and its return symbolic dynamics is
\begin{equation}\label{eq:SymbolicDynamicsOut}
(1,1,1,(0,1)^{2m-1},1,1).
\end{equation}
The fixed set of $L_m^{\rm out}$ is a segment connecting $P_0(m)$ to the 
opposite vertex $P_1(m-1)$. 
In addition, there is a (possibly infinite) set of irregular atoms 
of total area $O(\lambda^3)$.
\end{theorem}

\proof
From lemma \ref{lemma:Mismatch}, it suffices to consider the pre-images of
$\Upsilon_1$ and $\Lambda'$, and from the rightmost equation in 
(\ref{eq:IntersectionSymmetricDifference}), we may replace $\Lambda'$ by 
$\Upsilon_2$.

From the proof of lemma \ref{lemma:Mismatch}, we see that the part
of $\Upsilon_1$ that maps to $\Lambda$ under $F^3$ is the 
$H$-symmetrical quadrilateral
$$
\Upsilon_1\cap H(\Upsilon_1)\,=\, 
      \langle\, {\bf -3},{\bf 5}^1,{\bf 4},{\bf -2}\,\rangle\, \subset\, \Xi.
$$
To obtain the first atom $L_1^{\rm out}$, we apply $G\circ F^{-4}$ to it,
which maps it back to $\Lambda$. We obtain
\begin{equation}\label{eq:FirstAtom}
\Lambda_1^{\rm out}\in 
G\left(\langle\, {\bf -7},{\bf 1},{\bf 0},{\bf -6}\,\rangle \cap G(\Lambda)\right)=
   \langle\, {\bf 8},{\bf 7},{\bf 1},{\bf 0}\,\rangle \cap \Lambda,
\end{equation}
as desired. 

We have $\Upsilon_2\subset \Xi$, and so $G\circ F^{-4}(\Upsilon_2)\subset \Lambda$.
From the proof of lemma \ref{lemma:Mismatch}, we see that the part of $\Upsilon_2$ 
that maps to $\Lambda$ under $F^7$ is the $F^{-4}\circ H$-symmetrical hexagon
\begin{equation}\label{eq:Atom2}
\Upsilon_2\cap F^{-4}\circ H(\Upsilon_2)
= \langle{\bf 0},{\bf -7}^1,{\bf 5}^1,{\bf -2},{\bf 4},{\bf -6}^2\rangle
\subset \langle{\bf 0},{\bf -2},{\bf 4},{\bf -6}^2\rangle=F^{-2}(\Lambda_1')
\end{equation}
where $\Lambda_1'$ is the first atom of the 2-sector map (lemma \ref{lemma:Lin2}).
The atom $\Lambda_2^{\rm out}$ of $L^{\rm out}$ is the image under 
$G\circ F^{-4}$ of the symmetrised hexagon in (\ref{eq:Atom2}). We obtain 
\begin{equation}\label{eq:SecondAtom}
\Lambda_2^{\rm out}\,=\,
    \langle\, {\bf 0},{\bf 12},{\bf 5},{\bf 11},{\bf 1},{\bf 7}\rangle\cap \Lambda.
\end{equation}

The rest of $\Upsilon_2$, namely the triangle 
$$
F^{-2}(\Lambda'\setminus \Lambda_1') =\langle{\bf 0},{\bf -6}^2,{\bf 4}\rangle
\subset \Xi
$$
is covered, apart from a set of measure $O(\lambda^3)$, by the second pre-image
of the atoms $\Lambda_m', m=2,\ldots,M-1$ of the 2-sector map, from lemma 
\ref{lemma:Lin2}. The above inclusion, together with
lemmas \ref{lemma:Continuity} and \ref{lemma:ContinuityII}, show that 
$F^{-6}$ is continuous on $\Lambda'\setminus\Lambda_1'$.
Then, from lemma \ref{lemma:Lin2}, we have
$$
G\circ F^{-6}(\Lambda_{m-1}')\,=\, 
    \langle\, {\bf 1},{\bf 4m-1},{\bf 5},{\bf 4m+3} \rangle \cap \Lambda
$$
as desired

Verifying that the polygons defined in (\ref{eq:AtomsSidesII})
have the boundaries as given in (\ref{eq:AtomsBoundaries}) is a simple
exercise in orientation, recalling the conventions of table 1.

The code of each point in $\Lambda$ begins with the string $(1,1,1,0)$,
associated with the map $\Lambda^{\rm out}\to\Xi$ (lemma \ref{lemma:Continuity}).
One verifies directly that the points in $\Lambda_1^{\rm out}$ return to 
$\Lambda$ with the code $(1,1,1)$. For $m>1$, the $m$-th atom $\Lambda_m^{\rm out}$
follows the itinerary:
\begin{quote}
2 iterations to reach $\Lambda'$, with code $(1,0)$
(lemma \ref{lemma:ContinuityII});\newline
$4(m-1)$ iterations to reach $H(\Lambda')$, with code $(1,0)^{2m-2}$ 
(lemma \ref{lemma:Lin2});\newline
4 iterations to reach $\Xi'$, with code $(1,0,1,1)$
(lemma \ref{lemma:ContinuityII});\newline
1 iteration to reach $\Lambda$, with code $(1)$.
\end{quote}
Putting everything together, we obtain the code (\ref{eq:SymbolicDynamicsOut}).
In particular, the transit time of the $m$-th atom is odd, and equal to $4m+3$. 
From lemma \ref{lemma:FixedSets}, we then have 
\begin{equation}\label{eq:FixedSetsOut}
{\rm Fix}\, L_m^{\rm out}=F^{2m+1} ({\rm Fix}\,H)\cap \Lambda_m^{\rm out},
\end{equation}
and it remains to show that all these intersections are non-empty. 

We identify the fixed set of an involution within each regular atom.
The atom $\Lambda_{1}^{\rm out}$ is the image under $G\circ F^{-4}$ of the 
set $\Upsilon_1\cap H(\Upsilon_1)$, which intersects ${\rm Fix}\,H$.
Hence $\Lambda_1^{\rm out}$ intersects $G\circ F^{-4}({\rm Fix}\, H)=F^3({\rm Fix}\, H)$,
in accordance with (\ref{eq:FixedSetsOut}).
Likewise, the atom $\Lambda_{2}^{\rm out}$ is the image under $G\circ F^{-4}$ of
the set $\Upsilon_2\cap F^{-4}\circ H(\Upsilon_2)$ which intersects 
$F^{-2}({\rm Fix}\, H)$. 
Hence $\Lambda_2^{\rm out}$ intersects 
$G\circ F^{-6}({\rm Fix}\, H)=F^5({\rm Fix}\, H)$, as desired.

The symmetry properties of the atoms $\Lambda_m^{\rm out}$ for $m\geq 3$ are
proved with an argument analogous to that used in the proof of
theorem \ref{theorem:Lin}. Specifically, in place of (\ref{eq:PointsInFixG})
and (\ref{eq:ExtraPointsInFixG}), we have the symmetric points
\begin{equation}\label{eq:PointsInFixH}
\gamma_m'\cap H(\gamma_m')\,\in\, {\rm Fix}\, H\cap \Sigma' \qquad 
   m=0,\ldots,\lfloor (M-1)/2\rfloor
\end{equation}
as well as
\begin{equation}\label{eq:ExtraPointsInFixH}
F^{4m-2}(\gamma_0')\cap H\circ F^{4m-2}(\gamma_0')\,\in\,{\rm Fix}\,H \qquad
   m=1,\ldots,\lfloor (M-1)/2\rfloor.
\end{equation}
These symmetry lines intersect images of atoms, and by application of 
$G\circ F^{-t}$, for a suitable $t$, one ensures that the intersection
in (\ref{eq:FixedSetsOut}) is non empty also for $m>2$. 
We omit the details, but we note that, by construction, these sets intersect
one vertex of the atom, and hence also the opposite vertex, due to symmetry.

All regular atoms are maximal, being bound by images of the discontinuity line.
Then, from lemmas \ref{lemma:MaximalAtoms} and \ref{lemma:Involutivity}, they 
are invariant under the map $L^{\rm out}$.
\endproof

\section{Quantitative results}\label{section:QuantitativeResults}

In this section we supplement the geometrical approach to return-map dynamics 
developed in section \ref{section:ReturnMap}, with a predominantly algebraic one.
To illustrate the technique, we will prove a version of theorem \ref{theorem:Lout}.
As we shall see, the proof is conceptually simple, and the details relatively easy 
to carry out, provided that one is willing to rely heavily on computer assistance 
to manipulate complicated expressions. For the sake of brevity, the polygons 
treated in this section do not include their boundaries.

\subsection{Calculation of atom vertices\label{subsection:atomvertices}}
The elementary ingredients of our algebraic calculations are the regular atoms of 
$L^{\rm in}$ and $L^{\rm out}$ expressed in terms of their vertices. To obtain explicit 
expressions for the latter as rational function in $\lambda$, we rely on formulae 
(\ref{eq:AtomsLin}) and (\ref{eq:AtomsSidesII}), which express the polygonal atoms 
in terms of bounding lines. Our strategy is to express the latter in terms of the 
cartesian coordinates of their endpoints, then calculate the intersection of 
neighbouring lines to obtain the corresponding polygonal vertex. All quantities
are restricted to $\Omega$, that is, we don't make use of the periodicity of the torus.

Before proceeding to the calculation of vertices, let us introduce the convenient notation
$$
F_\iota=F_{(\iota_0 \iota_1\cdots \iota_{k-1})} 
   = F_{(\iota_{k-1})}\circ F_{(\iota_{k-2})}\circ \cdots \circ F_{(\iota_0)}
$$
where $\iota\in\{0,1\}^k$ is an (arbitrary) codeword, and 
$$
F_{(i)}:\R^2\to\R^2\hskip 40pt
\left(\begin{array}{c}x\\y\end{array}\right)\,\mapsto\, 
C\cdot\left(\begin{array}{c}x\\y\end{array}\right) + \left(\begin{array}{c} i\\0\end{array}\right)
$$
where the matrix $C$ was defined in (\ref{eq:Matrix}).
From (\ref{eq:Map}) we see that our piecewise isometry $F$ acts as $F_{(0)}$ (generalized 
rotation about the origin) on the atom $\Omega_0$ and as $F_{(1)}$ (generalized rotation 
about the fixed point $(\frac{1}{2-\lambda},\frac{1}{2-\lambda})$) on $\Omega_1$.
It is crucial to distinguish between the expressions $F_{(\iota_0 \iota_1\cdots \iota_{k-1})}$,
an isometry of the plane, which acts as a generalized rotation about some fixed point,
and $F^k$, a piecewise isometry of the unit square, acting differently on different domains.

Let us calculate the vertices $P_0(m)$ and $P_1(m)$ of the regular atoms of $L^{\rm out}$, 
see equation (\ref{eq:AtomsVerticesII}), where, based on the geometric analysis of section 
\ref{section:ReturnMap}, 
$$
P_0(m)=\xi_m\cap {\bf 1},\qquad P_1(m)=\xi_m\cap {\bf 5}^1,
$$
with
$$
 \xi_m=F_{(1^3 (0,1)^{(2m-1)} 1^2)}([(1,0),(1,1)]).
$$
We define 
\begin{equation}\label{eq:alphatheta}
\alpha=2\pi\rho=\cos^{-1}(\lambda/2)
\hskip 40pt
\theta=\frac{\pi}{2}-\alpha=\sin^{-1}(\lambda/2).
\end{equation} 
From (\ref{eq:Matrix}) we see that $\alpha$ is the angle of rotation
of the matrix C, while $\theta=2\pi(\frac{1}{4}-\rho)$. Hence
$\theta$ is the angular departure from the $\lambda=0$ rotation, that is,
$4\theta$ is the (positive) angle between successive regular components of
the sector map on $\Sigma$, in the ${\cal Q}$-metric.

To evaluate $\xi_m$ explicitly for general $m$, we regard the isometry 
$F_{(1^3 (0,1)^{(2m-1)} 1^2)}$ as a rotation by $3\alpha$ about 
$z_1=(\frac{1}{2-\lambda},\frac{1}{2-\lambda})$, then a rotation by $(4 m-2)\alpha$ 
about $z_2=\left(\frac{2}{4-\lambda^2},\frac{\lambda}{4-\lambda^2}\right)$, and finally a rotation by $2\alpha$ about $z_1$.
A multiple-angle rotation corresponds to a power of the matrix $C$ and can be executed in a 
single step, thanks to the formula (a consequence of the Jordan decomposition of $C$)
\beq\label{eq:Ckgeneral}
C^k = \cos\left(k \alpha\right)\left(\begin{array}{cc}1&0\\0&1\end{array}\right)
+\frac{\sin\left(k \alpha\right)}{\sqrt{4-\lambda^2}} 
\left(\begin{array}{cc}\lambda&-2\\2&-\lambda\end{array}\right).
\eeq
Our geometric constructions of the atom boundaries typically involve powers of $C^4$, in which case (\ref{eq:Ckgeneral}) can be conveniently written in terms of the $O(\lambda)$ angle $\theta$ defined in (\ref{eq:alphatheta}):
\beq\label{eq:C4l}
C^{4\,l} = \cos\left(4\, l \theta\right)\left(\begin{array}{cc}1&0\\0&1\end{array}\right)
-\frac{\sin\left(4\, l \theta\right)}{\sqrt{4-\lambda^2}} 
\left(\begin{array}{cc}\lambda&-2\\2&-\lambda\end{array}\right).
\eeq
Once the segments $\xi_m$ have been determined, the indicated intersections are 
calculated by elementary algebra (a tedious exercise, best done with computer 
assistance\cite{ESupplement}). 
With $M$ as given in (\ref{eq:M}), we have
\begin{eqnarray}
P_0(m)&=&\left(\frac{1}{2}\left(1+2\lambda+ 
   \frac{\tan{\theta}}{\tan(2 m +1)\theta}\right),1\right),\quad 1\leq m\leq M,
\label{eq:P0}\\
P_1(m)&=&\left(\frac{1}{2}\left(1+2\lambda-\lambda^2+\lambda^4
  + (1-3\lambda^2+\lambda^4)\frac{\tan\theta}{\tan(2 m -1)\theta}\right),\right.\nonumber\\
&&\left. 1+\frac{\lambda^3}{2}+\frac{\lambda}{2}(\lambda^2-2)
    \frac{\tan\theta}{\tan(2m-1)\theta}\right),\quad 2\leq m \leq M.
\label{eq:P1}
\end{eqnarray}

We recall that for $m=1,2$, the symmetry axis of $\Lambda_m^{\rm out}$ 
intersects ${\bf 0}$ instead of ${\bf 5}^1$, and $P_1(m)$ was defined 
to be the right-hand endpoint of the symmetry axis (see remarks preceding 
equation (\ref{eq:AtomsVerticesII}), and theorem \ref{theorem:Lout}).
The atom $\Lambda_1^{\rm out}$ is a quadrilateral 
(see equation \ref{eq:AtomsSidesII}), with vertices
\beq\label{eq:mequal1}
\Lambda_1^{\rm out}= [P_0(1),(1,1),P_1(0),P'],
\eeq
where
$$
P_1(0)=\left(1,\frac{2-\lambda-\lambda^2}{2-4\lambda^2+\lambda^4}\right),
\qquad
P'=(1-2\lambda^2+6\lambda^3+\lambda^4-5\lambda^5+\lambda^7,1-\lambda+3\lambda^2+\lambda^3-4\lambda^4+\lambda^6).
$$
The atom $\Lambda_2^{\rm out}$ is the hexagon
\beq\label{eq:mequal2}
\Lambda_2^{\rm out}=[P_0(2),P_0(1),P'',P_1(1),P''',P_1(2)],
\eeq
where
\begin{eqnarray*}
P_1(1)&=&\left(1,\frac{2-2\lambda-4\lambda^2+\lambda^3+\lambda^4}{2-9\lambda^2+6\lambda^4-\lambda^6}\right),\qquad
P''\,=\,\left(1,\frac{1-\lambda-3\lambda^2+\lambda^3+\lambda^4}{1-6\lambda^2+5\lambda^4-\lambda^6}\right),
\\
P'''&=&\left(\frac{1-6\lambda^2-3\lambda^3+11\lambda^4+4\lambda^5-12\lambda^6-\lambda^7+6\lambda^8-\lambda^{10}}{1-6\lambda^2+5\lambda^4-\lambda^6},\right.\\
&&\quad \left.\frac{1-\lambda-4\lambda^2+3\lambda^3+4\lambda^4-7\lambda^5-\lambda^6+5\lambda^7-\lambda^9}{1-6\lambda^2+5\lambda^4-\lambda^6}\right).
\end{eqnarray*}

To calculate the vertices $Q_0(n)$ and $Q_1(n)$ of the regular atoms of $L^{\rm in}$, we apply to the line $\gamma_0={\bf 1}$ a $4n$-fold $C$-rotation about the centre of ${\cal E}$, followed by a $G$-reflection and, finally, intersection with the lines  $\gamma_0$ and $\gamma_1={\bf 5}^1$, respectively.

The results are \cite{ESupplement}
\begin{eqnarray}
Q_0(n)&=&\left(\frac{2-\lambda-\lambda^2 - s_{n}(\lambda) +(2+\lambda) c_{n}(\lambda)}{-\lambda s_{n}(\lambda)+(4-\lambda^2) c_{n}(\lambda)},1\right)
\label{eq:Q0}\\
Q_1(n)&=&\left(\frac{a(\lambda) + b(\lambda) s_n(\lambda) +c(\lambda) c_n(\lambda)}{d(\lambda) s_n(\lambda)+e(\lambda)  c_n(\lambda)},\frac{f(\lambda) + g(\lambda) s_n(\lambda) +h(\lambda)  c_n(\lambda)}{d(\lambda) s_n(\lambda)+e(\lambda)  c_n(\lambda)}\right)
\label{eq:Q1}
\end{eqnarray}

\noindent
with

\hfil\vtop{\baselineskip=15pt\halign{$#$\hfil&\qquad $#$\hfil\cr
a(\lambda)=2 - \lambda - 7 \lambda^2 + 3 \lambda^3 + 5 \lambda^4 - \lambda^5 - \lambda^6
&
b(\lambda)=-1 - 2 \lambda - \lambda^2 + 2 \lambda^3 + \lambda^4
\cr
c(\lambda)=2 + \lambda - 6 \lambda^2 - 3 \lambda^3 + 2 \lambda^4 + \lambda^5
&
d(\lambda)=-\lambda (5 - 5 \lambda^2 + \lambda^4)
\cr
e(\lambda)=(4 - \lambda^2)\,(1 - 3 \lambda^2 + \lambda^4)
&
f(\lambda)=-\lambda (4 - 2 \lambda - 4 \lambda^2 + \lambda^3 + \lambda^4)
\cr
g(\lambda)=\lambda (-3 - \lambda + 2 \lambda^2 + \lambda^3)
&
h(\lambda)=4 - 7 \lambda^2 - 3 \lambda^3 + 2 \lambda^4 + \lambda^5.
\cr
}}\hfil

\noindent
and\newline
\begin{eqnarray*}
c_k(\lambda)&=&\cos(4 k \theta)=T_{4 k}\left(\frac{\lambda}{2}\right),
\\
s_k(\lambda)&=&\sqrt{4-\lambda^2}\sin(4 k \theta)=\left(\frac{\lambda^2-4}{2}\right) U_{4 k -1}\left(\frac{\lambda}{2}\right).
\end{eqnarray*}
Here $T_k$ and $U_k$ are the Chebyshev polynomials of the first and second kinds, 
respectively\cite[p 1032]{GradshteynRyzhik}.

The use of the rotation angle $4 n \theta$ in (\ref{eq:C4l}), as well as in the formulae 
for $Q_0(n)$ and $Q_1(n)$, allows for a natural interpolation of these quantities in the 
range $0\leq n \leq N-1$.  For continuously varying $n$, the matrix $C^{4n}$ defined by 
(\ref{eq:C4l}) (not (\ref{eq:Ckgeneral})) represents a family of $C$-rotations 
whose angle increases monotonically from $0$ to $4 N\theta = \pi/2-O(\lambda)$.  
Moreover, as $n$ increases from $0$ to $N-1$, $Q_0(n)$ and $Q_1(n)$ smoothly 
trace out, in order-preserving fashion, the upper and lower boundaries of the regular part of $\Lambda$.  

Note that formulae (\ref{eq:Q0}) and (\ref{eq:Q1}) explicitly specify the invertible 
transformations between parameter $n$ and points $X=(x,y)$ on the lines $\gamma_0$ and 
$\gamma_1$, respectively.  This gives precise meaning to the functional notation 
$n(X), X\in\Sigma$, which we will employ below. In section \ref{section:PeriodicPoints} 
we will exploit this parametrization to study the relative positions on the two boundary 
lines of the points $P_0(m)$ and $Q_0(n)$ (resp. $P_0(m)$ and $Q_0(n)$).

\subsection{A preliminary lemma} 
The following lemma plays an important role in the theorem which follows.

\begin{lemma}\label{lemma:triangles}
The map $F^4$ is an isometry on the triangle $[(0,1),(1,1),(\lambda^2,1-\lambda)]$, with itinerary
\begin{eqnarray*}
&&[(0,1),(1,1),(\lambda^2,1-\lambda)]\stackrel{F_1}{\rightarrow}
[(0,0),(0,1),(\lambda,\lambda^2)]\stackrel{F_1}{\rightarrow}
[(0,0),(1,\lambda),(1,0)]
\\
&&\quad \stackrel{F_0}{\rightarrow} [(0,0),(0,1),(\lambda,1)]\stackrel{F_1}{\rightarrow}
[(0,0),(\lambda^2,\lambda),(1,0)].
\end{eqnarray*}
The triangles $F^k([(0,1),(1,1),(\lambda^2,1-\lambda)]),\; k=1,\ldots,4$, do not intersect the 
sector $\Sigma$.
\end{lemma}

\proof The proof is by direct mapping of the vertices, accompanied by straightforward 
verification that the interior each triangle is contained in a single atom of $F$, 
the latter being correctly identified in the statement of the lemma. 
The empty intersection of the four image triangles with $\Sigma$ is trivially verified.
\subsection{Algebraic alternative to theorem \ref{theorem:Lout}}
We now state an alternative to theorem \ref{theorem:Lout} in which the polygons 
$\Lambda_m^{\rm out}$ are defined explicitly in terms of their respective vertices, 
as in section \ref{subsection:atomvertices}. That these are symmetric atoms of 
$L^{\rm out}$ is proved by an inductive calculation of their return paths.
\begin{theorem}
For $1\leq m \leq M$, we have
$$
L^{\rm out}(\Lambda^{\rm out}_m) =F^{4m + 3}\circ G (\Lambda^{\rm out}_m) = F_{(1^3(0,1)^{2m-1}1^2)}\circ G (\Lambda^{\rm out}_m)=\Lambda^{\rm out}_m.
$$
Moreover,
$$
[P_0(m),P_1(m-1)]={\rm Fix}\,L_m^{\rm out}.
$$
\end{theorem}
\proof We begin by defining some auxiliary polygons related to $\Lambda^{\rm out}_m$:
$$
\begin{array}{ll}
\hat{\Lambda}^{\rm out}_m&\stackrel{\rm def}{=}\,\bigcup_{k\geq m}^M \Lambda^{\rm out}_m \cup
 [P_0(M),P_1(M),(g,1)],\\ 
\\
\check{\Lambda}^{\rm out}_m&\stackrel{\rm def}{=}\,
F_{(1^3(0,1)^{2m-1})}\circ G (\hat{\Lambda}^{\rm out}_m).
\end{array}
$$
We want to prove that $G( \hat{\Lambda}_m)$ is mapped isometrically onto $\check{\Lambda}_m$ by $F^{4m+1}$.

The proof is by induction on $m$. For $m=1$, the isometric property follows from lemma 
\ref{lemma:triangles} by showing that $F(G \hat{\Lambda}^{\rm out}_1)$ is contained in 
$[(0,1),(1,1),(\lambda^2,1-\lambda)]$. To verify this, we note that
$$
F(G \hat{\Lambda}^{\rm out}_1)\subset F( [(g'',1),(1,1),(1,g)])=[(\lambda g'',g''), (\lambda,1),(1+\lambda-g,1)].
$$
It is a simple exercise \cite{ESupplement} (sec. E.5.1)to show that for $0<\lambda<1$, all three vertices of the 
image triangle are within or on the boundary of $[(0,1),(1,1),(\lambda^2,1-\lambda)]$.

We now assume that the induction hypothesis $F^{4m+1}\circ G (\hat{\Lambda}_m) =\check{\Lambda}_m$  
holds for a given $m<M$. Note that $G( \hat{\Lambda}_m)$ is partitioned by the line 
segment $G ([P_0(m),P_1(m)])$ ($m>1$) or $G([P_0(1),P'])$ ($m=1$)
into adjacent polygons $G (\Lambda^{\rm out}_m)$ and $G(\hat{\Lambda}_{m+1})$. We denote the images 
of these polygons $\check{\Lambda}_m^+$ and $\check{\Lambda}_m^-$,respectively.
One checks \cite{ESupplement} (secs. E.5.2, E.5.3) that the boundary separating the two is a subset of 
$[(0,0),(1,\lambda)]$, so that 
$$
\check{\Lambda}_m^-\subset\Omega_0, \qquad 
\check{\Lambda}_m^+\,\subset\, ([(0,0),(\lambda^2,\lambda),(1,0)]\setminus \Omega_0)\,\subset\, 
(\Omega_1\setminus\Sigma).
$$

To complete the proof of 
\beq\label{eq:F4m+1}
F^{4m+1} \circ G( \hat{\Lambda}_{m+1}) =\check{\Lambda}_{m+1},
\eeq
we need to show that $F^4$ maps $\check{\Lambda}_m^-$ isometrically onto $\check{\Lambda}_{m+1}$, 
with itinerary (0,1,0,1).  Using lemma \ref{lemma:triangles} and $\check{\Lambda}_m^-\subset\Omega_0$, 
we have 
\beq\label{eq:inclusionA}
F^2(\check{\Lambda}_m^-) = F_{(01)}(\check{\Lambda}_m^-) \subset \Omega_0,
\eeq
where the inclusion is verified explicitly \cite{ESupplement} (sec. E.5.4). Again using lemma \ref{lemma:triangles}, 
we have 
$$
F^4(\check{\Lambda}_m^-) = F_{(0101)}(\check{\Lambda}_m^-).
$$
The right-hand side is identical to $\check{\Lambda}_{m+1}$ by construction.

Having established (\ref{eq:F4m+1}), we need only show that $F^2$  maps $\check{\Lambda}_m^+$ isometrically onto $\Lambda^{\rm out}_m$ with itinerary $(1,1)$. For the isometry, we verify the inclusion \cite{ESupplement} (sec. E.5.5)
\beq\label{eq:inclusionB}
F_1(\check{\Lambda}_m^+)\subset(\Omega_1\setminus \Sigma).
\eeq
The rest is an explicit calculation\cite{ESupplement} (sec. E.6) of the image points of the vertices, 
using elementary algebra, 
In particular, we find that the map $L_m^{\rm out}=F^{4m+3} \circ G$ acts
on the vertices of the polygons $\Lambda^{\rm out}_m$ as follows 
(notation of section \ref{subsection:atomvertices}):
\begin{eqnarray*}
&&P_0(m)\mapsto P_0(m), \quad P_1(m-1)\mapsto P_1(m-1), \quad 1\leq m\leq M,
\\
&&P_0(m-1)\leftrightarrow P_1(m),  \quad 2\leq m\leq M,
\\
&&(1,1)\leftrightarrow P',\quad  m=1,\hskip 40pt  P''\leftrightarrow P''',\quad m=2.
\end{eqnarray*}

The final identification of $F^{4m+3}\circ G$ as the first-return map $L^{\rm out}$, with regular 
atoms $\Lambda_m^{\rm out}$ requires verification that all polygons $F^k\circ G(\Lambda_m^{\rm out})$, 
$1\leq k<4m+3$, are disjoint from $\Sigma$. We have taken care to establish this property at each 
stage of the proof (specifically, in lemma \ref{lemma:triangles} and equations (\ref{eq:inclusionA}) 
and (\ref{eq:inclusionB})).
  
Recalling that the map $L^{\rm out}$ acts linearly on $\Lambda^{\rm out}_m$, the vertex mappings 
are sufficient to establish that the atom is  $L^{\rm out}$-symmetric with fixed segment $[P_0(m),P_1(m-1)]$.
\endproof

\section{Regular atoms and fixed points of $L$}
\label{section:PeriodicPoints}
\subsection{Overview}
In section \ref{section:ReturnMap} we decomposed the return map $L$ into a product of 
involutions $L^{\rm in}$ and $L^{\rm out}$.  
According to theorem \ref{theorem:Lin}, $L^{\rm in}$ is endowed with a sequence of 
regular atoms $\Lambda_n^{\rm in}, \,n=1,\ldots,N-1$, each of which is an 
$L^{\rm in}$-symmetric polygon.  
Theorem \ref{theorem:Lout} establishes an analogous sequence of regular atoms, 
$\Lambda_m^{\rm out},\, m=1,\ldots,M$ for $L^{\rm out}$.  
We now define {\bf regular atoms} of $L$ to be the nontrivial intersections of the 
regular atoms of the two involutions:
\begin{equation}\label{eq:AtomsL}
\Lambda_{m,n}=\Lambda_m^{\rm out}\,\cap\,\Lambda_n^{\rm in}.
\end{equation}

For some, but not necessarily all, of these atoms, the symmetry lines of $L^{\rm in}$ and
$L^{\rm out}$ in $\Lambda_n^{\rm in}$ and $\Lambda_m^{\rm out}$, respectively, intersect
within $\Lambda_{m,n}$, at a symmetric return-map fixed point $Z(m,n)$.
Such fixed points will be called {\bf regular,} and will be the main focus of this and
succeeding sections. They correspond to a 2-parameter family of periodic orbits of $F$,
with periods
\begin{equation}\label{eq:Periods}
t(m,n)=4(m+n)-1
\end{equation}
and symbolic codes
\begin{equation}\label{eq:Codes}
\iota^{(n,m)}=\bigl(\overline{1^{4n-1},(0,1)^{2m-1},1,1}\bigr),
\end{equation}
the latter obtained by concatenating the codes of the involutions $L^{\rm in,out}$
given in theorems \ref{theorem:Lin} and \ref{theorem:Lout}.

\begin{figure}
\begin{center}
\epsfig{file=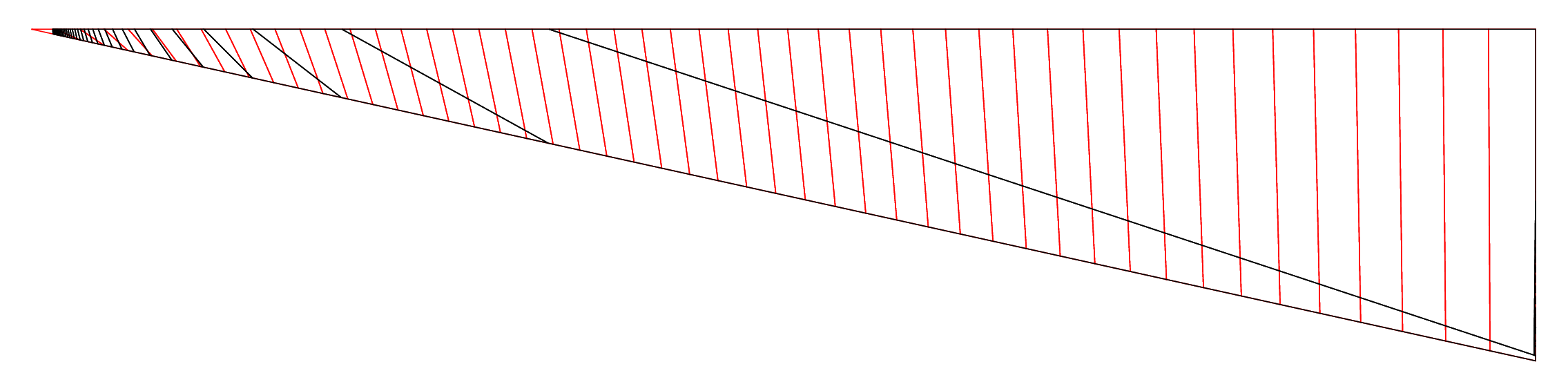,width=12cm}
\end{center}
\begin{center}
\epsfig{file=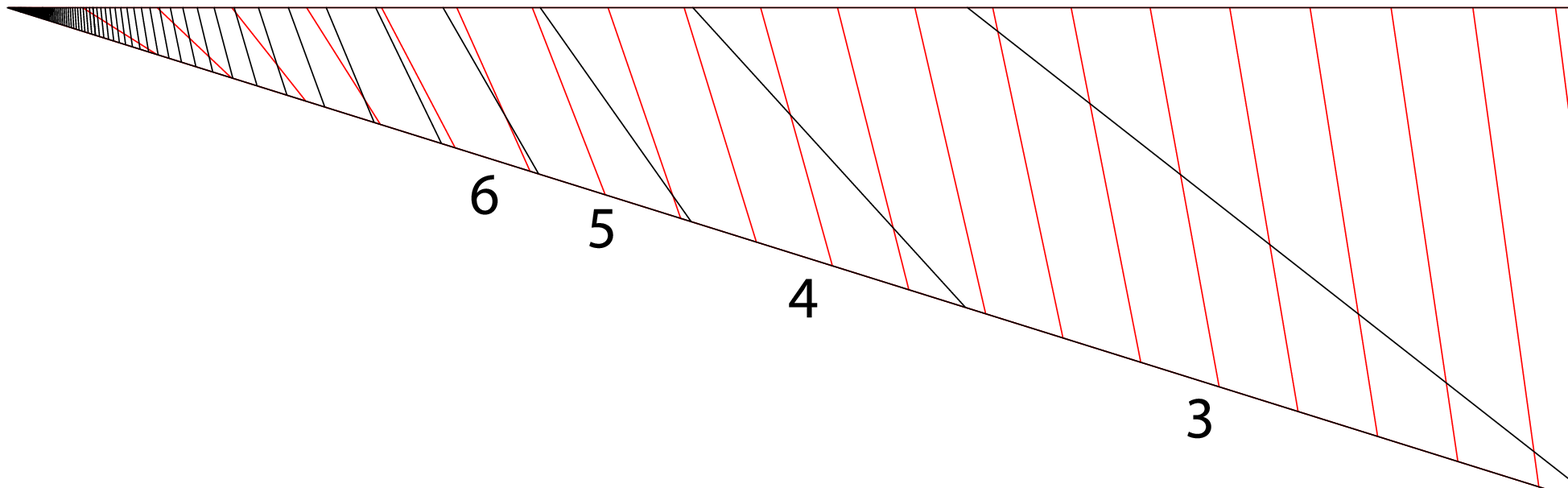,width=12cm}
\end{center}
\caption{\label{fig:Atoms}
\small Top: Superposition of the atoms of the domain $\Lambda$ into atoms. 
The atoms of $L^{\rm in}$ are in red; those of $L^{\rm out}$ in black.
Bottom: Magnification of the region near the tip of $\Lambda$, showing
the crossover phenomenon. The integers denote values of $m$, and crossover
occurs for $m=6$.}
\end{figure}
\begin{figure}[t]
\hfil\epsfig{file=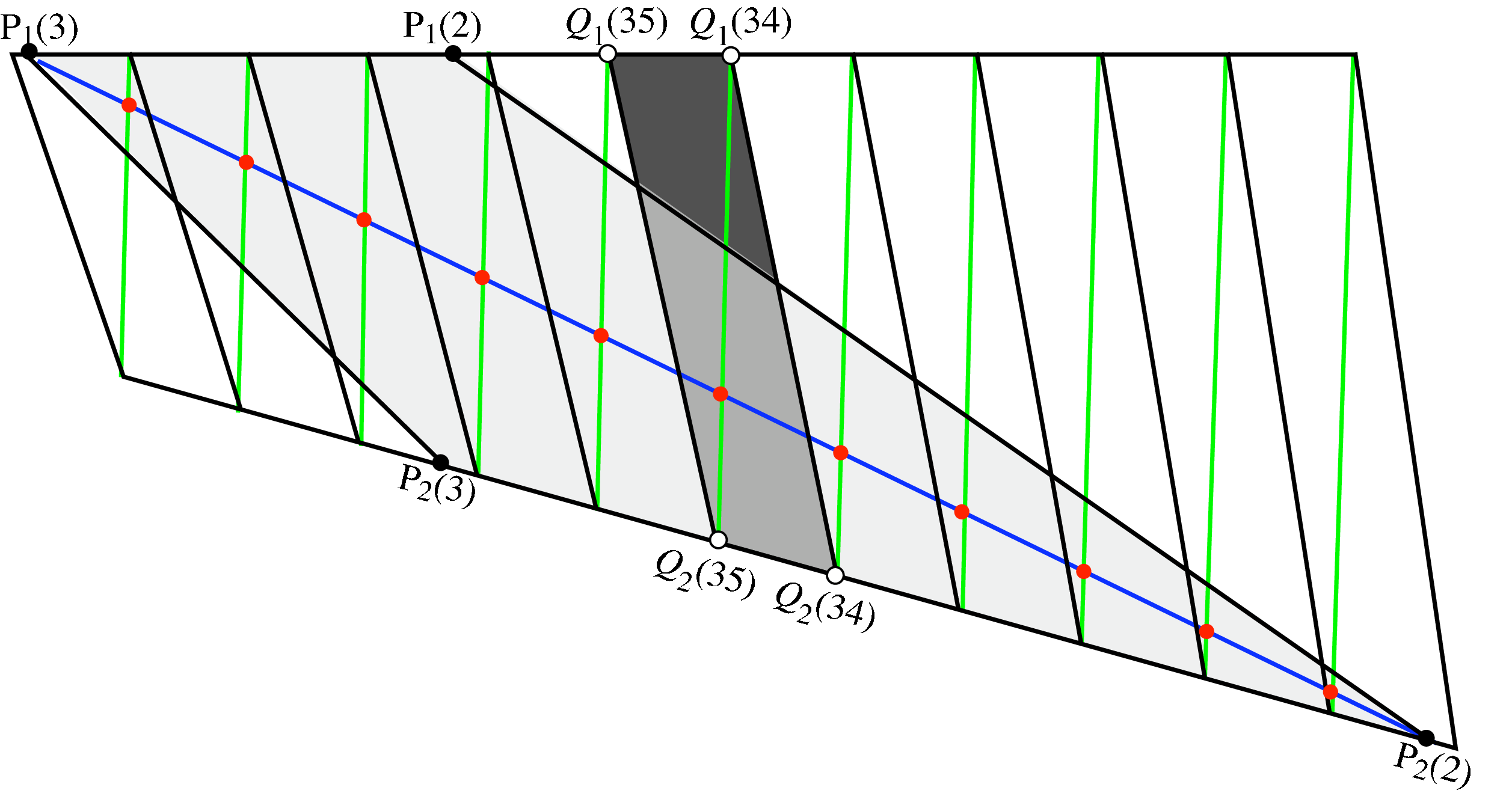,width=8cm}\hfil
\caption{\label{fig:FixedPoints}
\small The fixed points of the return map $L$ are the intersection of the 
fixed lines of the involutions $L^{\rm in,out}$, which are segments connecting
opposite vertices of the atoms $\Lambda_n^{\rm in}$ 
(vertices $Q_0(n),Q_0(n-1),Q_1(n-1),Q_1(n)$),
and $\Lambda_m^{\rm out}$
(vertices $P_0(m),P_0(m-1),P_1(m-1),P_1(m)$).
The example focuses on $m=3, n=34$ for $\lambda=\frac{1}{64}$. }
\end{figure}

Before examining in detail the geometric properties of the regular atoms of $L$, it is 
useful to introduce a broad classification according to the size of the index $m$ 
relative to its maximum value $M=M(\lambda)$, given by equation (\ref{eq:M}).
A typical example is displayed in figures \ref{fig:Atoms} and 
\ref{fig:FixedPoints}, where $\lambda$ has been assigned the value $\frac{1}{64}$. 
Figure \ref{fig:Atoms} shows the regular atoms of $L^{\rm in}$ (lighter boundaries) and 
$L^{\rm out}$ (darker boundaries) in all of $\Lambda$ (note that the aspect ratio has 
purposely been distorted to make the structural relationships visible).  
The bottom figure \ref{fig:Atoms} zooms in on the left-hand end of $\Lambda$, while figure 
\ref{fig:FixedPoints} zooms in on the atoms with $m=3,\, 30\leq n \leq 40$, showing 
explicitly ${\rm Fix}\,L^{\rm in}$ and ${\rm Fix}\,L^{\rm out}$ and the regular 
fixed points of the return-map where these symmetry lines intersect.  

We observe some interesting features of the regular atoms and their symmetry lines. 
The orientations of the `vertical' boundary lines change substantially as one traverses $\Lambda$.
The $\Lambda_n^{\rm in}$ boundaries are vertical at the right-hand edge but suffer significant 
counterclockwise rotations as one moves toward the left-hand vertex.  
The $\Lambda_n^{\rm out}$ boundaries, on the other hand, rotate in the opposite sense.  
Somewhere in the middle, there is a crossover region where the boundaries of the two 
sets of atoms are approximately aligned (this corresponds to $m=6$ in figure 
\ref{fig:Atoms}).  
To the right of the crossover region, each $\Lambda_m^{\rm out}$ intersects a number 
of successive $\Lambda_n^{\rm in}$ atoms, and each of the intersections $\Lambda_{m,n}$ 
(the regular atoms of $L$) contains a unique regular fixed point $Z(m,n)$.   
To the left of the crossover region, on the other hand, each $\Lambda_m^{\rm in}$ 
intersects a number of successive $\Lambda_n^{\rm out}$ atoms, and now there are 
regular atoms $\Lambda_{m,n}$ which do not contain a regular fixed point.  
This can be seen quite clearly in the case $m=8$, which has a triangular intersection 
region in its northeast corner which is not intersected by the diagonal 
${\rm Fix}\,L_8^{\rm out}$.  

In appendix A.7 we give a precise definition of the crossover point $m=m^*, n=n^*$ and 
obtain a perturbative formula for $m^*$ and $n^*$ as functions of $\lambda$.  
In particular, we will find
$$
m^*=\frac{1}{\sqrt{2\lambda}} \left(1 + \frac{\lambda}{3} +\cdots\right),\qquad
n^*=\frac{\pi}{4\lambda}-\frac{1}{\sqrt{2\lambda}}-\frac{3}{4} +\cdots.
$$
We note that in the limit $\lambda \rightarrow 0^+$, the crossover value of $m$ tends to infinity.  
Thus, if we are interested in the asymptotic behaviour of regular fixed points and their 
stability ellipses for fixed positive integer $m$, we can limit our attention to those 
regular atoms which are well to the right of the crossover, and avoid the complexities 
of the crossover phenomenon. This is achieved by introducing a cut-off condition
\begin{equation}\label{eq:CutOff}
m\leq \lambda^{-\nu_0}, \quad\mbox{\rm where}\quad \quad 0<\nu_0<\frac{1}{2}.
\end{equation}
We will prove that the regular atoms satisfying this condition have a number of nice 
properties, which allow us to obtain explicit formulae for their total area, 
whereas those with $m$ exceeding the cut-off have a total area tending to 
zero as $\lambda^{\frac{1}{2}-\nu_0}$, and can therefore be neglected.

In the remainder of this section, we study in detail the properties of 
the regular atoms of $L$.
Our starting point is the list of explicit formulae for the atomic
vertices given in section \ref{subsection:atomvertices}. 
The vertex coordinates are rational functions of $\lambda$; as $\lambda\to 0$, 
we will consider Taylor expansions in $\lambda$, with suitable remainder. 
Our most complete description of the atoms $\Lambda_{m,n}$ will be restricted to 
those satisfying cut-off conditions $\lambda\leq \lambda_0, \, m\leq \lambda^{-\nu_0}$.
Our choice of cut-off parameters is (cf.~(\ref{eq:CutOff}))
\begin{equation}\label{eq:CutOffParameters}
\lambda_0=10^{-4}\quad \mbox{\rm and}\quad \nu_0=\frac{1}{4}.
\end{equation}
Even though a specific choice is required to obtain explicit estimates, our qualitative 
results are insensitive to these special values, and the detailed calculations included 
in the Electronic Supplement \cite{ESupplement} can readily be adapted to other cut-off values.
For the rest of this section, we shall assume that $\lambda_0$ and $\nu_0$ are given by 
(\ref{eq:CutOffParameters}).

The organization of the regular atoms of $L^{\rm in, out}$ and their intersections 
$\Lambda_{m,n}$ is determined by the order and spacing of the 
vertices along the lines ${\bf 1}$ and ${\bf 5}^1$. 
Often it will be convenient to 
regard $n$ as a continuous variable on the interval $[0,N]$, as discussed in section 
\ref{subsection:atomvertices}, and to express our results in terms of auxiliary variables 
\begin{equation}\label{eq:etatau}
\eta=\frac{\pi}{2}-2 n \theta,\qquad  \tau=\tan(\eta).
\end{equation}
Note that
as $n$ increases from 0 to $N$, the quantity $\tau$ decreases 
monotonically from 1 to $\tan(\pi/4-2 N \theta)=O(\lambda)$.

In section \ref{subsection:atomvertices}, we introduced the notation $n(X)$ to designate 
the real parameter $n$ assigned to a point $X$ on one of the lines ${\bf 1}$ or ${\bf 5}^1$. 
In the same spirit, we will employ $\tau(X)$ to denote the corresponding $\tau$-variable, 
namely $\tan(\pi/4-2 n(X)\theta)$.

\subsection{Expansions for the atom vertices}
We begin our structural analysis by determining the positions of the vertices
$P_{0,1}(m), Q_{0,1}(n)$, correct to next-to-leading order in powers of $\lambda$.  
The following lemma establishes, among other things, that  the horizontal positions and 
widths of the atoms $\Lambda_m^{\rm out}$ become fixed as $\lambda$ tends to zero, 
while the vertical dimension shrinks proportional to $\lambda$.  
For the atoms $\Lambda_n^{\rm in}$, on the other hand, both dimensions shrink 
proportional to $\lambda$.  In addition, we see that the $x$-coordinates of $P_0(m-1)$ 
and $P_1(m)$ coincide up to first order in $\lambda$.  
These properties are apparent in the figures. 
The proof of the lemma is given in appendices A.2 and A.3.

\begin{lemma}\label{lemma:PQexpand}
For  $\lambda\leq\lambda_0$, and $2\leq m\leq \lambda^{-\nu_0}$, 
with $\lambda_0, \nu_0$ given in (\ref{eq:CutOffParameters}), the cartesian 
coordinates of the vertices of the regular atoms $\Lambda_m^{\rm out}$ and 
$\Lambda_n^{\rm in}$ are given by
\begin{eqnarray}
P_0(m)&=& \left(\frac{m+1}{2 m+1} + \lambda -\frac{1}{2m+1}\Big(\frac{m(m+1)}{6}\lambda^2+ r_{0x} \lambda^3\Big),\, 1\right)
\nonumber
\\
P_1(m)&=&\left(\frac{m}{2 m-1} +\lambda -\frac{1}{2m-1}\Big(\frac{(m+2)(m+3))}{6}\lambda^2 +r_{1x}\lambda^3\Big),\right.
\\
&& \quad \left.1+\frac{1}{2m-1}\Big(- \lambda + \frac{m(m+2)}{3}\lambda^3+r_{1y}\lambda^4\Big)\right)
\nonumber
\\
Q_0(n)&=&\left(\frac{1}{2}(1+\tau) +\frac{\lambda}{8}(3-2\tau-\tau^2) +\frac{\lambda^2}{32}(2-5\tau+2\tau^2+\tau^3)+r_{2x}\lambda^3,1\right)
\label{eq:PQexpand}
\\
Q_1(n)&=&\left(\frac{1}{2}(1+\tau) +\frac{\lambda}{8}(7-2\tau-5 \tau^2) +\frac{\lambda^2}{32}(-6-29\tau+10\tau^2+25\tau^3)+r_{3x}\lambda^3,\right.
\nonumber
\\
&& \quad \left. 1-\tau\lambda+\frac{\lambda^2}{4} (1+2 \tau+5 \tau^2)-\frac{\lambda^3}{16}(2+11\tau+10\tau^2+25\tau^3)+r_{3y} \lambda^4\right),
\nonumber
\end{eqnarray}
where $\tau$ is defined in (\ref{eq:etatau}),
and the remainders $r$ are rational functions of $\lambda^{1/2}$ that
admit the following bounds
$$
\vcenter{\halign{
$\displaystyle #$\hfil\quad&
$\displaystyle #$\hfil\quad&
$\displaystyle #$\hfil\cr
r_{0x}\in \left[-\frac{1697}{10^8},\frac{1507}{10^5}\right], &
r_{1x}\in \left[-\frac{8523}{10^6},\frac{1243}{10^5}\right],&
r_{1y}\in\left[-\frac{5449}{10^6},\frac{1203}{50000}\right],\cr
r_{2x}\in\left[-\frac{2051}{5000},\frac{4499}{10000}\right],&
r_{3x}\in\left[-\frac{23}{20},\frac{389}{250}\right],&
r_{3y}\in\left[\frac{519}{5000},\frac{539}{100}\right].\cr
}}
$$
In addition, we have
$$
\begin{array}{lll}
P_0(1)=\Big(\frac{2}{3} + \lambda -\frac{1}{9}\lambda^2+ r_{01x} \lambda^4,\, 1\Big),
&\quad&
r_{01x}\in \left[-\frac{463}{12500},-\frac{3703}{100000}\right]
\\
P_1(0)= \Big(1, 1-\frac{1}{2}\lambda+\frac{3}{2}\lambda^2-\lambda^3+r_{10y}\lambda^4\Big),
&\quad&
r_{10y}\in\left[\frac{1249}{500},\frac{313}{125}\right]
\\
P_1(1)= \Big(1,1-\lambda+\frac{5}{2}\lambda^2-4\lambda^3 + r_{11y}\lambda^4\Big),
&\quad&
r_{11y}\in \left[\frac{8739}{1000},\frac{8763}{1000}\right].
\end{array}
$$
\end{lemma}

\begin{corollary}
For  $\lambda\leq\lambda_0$ and $2\leq m\leq \lambda^{-\nu_0}$, we have
$$
P_0(m)_x< P_1(m)_x<P_0(m-1)_x<P_1(m-1)_x,
$$
with 
$$
\vcenter{\baselineskip 18pt\halign{
\quad 
 $#$ \hfil&\qquad
 $#$ \hfil\cr
P_1(m)_x-P_0(m)_x=O(1), & P_1(m)_x-P_0(m-1)_x=O(\lambda^2).\cr
}}
$$
For $m=1$, 
$$
P_0(1)<P_1(0),\qquad
P_1(0)-P_0(1)=O(1).
$$
\end{corollary}

In addition to the information about the horizontal ordering of cartesian coordinates, 
we will need analogous results with respect to $n$, considered as a continuous variable.
The proof is found in appendix A.5.
\begin{lemma}\label{lemma:nP}
For $\lambda\leq \lambda_0$ and $2\leq m\leq \lambda^{-\nu_0}$, we have
\beq\label{eq:nPinequalities}
n(P_0(m))> n(P_1(m))>n(P_0(m-1))>n(P_1(m-1)),
\eeq
and the quantities
$\lambda (n(P_0(m))-n(P_1(m)))$ and $ \lambda^{-1}(n(P_1(m))-n(P_0(m-1)))$
have positive lower and upper bounds which are independent of $\lambda$.
In the case $m=1$, we have
$$
n(P_0(1))>n(P_1(0))=0,
$$
and the quantity $\lambda (n(P_0(1))-n(P_1(0)))$ has positive lower and upper bounds which are independent of $\lambda$.
\end{lemma}

\subsection{Structural theorem}
We are now in a position to establish the most important structural properties 
of the regular atoms of the return-map $L$.
\begin{theorem}\label{theorem:atomintersections}
If the inequalities 
\begin{equation}\label{eq:Inequalities}
P_0(m)_x<Q_0(n-1)_x  \quad \mbox{ and } \quad P_1(m-1)_x>Q_1(n)_x 
\end{equation}
are satisfied for some integers $m,n$ in the range $1\leq m\leq M$, $1\leq n<N$,
then the regular atom $\Lambda_{m,n}=\Lambda_n^{\rm in}\cap\Lambda_m^{\rm out}$ 
has positive area, and contains a regular fixed point $Z(m,n)$.
The stability ellipse ${\cal E}_{m,n}$ surrounding $Z(m,n)$ is tangent to (at least) 
three sides of $\Lambda_{m,n}$, one of which is a subset of either $\bf{1}$ or $\bf{5}^1$. 

Given $0<\lambda\leq \lambda_0$ and $1\leq m\leq \lambda^{-\nu_0}$, then
there is a non-empty set of values of $n$ for which the inequalities 
(\ref{eq:Inequalities}) hold, given by
$$
n_-(m) \leq  n \leq n_+(m), \quad n_-(m)=\lceil n(P_1(m-1))\rceil, \quad
n_+(m)=\lceil n(P_0(m))\rceil.
$$
As $\lambda\to 0$, such a set becomes infinite.
\end{theorem}

\proof From equations (\ref{eq:AtomsVertices}) and (\ref{eq:AtomsVerticesII})
we see that the points $P_0(m)$ and $Q_0(n-1)$ belong to the 
segment ${\bf 1}$, which is parallel to the $x$-axis. 
Suppose that $m>2$. Then the points $P_1(m-1)$ and $Q_1(n)$ belong 
to ${\bf 5}^1$, which is quasi-parallel to the $x$-axis. 
Then the inequalities (\ref{eq:Inequalities}) imply a corresponding 
ordering of the $P$s and $Q$s along the respective segments. 
It follows that the segments $[P_0(m),P_1(m-1)]$ and $[Q_0(n-1),Q_1(n)]$ 
intersect transversally at a point $Z(m,n)$, which lies in the interior 
of both atoms, hence in the interior of their intersection. 
This establishes the existence of a regular atom $\Lambda_{m,n}$ of positive measure.  
From theorems \ref{theorem:Lin} and \ref{theorem:Lout} we recognize 
$[P_0(m),P_1(m-1)]$ as the symmetry axis for $L_m^{\rm out}$, and 
$[Q_0(n-1),Q_1(n)]$ as the symmetry axis for $L_n^{\rm in}$, so that the 
intersection point $Z(m,n)$ is indeed a regular fixed point of $L$.
If $m=1,2$, the symmetry axis of $L_{1,2}^{\rm out}$ connects $P_0(m)$
to ${\bf 0}$; the second inequality in (\ref{eq:Inequalities}) is
automatically satisfied, and the first inequality in (\ref{eq:Inequalities})
suffices to establish transversal intersection of symmetry lines.

Because the atoms $\Lambda_m^{\rm out}$ and $\Lambda_n^{\rm in}$ are maximal, 
so is $\Lambda_{m,n}$. Since the fixed point $Z(m,n)$ lies in the interior
of $\Lambda_{m,n}$, the associated ellipse ${\cal E}_{m,n}$ is necessarily 
tangent to the boundary of $\Lambda_{m,n}$. In particular, $Z(m,n)$ is stable.

We must now establish the number of points of tangency.
Let $T$ be a point at which $\Lambda_{m,n}$ is tangent to ${\cal E}_{m,n}$.
To make our argument independent from the details of the boundary, we
use the expression $L_{m,n}^{\rm in, out}(T)$ to denote the limit of 
$L_{m,n}^{\rm in, out}(z)$ as $z$ approaches $T$ from the interior of the atom. 
Now, any $L$-invariant ellipse in the interior of ${\cal E}_{m,n}$ is also 
invariant under the individual involution $L_n^{\rm in}$ and $L_m^{\rm out}$.
To see this, note that, generically (irrational rotation number), such
an ellipse contains a symmetric dense orbit, which is mapped into itself 
by one involution, hence by the other. By continuity, this
result then extends to the non-generic case (rational rotations).
Because of this property, if $T$ is tangent to ${\cal E}$, so are
$L_n^{\rm in}(T)$ and $L_m^{\rm out}(T)$.

Let us assume that $(m,n)\not=(1,1)$ or $(2,1)$, which are special 
cases to be dealt with separately.
We distinguish three possibilities, depending on whether $\Lambda_{m,n}$ is 
tangent to ${\bf 5}$, to ${\bf 1}$, or to both\footnote{To lighten up the 
notation, we'll omit all superscripts identifying regular components.}
(see figure \ref{fig:FixedPoints}).
If $\Lambda_{m,n}$ is tangent only to ${\bf 5}$, 
then $T\not\in {\bf -4n+4}$, because $L_n^{\rm in}({\bf -4n+4})={\bf 1}$, which 
lies outside the closure of the atom.
The two involutions exchange the remaining three boundaries as follows
$$
{\bf -4n}\hskip 10pt\stackrel{L_n^{\rm in}}{\longleftrightarrow}\hskip 10pt {\bf 5}
\hskip 10pt\stackrel{L_m^{\rm out}}{\longleftrightarrow}\hskip 10pt {\bf 4m-1}
$$
and since the images of $T$ lie on the boundary of ${\cal E}$, we get one
tangency point on each of the above three sides.
In particular, if $\Lambda_{m,n}$ is a triangle, then ${\cal E}_{m,n}$
is tangent to all its sides.
In the case $(m,n)=(2,1)$, the atom $\Lambda_{2,1}$ has the
additional involutory relation ${\bf 12}=L_1^{\rm out}({\bf 0})$.
Since $T$ cannot be tangent to ${\bf 0}$, then it cannot be 
tangent to ${\bf 12}$. But then $\Lambda_{2,1}$ is necessarily a 
pentagon (i.e., ${\bf 5}$ is a side of the atom), for otherwise 
the remaining two sides would map outside the atom, and we would have
no tangency, contradicting maximality. So we still get three points
on ${\bf 5}$, ${\bf 7}$ and ${\bf 12}$.

An analogous argument applies to the case in which $\Lambda_{m,n}$ is 
tangent to ${\bf 1}$ but not to ${\bf 5}$.
Now the ellipse ${\cal E}$ cannot be tangent to ${\bf -4n}$, and 
the involutions act on the remaining sides as follows
$$
{\bf -4n+4}\hskip 10pt\stackrel{L_n^{\rm in}}{\longleftrightarrow}\hskip 10pt {\bf 1}
\hskip 10pt\stackrel{L_m^{\rm out}}{\longleftrightarrow}\hskip 10pt {\bf 4m+3}
$$
giving again one tangency point on each of the above three sides.
In the special cases $(m,n)=(1,1)$, we have the additional involutory 
relation ${\bf 8}=L_1^{\rm out}({\bf 0})$.
Again, $\Lambda_{1,1}$ is necessarily a pentagon (lest we would have no tangency).
This time the ellipse ${\cal E}_{1,1}$ is tangent to ${\bf 0}$,
and so we obtain a fourth point of tangency to ${\bf 8}$.

It remains to consider the case in which $\Lambda_{m,n}$ is tangent to 
both ${\bf 1}$ and ${\bf 5}$. We cannot have the quadrilateral atom 
$\Lambda_{m,n}\not=\Lambda_n^{\rm in}$, because its boundaries 
${\bf 1}$ and ${\bf 5}$ are mapped outside the closure of the 
atom by $L_m^{\rm out}$, and mapped to the other two boundaries 
by $L_n^{\rm in}$. But this would imply the absence of tangency points, 
contrary to the maximality of the atom.
For a similar reason $\Lambda_{m,n}\not=\Lambda_m^{\rm out}$.
Thus $\Lambda_{m,n}$ is either a pentagon or a hexagon. In the former 
case, symmetry still forces the same three tangency points as above,
one of which is either ${\bf 1}$ or ${\bf 5}$. However, ${\cal E}$
cannot be tangent to both, since one of these sides is mapped outside the atom
by $L_m^{\rm out}$.
In the hexagonal case there are two possible configurations with
three tangency points, of which at least one must be realized.
Each configuration must be realized in some parametric interval, and 
by changing parameter the two configurations exchange role.
Then there must exist a parameter value corresponding to six tangency points, 
the most symmetric configuration. 

Let us now assume that the parameters $\lambda$, $m$, and $n$ belong to the given range.
Suppose $P_0(m)_x\geq Q_0(n-1)_x$. Then $\tau(P_0(m))\geq \tau(Q_0(n-1))$.
For $m=1$, the fact that $Q_1(n-1)$ lies on ${\bf 5}$ implies that the right-hand edge 
$[Q_0(n-1),Q_1(n-1)]$ of $\Lambda_n^{\rm in}$ lies 
to the left of the interior of $\Lambda_1^{\rm out}$ and there is no nontrivial intersection.   
Assume therefore $m\geq2$.
Since $\tau(Q_0(n-1))=\tau(Q_1(n-1))$, and, using (\ref{eq:nPinequalities})
and monotonicity, $\tau(P_1(m))>\tau(P_0(m))$, we conclude that 
$\tau(P_1(m))\geq\tau(Q_1(n-1))$.  
Thus $[Q_0(n-1),Q_1(n-1)]$  lies 
to the left of the interior of $\Lambda_m^{\rm out}$ and there is no nontrivial intersection.

If $m=1$, the second inequality in (\ref{eq:Inequalities}) is automatically satisfied, 
since $P_1(0)_x=1$.  Assume therefore $m\geq2$, and
suppose  $P_1(m-1)_x \leq Q_1(n)_x$. Then $\tau(P_1(m-1))\geq \tau(Q_1(n))$. 
But $\tau(Q_1(n))=\tau(Q_0(n))$ and, using (\ref{eq:nPinequalities})
and monotonicity, $\tau(P_1(m-1))>\tau(P_0(m-1))$.  
Thus  $\tau(P_0(m-1))\leq\tau(Q_0(n))$.  
We conclude that the left-hand edge $[Q_0(n),Q_1(n)]$ of $\Lambda_n^{\rm in}$ lies 
to the right of the interior of $\Lambda_m^{\rm out}$ and there is no nontrivial 
intersection.

Defining $n_+(m)=\lceil n(P_0(m)) \rceil$, we know that $n_+(m)-1<n(P_0(m))$ and hence, by monotonicity, the first inequality of  (\ref{eq:Inequalities}) holds, not only for $n=n_+(m)$, but for all $n\leq n_+(m)$.   Similarly, if $n_-(m)=\lceil n(P_1(m-1) \rceil$, we know that $n_-(m)<n(P_1(m-1))$ and hence, by monotonicity, the second inequality of  (\ref{eq:Inequalities}) holds, not only for $n=n_-(m)$, but for all $n\geq n_-(m)$. We conclude that both inequalities hold for $n_-(m)\leq n \leq n_+(m)$.

From lemma \ref{lemma:nP}, we know that $\lambda(n(P_0(m))-n(P_1(m-1))$ is bounded below by a positive number.  Hence $n_+(m)-n_-(m)$, which by definition exceeds $n(P_0(m)-n(P_1(m-1))-1$, tends to infinity in the limit of vanishing $\lambda$.
\endproof



\subsection{Expansion of $Z(m,n)$}
The point $Z(m,n)$ lies at the intersection of the line segments $[P_0(m), P_1(m-1)]$ 
and $[Q_0(n-1), Q_1(n)]$. Using linear algebra, we have
\begin{equation}\label{eq:Zmnformula}
Z(m,n)=
\frac{1}{w\times z}
\left(\begin{array}{cc}
-z_x & w_x\\
-z_y & w_y
\end{array}\right)
\left(\begin{array}{c}
u\times w\\
v\times z
\end{array}\right)
\end{equation}
where
$$
u=P_0(m),\quad v=Q_0(n-1),\quad w=P_1(m-1)-P_0(m),\quad z=Q_1(n)-Q_0(n-1).
$$
and, for any pair of 2-vectors $a$ and $b$, the cross product is defined as
$$
a\times b= a_x b_y  -a_y b_x.
$$
Inserting the estimates of lemma \ref{lemma:PQexpand} and corollary
\ref{corollary:mtauexpand}, and once again applying the rules of interval 
arithmetic, we now obtain similar estimates for $Z(m,n)$.

\begin{lemma}\label{lemma:Zmnexpand}
For $0<\lambda\leq \lambda_0$ and $1\leq m\leq \lambda^{-\nu_0}$,
\begin{eqnarray}
Z(m,n)&=&\left( \frac{1}{2} (1+\tau) + \frac{\lambda}{8}\left(7+(1-4m)\tau^2\right)\right.
\label{eq:Zmnexpand}\\&&
\left.+\frac{\lambda^2}{32}\left(4-(5+16 m)\tau+(1-8 m+16 m^2)\tau^3\right) +r_{Zx} \lambda^{5/2}, \right.
\nonumber\\&&
\left.1+\frac{\lambda}{4}(1-(2 m+1)\tau) +\frac{\lambda^2}{16}\left(1+2 m+(8 m^2+2 m-1)\tau^2\right)\right.
\nonumber\\&&
\left.+\frac{\lambda^3}{192}\left(-20+48 m-16 m^2 +(23+6 m-80 m^2)\tau -(3-18 m+96 m^3)\tau^3\right) +
r_{Zy}\lambda^{7/2}\right),
\nonumber
\end{eqnarray}
where  $\tau=\tau(Q_0(n))$ and, for $m\geq 3$,
$$
r_{Zx} \in \left[-\frac{1033}{2000},\frac{2671}{5000}\right],\quad
r_{Zy} \in \left[-\frac{5563}{1000},\frac{7373}{1000}\right].
$$
For $m=1,2$, we have
\begin{eqnarray}
Z(1,n)&=&\left(\frac{1}{2}(1+\tau)+\frac{\lambda}{8}(7-3\tau^2)+\frac{\lambda^2}{32}(4-21\tau+9\tau^3)+r_{Z1x} \lambda^3,\right.\nonumber\\
&& \left.1+\frac{\lambda}{4}(1-3\tau) +\frac{3\lambda^2}{16}(1+3\tau^2) + \frac{\lambda^3}{64}(4-17\tau-27\tau^3)+r_{Z1y} \lambda^4 \right),\nonumber
\end{eqnarray}
with
$$
r_{Z1x} \in \left[-\frac{1011}{100},\frac{334}{25}\right],\quad
r_{Z1y} \in \left[-\frac{4211}{100},\frac{437}{10}\right].
$$
and
\begin{eqnarray}
Z(2,n)&=&\left(\frac{1}{2}(1+\tau)+\frac{7\lambda}{8}(1-\tau^2)+\frac{\lambda^2}{32}(4-37\tau+49\tau^3)+r_{Z2x}\lambda^3,\right.\nonumber\\
&& \left. 1+\frac{\lambda}{4}(1-5\tau)+\frac{5\lambda^2}{16}(1+7\tau^2)+\frac{\lambda^3}{64}(4-95\tau-245\tau^3)+r_{Z2y}\lambda^4 \right),\nonumber
\end{eqnarray}
with
$$
r_{Z2x} \in \left[-\frac{506}{25},\frac{331}{20}\right],\quad
r_{Z2y} \in \left[-\frac{9517}{100},150\right].
$$
\end{lemma}
The proof is given in appendix A.6.

\subsection{Disk areas }
In this section we complete the proof of theorem B. 
The regular fixed point $Z(m,n)$ of $L$ is a periodic point of
the map $F$, with period  $t(m,n)=4(m+n)-1$, see equation (\ref{eq:Periods}).
Let $A_{m,n}$ be the area of the cell (stability ellipse) of $Z(m,n)$;
furthermore, for fixed $m$, let $A_m$ be the total area of the ellipses
associated with orbits of $Z(m,n), \;n=n_-(m),\ldots,n_+(m)$. Thus
$$
A_m=\sum_{n=n_-(m)}^{n_+(m)}A_{m,n} t(m,n)\qquad t(m,n)=4(m+n)-1.
$$
Finally, we denote by $A^{\rm reg}$ the total area associated with the
regular orbits:
$$
A^{\rm reg}=\sum_{m=1}^{M} A_m.
$$

Our goal in this section is to prove that the $\lambda\rightarrow 0^+$ limits of 
$A_m$ and $A^{\rm reg}$ exist, and to derive explicit expressions for the limiting areas.

From theorem \ref{theorem:atomintersections}, we know that the ellipse associated with $Z(m,n)$
is tangent to one of the lines ${\bf 1}$ or ${\bf 5}^1$, and so the area of the ellipse is just
$\pi$ times the square of the $\cal Q$-metric distance from $Z(m,n)$ to the tangent point on the
closer line. Using the inner product ${\cal Q}$ defined in equation (\ref{eq:Metric}),
the distance between a point $Z$ and a line segment $[P,Q]$ is given by
$$
\sqrt{ {\cal Q}(U,U) - \frac{{\cal Q}(U,V)}{{\cal Q}(V,V)}},\qquad  U=Z-P,\hskip 5pt V=Q-P.
$$
For the area of disks tangent to the line ${\bf 1}$, we apply this formula
with $Z=Z(m,n), P=(0,1), Q=(1,1)$ to obtain the expression
\beq\label{eq:Amn0}
A^{(0)}_{m,n}=\pi\left(1-\frac{\lambda^2}{4}\right) (Z(m,n)_y -1)^2.
\eeq
For the disks tangent to ${\bf 5}^1$, on the other hand, we use
$P=(g,1),\, Q=(1,g'')$ (cf. equations (\ref{eq:ggprime}) and (\ref{eq:gdoubleprime})),
to obtain
(writing $Z=(Z_x,Z_y)$)
\beq\label{eq:Amn1}
A^{(1)}_{m,n}=\pi \left(1-\frac{\lambda^2}{4}\right) \left(1-Z_y + (1-2 Z_x) \lambda -(1-3 Z_y)\lambda^2-(1-Z_x) \lambda^3 -Z_y \lambda^4\right)^2.
\eeq
The above formulae give us a straightforward way of calculating
$$
A_{m,n}= \min\{ A^{(0)}_{m,n},A^{(1)}_{m,n}\}
$$
for given $m,n,\lambda$.  Supplemented  by Taylor expansion with rigorous bounds on the remainders,
the same formulae also provide the machinery for our asymptotic calculations in the
limit $\lambda\rightarrow 0^+$.

\subsection{Bounds for the disk areas}
To obtain formulae for the disk areas $A^{(0)}(m,n)$ and $A^{(1)}(m,n)$, we insert $Z(m,n)$ from lemma 
\ref{lemma:Zmnexpand} into (\ref{eq:Amn0}) and (\ref{eq:Amn1}), respectively, combining the remainder 
intervals using the rules of interval arithmetic \cite{ESupplement} (secs. E.10,E.10A).  We state the results as a lemma.
\begin{lemma}\label{lemma:Amnexpand}
For $0<\lambda\leq \lambda_0$ and $1\leq m\leq \lambda^{-\nu_0}$,
\begin{eqnarray}
A^{(0)}_{m,n}&=&\frac{\pi}{16}\left(1-(2 m+1)\tau\right)^2\lambda^2 +r_{A0}\, \lambda^{5/2},
\label{eq:Amn0expand}
\\
A^{(1)}_{m,n}&=&\frac{\pi}{16}\left(1-(2 m-3)\tau\right)^2\lambda^2 +r_{A1} \,\lambda^{5/2},
\label{eq:Amn1expand}
\end{eqnarray}
where $\tau=\tau(Q_0(n))$  and, for $m\geq 3$,
$$
r_{A0} \in\left[-\frac{6347}{100000},\frac{183}{6250}\right],\quad
r_{A1}\in \left[-\frac{7091}{100000},\frac{2857}{50000}\right],\quad
r_{A1}-r_{A0}\in \left[-\frac{859}{50000},\frac{3761}{100000}\right].
$$
For $m=1$, only $A^{(0)}(1,n)$ is relevant, with
$$
r_{A0} \in\left[-\frac{221}{6250},\frac{1179}{100000}\right],
$$
and, for $m=2$,
$$
r_{A0} \in\left[-\frac{393}{2000},\frac{787}{20000}\right],\quad
r_{A1}\in \left[-\frac{7859}{1000000},\frac{7859}{1000000}\right],\quad
r_{A1}-r_{A0}\in \left[-\frac{63}{2000},\frac{943}{5000}\right].
$$
\end{lemma}

\begin{figure}[t]
\begin{center}
\epsfig{file=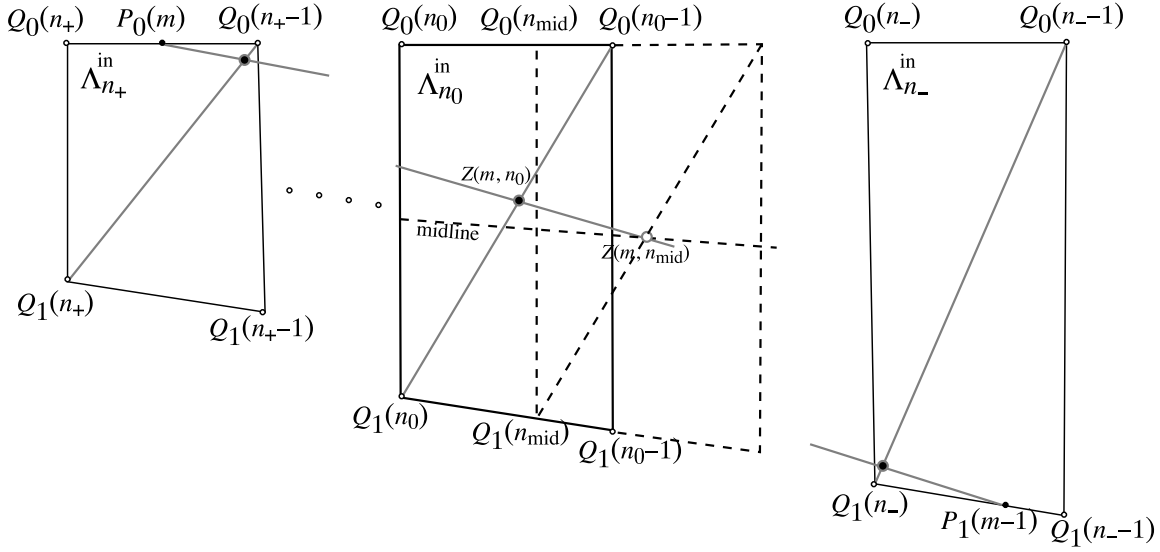,width=6in}
\end{center}
\caption{\label{fig:endpoints}
\small Sketches of the atoms $\Lambda_{n_+}^{\rm in}$, $\Lambda_{n_0}^{\rm in}$, 
and $\Lambda_{n_-}^{\rm in}$. 
The real number $n_{\rm mid}$ is defined by the condition that $Z(m,n_{\rm mid})$ 
lies on the `mid-line' equidistant (in the ${\cal Q}$-metric) from the upper 
and lower boundaries.}
\end{figure}

For fixed $m$, the fixed points $Z(m,n)$ are lined up along the line segment 
$[P_0(m),P_1(m-1)]$ in order of decreasing $n$, with $Z(m,n_+(m))$ just to the right 
of $P_0(m)$ and $Z(m,n_-(m))$ just to the left of $P_1(m-1)$ (see figure \ref{fig:endpoints}).  
Regarded as a continuous function of $n$, $A^{(0)}(m,n)$ increases from zero monotonically 
with decreasing $n$ while $A^{(1)}(m,n)$ decreases to zero monotonically. 
Thus there exists a value $n=n_{\rm mid}(m)$ such that
$$
A^{(0)}(m,n_{\rm mid}(m))=A^{(1)}(m,n_{\rm mid}(m)),
$$
and hence
$$
\min\{A^{(0)}(m,n),A^{(1)}(m,n)\}=\cases{
A^{(1)}(m,n) & $n_-(m)\leq n \leq n_0(m)-1$\cr
\noalign{\vskip 5pt}
A^{(0)}(m,n) & $n_0(m)\leq n \leq n_+(m)$\cr}
$$
where 
\begin{equation}\label{eq:n0}
n_0(m)=\lceil n_{\rm mid}(m)\rceil.
\end{equation}

To locate $\tau_{\rm mid}=\tan(\frac{\pi}{4}-2n_{\rm mid}\theta)$
(see equation (\ref{eq:etatau})), we equate $A^{(0)}(m,n)$ and $A^{(1)}(m,n)$ in 
\ref{eq:Amn0expand}) and (\ref{eq:Amn1expand}) and solve the resulting 
quadratic equation in $\tau$ to get (discarding the irrelevant root)
$$
 \tau_{\rm mid}=\frac{1+\sqrt{1+\frac{8}{\pi}(2m-1)\lambda^{1/2}(r_{A1}-r_{A0})}}{2 (2 m-1)}.
$$
Next we have
\begin{lemma}\label{lemma:taumid}
The quantity $A^{(0)}(m,n)-A^{(1)}(m,n)$, considered as a function of $\tau$ for fixed $m\in [2,M]$, vanishes at
$$
\tau_{\rm mid}=\frac{1}{2 m-1}+ r_{\rm mid} \lambda^{1/2},  
$$
where, for $m\geq 3$,
$$
r_{\rm mid}\in \left[-\frac{2189}{100000},\frac{599}{25000}\right],
$$
and, for $m=2$,
$$
r_{\rm mid}\in \left[-\frac{1003}{25000},\frac{601}{5000}\right].
$$
\end{lemma}
Lemma \ref{lemma:taumid} is easily proved using the inequalities 
(valid for $0\leq x,y,\leq 1$ and $-x\leq z\leq y$)
$$
 1-x\leq\sqrt{1+z}\leq 1+y/2.
$$
An immediate (and useful) consequence of this lemma are the bounds
\begin{eqnarray}
&&\tau_{\rm mid}<\frac{2}{5},\label{eq:taumidbound}\\
&&| (2m-1\pm 2)\tau_{\rm mid}-1 | <1,
\label{eq:taumidinequality}
\\
&&|A^{(0)}(m,n)| < \frac{1}{4}\lambda^2,   \quad  n_{\rm mid}\leq n \leq n_+,
\label{eq:A0bound}
\\
&&|A^{(1)}(m,n)| < \frac{1}{4}\lambda^2,   \quad  n_-\leq n \leq n_{\rm mid}.
\label{eq:A1bound}
\end{eqnarray}

\subsection{Total disk area}
We are now in a position to write an exact formula for the sum 
of disk areas within a single atom $\Lambda^{\rm out}_m$:
\begin{eqnarray}
A_1&=&\sum_{n=1}^{n_+(1)} A^{(0)}(1,n) t(1,n),\label{eq:A1exact}\\
A_m&=&\sum_{n=n_-(m)}^{n_0(m)-1} A^{(1)}(m,n) t(m,n) +
\sum_{n=n_0(m)}^{n_+(m)} A^{(0)}(m,n) t(m,n),\label{eq:Amexact}
\end{eqnarray}
where the `midpoint' $n_0(m)$ was defined in (\ref{eq:n0}). 
To calculate the leading behaviour for $\lambda$ tending to zero, we note that the terms 
in (\ref{eq:Amexact}) proportional to $4 m-1$ are negligible in the limit. Specifically, 
\begin{eqnarray}\label{eq:A1discard}
3 A_1&<&3 N A^{(0)}(1,1)<12 N \lambda^2<\frac{3}{25} \lambda^{1/2},
\\ \label{eq:Amdiscard}
(4m-1)A_m &< & 4 m N A^{(0)}(m, n_{\rm mid})<\frac{1}{10}\lambda^{1/2},\quad m\geq 2,
\end{eqnarray}
where we have applied (\ref{eq:A0bound}) as well as 
$$
\lambda<10^{-4}, \quad m\leq \lambda^{-1/4}, \quad \lambda N\leq \frac{\pi\lambda}{8\theta}<1.
$$

The next step toward the derivation of an asymptotic formula is to insert the 
estimates (\ref{eq:Amn0expand}) and (\ref{eq:Amn1expand}) of lemma \ref{lemma:Amnexpand},
to obtain
\begin{eqnarray}\label{eq:Amexpand}
A_1&=&\frac{\pi\lambda^2}{4}\left(\sum_{n=1}^{n_+(1)} \left(1-3 \tau\right)^2 \, n\right)+ r_A \lambda^{1/2},
\\ \nonumber
A_m&=&\frac{\pi\lambda^2}{4}\left(\sum_{n=n_-(m)}^{n_0(m)-1} \left(1-(2 m-3) \tau\right)^2 \, n +
\sum_{n=n_0(m)}^{n_+(m)} \left(1-(2 m+1)\tau\right)^2\, n\right) + r_A \lambda^{1/2}
\end{eqnarray}
where $\tau=\tau(Q_0(n))$ and 
$$
|r_A|<\frac{1}{5}.
$$
Next, we replace the sum over $n$ by an integral over $\eta$ using the Euler-Maclaurin formula:
\begin{eqnarray}\nonumber
A_1&=&\frac{\pi \lambda^2}{4}\int_1^{n_+}dn\,n(1-3\tau)^2 +\frac{\pi \lambda^2}{8}\left(4+(1-3\tan\eta_+)^2n_+(1)\right)+
(r_{EM0}+r_A)\lambda^{1/2},
\\ \nonumber
A_m&=&\frac{\pi\lambda^2}{4}\left(\int_{n_-}^{n_0-1} \left(1-(2 m-3)  \tau\right)^2 n dn\right)
+\frac{\pi\lambda^2}{4}\left(\int_{n_0}^{n_+} \left(1-(2 m+1)  \tau\right)^2 n dn\right)\\
&&\quad+\frac{\pi\lambda^2}{8}\left(
\left(1-(2 m-3) \tan\eta_-\right)^2 n_- +\left(1-(2 m-3) \tan(\eta_0-2\theta)\right)^2 (n_0-1)\right)\nonumber\\
&&\quad+\frac{\pi\lambda^2}{8}\left( \left(1-(2 m+1) \tan\eta_0\right)^2 n_0+\left(1-(2 m+1) \tan\eta_+\right)^2 n_+\right)\label{eq:AmEM}\\ \nonumber
&&\quad+(r_{EM0}+r_{EM1} + r_A )\lambda^{1/2},\nonumber
\end{eqnarray}
where, for $1\leq m\leq \lambda^{-1/4}$,
$$
\begin{array}{ll}
r_{EM0} &\leq\, \frac{\pi\lambda^{3/2}}{8}\int_{n_0}^{n_+} dn\left| \frac{d}{dn}\left(n \left(1-(2 m+1) 
    \tan(\frac{\pi}{4}-2 n \theta)\right)^2\right) \right|\\
&\leq\, \frac{\pi\lambda^{3/2}}{8}\left[\int_{n_0}^{n_+} dn\left(1-(2 m+1) \tau\right)^2 +4\theta (2m+1)\int_{n_0}^{n_+} dn\, n\left|1-(2 m+1) \tau\right| \sec^2\eta\right]\\
&\leq\,\frac{\pi\lambda^{3/2}}{8}(N+4\theta N^2(2m+1))<\frac{1}{7}
\end{array}
$$
and an analogous expression for $r_{EM1}$ for $m>1$.  

As a final simplification, we change the integration variable from $n$ to $\eta$, defined in
equation (\ref{eq:etatau}), changing the integration limits to the values 
$\eta=\cot^{-1}(2 m+1)$, $\cot^{-1}(2 m-1)$, and $\cot^{-1}(2 m-3)$.   
We state the result in the form of a lemma, which refers to the cut-off parameter 
values (\ref{eq:CutOffParameters}).
\begin{lemma}\label{lemma:Am}
For $0<\lambda\leq \lambda_0$, we have
\begin{eqnarray}
A_1&=&I(2,1,2)+ r \lambda^{1/2},\nonumber \\
A_m&=& I(m+1,m,m+1) + I(m,m-1,m-1) + r \lambda^{1/2},\qquad 2\leq m\leq \lambda^{-\nu_0},\nonumber
\end{eqnarray}
where
$$
I(a,b,c)\stackrel{\rm def}{=}\frac{\pi}{4}
  \int_{\cot^{-1}(2 a -1)}^{\cot^{-1}(2 b -1)}((2 c-1)\tan\eta-1)^2 (\pi/4-\eta) d \eta
$$
and $|r| < \frac{4}{5}$.
\end{lemma}
\proof
The contribution to $r$ of the remainder terms in (\ref{eq:AmEM}) is bounded in magnitude by $1/2$.  Moreover, the endpoint terms in (\ref{eq:AmEM}), for both $m=1$ and $m>1$, are bounded by $\frac{\pi}{2} N \lambda^2$ and so (using $N\lambda<1$) the magnitude of their contribution to $r$ is bounded by $1/50$.

The change of variables produces a factor 
\beq\label{eq:variablechange}
\frac{\lambda^2}{4\theta^2} =\frac{\sin^2(\theta)}{\theta^2}\in[1-\theta^2/3,1]\subset[1-10^{-4}\lambda^{1/2},1].
\eeq
Since each of the integrals over $\eta$ is bounded by 
$\frac{1}{2}\left(\frac{\pi}{4}\right)^3 <\frac{1}{2}$, we see that the total
contribution of the variable change to $r$ is bounded in magnitude by $10^{-4}$.

Finally, we can estimate the shifts in the integration limits using
$$
\begin{array}{rcl}
\left|\eta_+-\tan^{-1}\left(\frac{1}{2m+1}\right)\right|&\leq&
  \left|\tan^{-1}(\tau((P_0(m))-\tan^{-1}\left(\frac{1}{2m+1}\right)\right| + 2\theta,\\
\left|\eta_--\tan^{-1}\left(\frac{1}{2m-3}\right)\right|&\leq&
  \left|\tan^{-1}(\tau(P_1(m-1))-\tan^{-1}\left(\frac{1}{2m-3}\right)\right| + 2\theta,\\
\left|\eta_0-\tan^{-1}\left(\frac{1}{2m-1}\right)\right|&\leq&
  \left|\tan^{-1}(\tau_{\rm mid})-\tan^{-1}\left(\frac{1}{2m-1}\right)\right| + 2\theta,\\
\left|\eta_0+2\theta -\tan^{-1}\left(\frac{1}{2m-1}\right)\right|&\leq&
  \left|\tan^{-1}(\tau_{\rm mid})-\tan^{-1}\left(\frac{1}{2m-1}\right)\right| + 2\theta.
\end{array}
$$
Using lemmas \ref{lemma:taumid} and \ref{lemma:tauP} and the inequality (valid for $x>0$, $y>0$)
$$
\tan^{-1}(x+y)-\tan^{-1}(x) = \tan^{-1}\left(\frac{y}{1+x(x+y)}\right)<\tan^{-1}(y)<y,
$$
it is straightforward to calculate bounds for the right-hand expressions.
We find that he contribution to $r$ of the shifts in the integration limits is bounded by $\frac{1}{4}$.  Adding up the various contributions to $|r|$, we get the stated upper bound.
\endproof

Now we are ready to sum up the disk areas from the orbits of all regular fixed points, 
thus completing the proof of theorem B. 

\begin{theorem}\label{theorem:Area}
In the limit $\lambda\rightarrow 0^+$, the total area $A^{\rm reg}(\lambda)$ 
of the islands of the regular periodic orbits of $F$ converges to the positive 
quantity
$$
A^{\rm reg}(0)=I(2,1,2) + \sum_{m=2}^{\infty}\left(I(m+1,m,m+1) + I(m,m-1,m-1)\right).
$$
\end{theorem}
\proof
It suffices to show that
$$
A^{\rm reg}(\lambda)=I(2,1,2) + \sum_{m=2}^{\mu(\lambda)}\left(I(m+1,m,m+1) 
   + I(m,m-1,m-1)\right) + O( \lambda^{1/4}).
$$
where $I$ is given by lemma \ref{lemma:Am} and $\mu(\lambda)=\lfloor\lambda^{-1/4}\rfloor$.
The partial sum 
$$
 \sum_{m=2}^{\mu(\lambda)}A_m(\lambda)
$$
is already in the desired form, since the integrals are $\lambda$-independent and we have 
$\lfloor\lambda^{-1/4}\rfloor$ remainder terms, each bounded by the same constant 
times $\lambda^{1/2}$.   

The remaining terms,
$$
 \sum_{m=\mu(\lambda)+1}^{M(\lambda)}A_m(\lambda)
$$
are bounded by the (${\cal Q}$-metric) area of the triangle 
$$
[P, P_0(\mu(\lambda)), P_1(\mu(\lambda))],\qquad 
     P=\left(\frac{1+2\lambda-\lambda^2-\lambda^3}{2-\lambda^2},1\right),
$$
multiplied by the longest return time, $4(M(\lambda)+N(\lambda))-1$. 
Inserting our estimates for the atomic vertices, the triangular area is 
\begin{eqnarray*}
&&\frac{1}{2}\sqrt{1-\lambda^2/4}\, \left(P_0(\mu(\lambda))_x-P_x\right)\, 
\left(1-P_1(\mu(\lambda))_y\right)\\
&&=
\frac{1}{2} 
\left( 1 + O(\lambda^2)\right) 
\left(\frac{1}{2 (2 M(\lambda)+1)}+O(\lambda^{3/2})\right) \left(\frac{\lambda}{2 M(\lambda)-1}+O(\lambda^{5/2})\right)=O(\lambda^{3/2}),
\end{eqnarray*}
while the upper bound on the period is $O(\lambda^{-1})$.  
Thus the upper bound on the sum of terms with $m>\mu(\lambda)$ is  $O(\lambda^{1/2})$.
\endproof

\subsection{Numerical evaluation of area formula}
Numerical integration gives the following results:
$$
\begin{array}{|c|l|}
m&\qquad\mbox{Area }A_m\\
\hline
1& 0.04394252102495575454\\
2& 0.027915747684440153071\\ 
3& 0.0043390573122902285760\\
4& 0.0011842200753144612564\\
5& 0.00044544206984605774866\\  
6& 0.00020336198417268636974\\
7& 0.00010566669827864850788\\
8& 0.0000602301047692083869367
\end{array}
$$
For asymptotically large $m$, the integrands can be expanded in powers of $\eta$ and integrated term by term to yield
$$
A_m=\frac{\pi^2}{48 m^4} +\frac{\pi(\pi-1)}{24 m^5} + O(m^{-6})
$$
Summing the numerical integrals over the first 2000 values of $m$ gives
$$
\sum_{m=1}^{2000}  A_m= 0.0783220277996.
$$
This is approximately 36 percent of the area in the square outside the fixed-point disk.

\section{Covering}\label{section:Covering}

In the previous sections we defined the return map $L$ of the domain $\Lambda$,
and then we identified a prominent sub-domain of $\Lambda$, 
namely the union of the regular atoms of $L$.
In the following theorem we show that the set of regular atoms is an appropriate 
surface of section for the orbits of $F$ outside the main island ${\cal E}$.
In so doing, we will complete the proof of theorem A in the introduction.

\begin{theorem}\label{theorem:Covering}
The images under $F$ of the regular atoms of the return map $L$ cover 
all of $\Omega\setminus {\cal E}$, apart from a set of area $O(\lambda^{2})$.
\end{theorem}

\proof  All metric considerations refer to the ${\cal Q}$-metric (\ref{eq:Metric}). 
We begin to deal with the area of the region complementary to ${\cal E}$.
The ${\cal Q}$-area of the square is
$$
{\cal A}^{\Omega}= \sqrt{1-\frac{\lambda^2}{4}} = 1 + O(\lambda^2).
$$
The area of the ellipse ${\cal E}$ is $\pi r^2$, where $r$ is the ${\cal Q}$-distance 
between the fixed point $(1/(2-\lambda), 1/(2-\lambda))$ and the point of tangency 
$T_0=((1+\lambda)/2,1)$.  This is
$$
{\cal A}^{\cal E}=\frac{2-3 \lambda + \lambda^3}{8-4 \lambda}= \frac{\pi}{4} (1-\lambda) + O(\lambda^2).
$$
Thus, the area outside ${\cal E}$ is
$$
{\cal A}^{\Omega\setminus{\cal E}}= 1 - \frac{\pi}{4} +\frac{\pi}{4} \lambda + O(\lambda^2).
$$

Now we consider the dynamics, reverting to the representation (\ref{eq:Lbar}). 
First, we calculate the area covered by the images of the regular atoms 
$\Lambda_m^{\rm out}$ under the outer map $\bar L^{\rm out}$. 
These atoms are polygons whose vertices are listed in appendix A.1. 

In the first-order area calculation, we can expand formulae (\ref{eq:mequal1}) and 
(\ref{eq:mequal2}) in Taylor series about $\lambda=0$, keeping only the lowest two terms and 
an $O(\lambda^2)$ remainder. 
For $m=1$, the quadrilateral $\Lambda_1^{\rm out}$ reduces to a right triangle 
$[(\lambda+{2}/{3},1),(1,1),(1,1-\lambda)]$ with area ${\lambda}/{6} +O(\lambda^2)$.  
Since the transit time under $F$ is 7 (theorem \ref{theorem:Lout}) we get a contribution to 
the overall area sum of 
$$
\frac{7}{6}\lambda+O(\lambda^2).
$$

For $m=2$, the hexagon $\Lambda_2^{\rm out}$ reduces to a quadrilateral with vertices 
$\{(\frac{3}{5}+\lambda,1)$, $(\frac{2}{3}+\lambda,1)$, $(1,1-\lambda)$, 
$(\frac{2}{3}+\lambda,1-\frac{\lambda}{3})\}$, and hence area $\lambda/15$.
Multiplying by the transit time of 11, we get
$$
\frac{11}{15}\lambda+O(\lambda^2).
$$

The remaining domains $\Lambda_m^{\rm out}$ are the quadrilaterals 
$[P_0(m),P_0(m-1),P_1(m-1),P_1(m)]$, $L^{\rm out}$-symmetric about the respective 
axes $[P_0(m),P_1(m-1)]$.  We use the formulae $(\ref{eq:P0}$) and (\ref{eq:P1}) to calculate the area of $\Lambda_m^{\rm out}$ as 
\begin{eqnarray}
{\cal A}_m^{\rm out}&=&
\sqrt{1-\lambda^2/4}\, \left(P_0(m)-P_1(m)\right)\times \left(P_1(m)-P_1(m-1)\right) \nonumber\\
&=&\frac{\lambda^4\cos((2m-1)\theta)}{8 \prod_{k=-1}^1\sin((2m+(2k-1))\theta)}.\nonumber
\end{eqnarray}
Expanding the right-hand side in powers of $\theta$, we have
$$
{\cal A}_m^{\rm out}=\frac{\lambda^4\left(1+m^2\theta^2 r(\theta,m)\right)}{8\theta^3 (2m-3)(2m-1)(2m+1)},
$$
where, as we prove in appendix A.8, $r(\theta,m)$ is uniformly bounded above and below for $3\leq m\leq M,\;0<\theta\leq \sin^{-1}(\lambda/2)$.
Multiplying by the transit time $4 m+3$ and summing over $m$ gives
$$
\sum_{m=3}^M \frac{4m+3}{(2m-3)(2m-1)(2m+1)}=\frac{7}{20} -\frac{\lambda}{\pi} +O(\lambda^2),
$$
$$
\sum_{m=3}^M \frac{m^2(4m+3)}{(2m-3)(2m-1)(2m+1)}=\frac{\pi}{4\lambda} +O(\log \lambda),
$$
and hence
$$
\sum_{m=3}^{M(\lambda)} (4m+3){\cal A}_m^{\rm out} = \frac{7}{20}\lambda+O(\lambda^2).
$$

Combining the results for all $m$, the area contribution is
$$
{\cal A}^{\rm out}=\frac{9}{4}\lambda +O(\lambda^2).
$$

Second, we compute the area covered by the orbits of the regular atoms in 
$\Lambda^{\rm in}$, under the inner map $\bar L^{\rm in}$.
Here again we start with exact formulae for the vertices of the regular atoms 
given in appendix A.1, and the expansions of lemma \ref{lemma:PQexpand}.
Each $\Lambda^{\rm in}_n$ is a quadrilateral $[Q_0(n),Q_0(n-1),Q_1(n-1),Q_1(n)]$ 
with $L^{\rm out}$ symmetry axis $[Q_0(n-1),Q_1(n)]$ and the area is calculated 
to be \cite{ESupplement}
\begin{eqnarray*}
{\cal A}_n^{\rm in}&=&
\sqrt{1-\lambda^2/4} \,(Q_0(n-1)-Q_0(n))\times (Q_0(n)-Q_1(n))
\\
&=&
\frac{1}{2}\tau(1+\tau^2)\lambda^2 -\frac{1}{8}(1+4\tau+4\tau^2+4 \tau^3+3\tau^4) +r_A \lambda^4,
\end{eqnarray*}
where $\eta$ and $\tau$ were defined in (\ref{eq:etatau}), and $r_A$ is a 
function of $\lambda$ and $\tau$ which satisfies $r_A<29$.
Next we sum the atom areas multiplied by the transit times $4(n-1)$ 
(theorem \ref{theorem:Lin}), and use the Euler-Maclaurin formula to turn 
it into an integral:
\begin{eqnarray*}
{\cal A}^{\rm in}&=&\sum_{n=1}^{N(\lambda)} 4(n-1){\cal A}_n^{\rm in}
\\
&=&\frac{1}{2} \int_0^{\frac{\pi}{4}} d\eta\,\left(\frac{\pi}{4}-\lambda-\eta\right)
   \left(4 (\tau+\tau^3)-(1+4\tau+4\tau^2+4 \tau^3+3\tau^4)\lambda\right)+O(\lambda^2)
\\
&=& 1-\frac{\pi}{4}+\frac{\pi-9}{4}\,\lambda+O(\lambda^2).
\end{eqnarray*}

Here we omit the detailed handling of error terms (Euler-Maclaurin endpoint and remainder terms, 
change of integration variable, adjustment of integration limits), which parallels that 
encountered in the calculation of disk areas in section \ref{section:ReturnMap}.
Details may be found in \cite{ESupplement}.

Now the regular atoms of $L$ are the intersections of the `in' and `out' atoms ---see 
equation (\ref{eq:AtomsL}). The mismatch between the union of the two sets
of atoms is $O(\lambda^3)$, from theorems \ref{theorem:Lin} and \ref{theorem:Lout},
while the respective transit times are $O(\lambda^{-1})$. Hence the area covered by the
images of the regular atoms of $L$ is ${\cal A}^{\rm in}+{\cal A}^{\rm out}+O(\lambda^2)$.
Combining the above area expressions gives us the desired identity
$$
{\cal A}^{\Omega}={\cal A}^{\cal E} + {\cal A}^{\rm out} + {\cal A}^{\rm in} + O(\lambda^2).
$$
\endproof

Although the irregular atoms of $L$ do contribute to the total area, 
the proof of theorem \ref{theorem:Covering} avoids the great difficulty 
of deriving bounds for their return times.

\section{Extension of our results}
\label{section:Extensions}

We consider, briefly and informally, some extensions of our results.

The regularity of the limit $\lambda\to0^+$ depends on the very simple 
dynamics at $\lambda=0$, for which there is a single atom. There is 
one dominant fixed point, while all other cycles disappear at $\lambda=0$,
and are arranged in a hierarchy that affords a perturbative analysis
in the small parameter $\theta$. Some asymptotic features depend on 
the product of $\theta$ and the return time of relevant domains,
so that large return times lead to an unavoidable non-uniformity in the
convergence to the limit. However, this problem remains tractable, due 
to the regular arrangement of the atoms.

Similar features are present in the limit $\lambda\to 0^-$, so we expect similar 
phenomena. For $-1<\lambda <0$, there are three symbols: $\iota(x,y)\in \{0,1,2\}$.
However, the atom $\Omega_0$ is a single point (the origin) and can be neglected.
The atom $\Omega_1$ is again a quadrilateral with area $O(1)$, while 
$\Omega_2$ is a triangle with area $O(\lambda)$, which plays the same 
role as the atom $\Omega_0$ in the positive $\lambda$ case.

For $-1 < \lambda < 0$, the basic cycles are
\begin{equation}\label{eq:BasicCyclesII}
\vcenter{\baselineskip 18pt\halign{
\quad 
 $#$ \hfil&\qquad
 $#$ \hfil&\qquad
 $\displaystyle #$ \hfil &\hskip 20pt 
 $\displaystyle #$ \hfil &\hskip 20pt 
 $\displaystyle #$ \hfil \quad
\cr
t &\mbox{code}  & \mbox{denominator} & \mbox{numerators}&\mbox{radius} \cr
\noalign{\vskip 4pt\hrule\vskip 5pt}
1&(\overline{1})  & 2-\lambda&\quad (1)&x_0(\lambda)\cr
2&(\overline{1,2}) & 4-\lambda^2&\quad (4+\lambda,2(1+\lambda))&1-x_0(\lambda)\cr
}}
\end{equation}
The fixed point is determined by the same rational function as for the 
$\lambda\geq 0$ case. The 2-cycle is different; however, changing the sign 
of $\lambda$, together with a reflection with respect to the centre of 
$\Omega$, reproduces the same behaviour to first order in $\lambda$. 
The symmetry properties of these cycles are also unchanged.

\begin{figure}[t]
\begin{center}
\epsfig{file=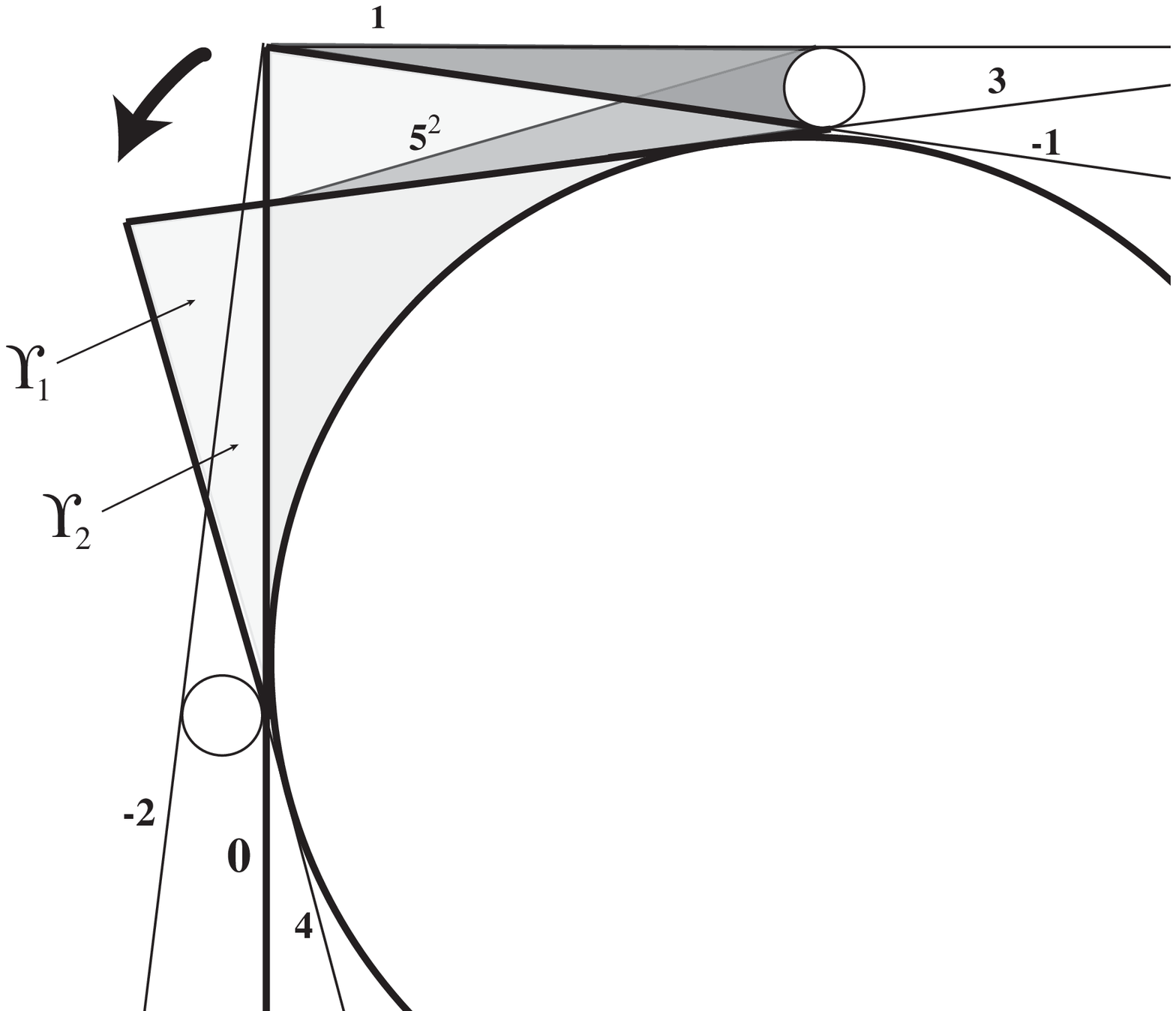,width=10cm}
\end{center}
\caption{\label{fig:TwoSectorsNeg}
\small The linked sector maps of the 1-cycle and the 2-cycle for a small 
negative value of the parameter $\lambda$. The corresponding sectors
$\Sigma$ and $\Sigma'$ are represented in light and dark grey, respectively.
Under the action of the map $F^4$, both sectors turns counterclockwise.
The points in $\Sigma$ that do not end up in $\Sigma'$ (under $F$)
comprise the triangle $\Upsilon_1$; those that do comprise the
quadrilateral $\Upsilon_2$.
}
\end{figure}

The linked sector maps of the 1- and 2-cycles for $\lambda<0$ are shown in figure 
\ref{fig:TwoSectorsNeg}. Compared to the case $\lambda>0$, there is now
a tighter connection between the two sectors, which results in 
shorter transit times.
As for $\lambda>0$, one constructs the transit maps $L^{\rm in, out}$, 
and identifies a two-parameter family of regular fixed points with periods
$$
t(m,n)=4(m+n)-3
$$
and symbolic codes
$$
\iota^{(n,m)}=\bigl(\overline{1^{4n-1},(2,1)^{2m-1}}\bigr)\qquad n\geq 1,\quad m\geq 1,
$$
to be compared with (\ref{eq:Periods}) and (\ref{eq:Codes}), respectively.
As before, there are anomalous cycles of total area $O(\lambda)$, which interact
with the above structure in regions of phase space of area $O(\lambda^2)$.
The main cycle of this kind is a fixed point located in the North-West corner
of the square $\Omega$. 

Returning to the limit $\lambda\to0^+$, we recall that the mismatch 
between the sector maps of the 1 and 2-cycles generates a collection of small 
irregular domains, which were neglected in our perturbative analysis. 
The dynamics over these domains is determined by the interaction among 
certain periodic orbits lying outside the sector $\Sigma$. 

As an example, let us consider the turnstile associated to the domain 
$\Upsilon_1$ (figure \ref{fig:Turnstiles}). 
Under the action of $G\circ F^{-4}$, the left triangle comprising the turnstile is 
mapped to the triangle $\Phi=\langle {\bf 0},{\bf 7},{\bf 8}\rangle$, of area
$O(\lambda^3)$, lying between the atom $\Lambda_1^{\rm out}$ and the boundary 
${\bf 0}$ of $\Lambda$ (cf.~the first equation in (\ref{eq:AtomsSidesII})).

\begin{figure}[h]
\begin{center}
\epsfig{file=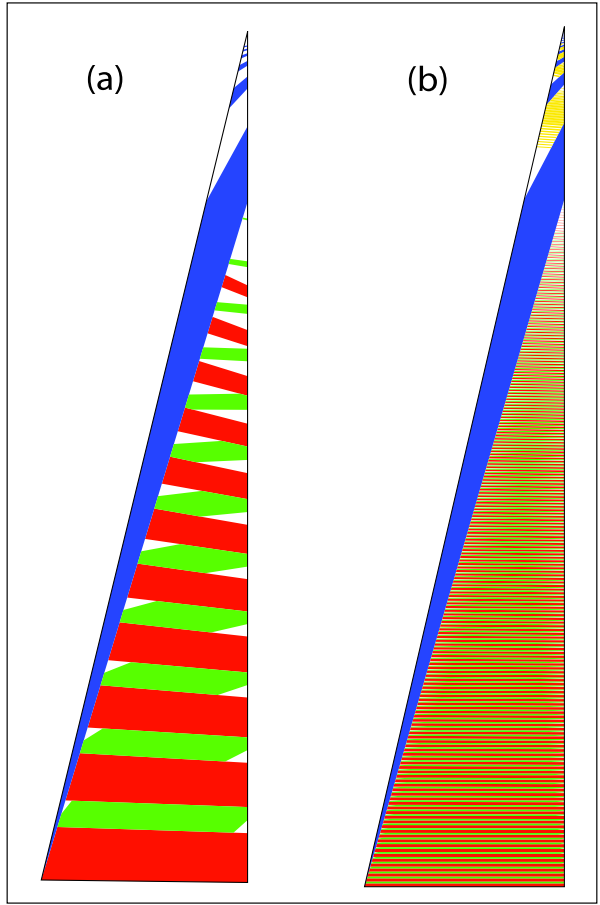,height=12cm}
\end{center}
\caption{\label{fig:IrregularDomain1}
\small (a) The irregular domain $\Phi$ on the eastern boundary 
of $\Lambda^{\rm out}_1$, showing three families of atoms. Here $\lambda=\frac{1}{64}$.  
(b) Same as (a), but with $\lambda=\frac{1}{1000}$.  
In both cases the horizontal and vertical dimensions have been rescaled 
to make the structure of the partition visible. } \end{figure}
The domain $\Phi$ belongs to the first regular atom of $L^{\rm in}$.
The partition of $\Phi$ into atoms under $\bar{L}^{\rm out}$ is shown 
in figure \ref{fig:IrregularDomain1} for two parameter values approaching
zero from above. The low-transit-time atoms (all symmetric) belong to three families, 
coloured blue, red, and green in the figure, and the picture suggests 
convergence to an asymptotic regime, although the convergence seems
slow. The atoms' respective codes and transit times are:
$$
\begin{array}{|c|c|c|}
\hline
\mbox{Colour} & \mbox{Code} & \mbox{Transit Time}\\ \hline\hline
\mbox{Blue} & \left(1 (1101101010)^{2m-1} 110111\right) & 20 m - 3\\ \hline
\mbox{Red} & \left(1 (110)^{4 m} 111\right) & 12 m+4 \\ \hline
\mbox{Green} &\left(1 (110)^{4 m-1} (1101101010) (110)^{4 m} 111\right) & 24 m+11\\
\hline
\end{array}
$$
We see that these codes shadow two periodic codes $(\overline{110})$ and 
$(\overline{1101101010})$. These correspond to a 3- and a 10-cycle lying 
outside $\Sigma$, whose cells have area $O(\lambda^2)$ and $O(\lambda^4)$, respectively.

In addition to $\lambda=0$, there is trivial dynamics (a single atom)
also for $\lambda=\pm1$ (rational rotation number $\rho=1/6,1/3$), and
for these values we expect the near-rational dynamics to be broadly similar 
to that considered in this paper. 
By contrast, the dynamics for all other rational rotation numbers is 
highly non-trivial. The simplest cases correspond to quadratic irrational 
$\lambda$-values (eight in all), for which the phase space is tiled by 
infinitely many periodic cells which admit an exact renormalization 
\cite{KouptsovLowensteinVivaldi}.
Preliminary investigations for $\lambda=(\sqrt{5}-1)/2$ ($\rho=1/5$) indicate 
that the construct of linked sector maps can only serve as a local model for 
near-rational dynamics; the global properties result from complex 
interactions between different local models.
We have observed a highly non-uniform convergence to the asympotic regime,
due to the presence of unbounded return times in the limiting rational
dynamics. This problem requires further investigation.

\section*{Appendix}
\subsection*{A.1 \quad  Rules of interval arithmetic with rational rounding}
Suppose we have rational lower and upper bounds on a set of real quantities 
$r_1,r_2,\ldots,r_n$ and we wish to obtain a similar bound on a rational 
function $R$ of these quantities. We can do so by assigning to the quantity
$r_i$ the closed interval $[a_i,b_i]$ where $a_i\leq r_i\leq b_i$ 
are the given bounds. If $r$ is a rational number, then we assign to it
the closed interval $[r,r]$.
We then proceed to deduce a bounding interval for $R$ by performing, in 
any convenient order, a sequence of elementary arithmetic operations, 
taken from the following list:
\begin{enumerate}
\item $ [a_1,b_1]+[a_2,b_2]=[a_1+a_2,b_1+b_2],$
\item $ [a_1,b_1] \times [a_2,b_2] = [\min\{a_1 a_2,\, a_1 b_2,\, b_1 a_2, \,b_1 b_2\},\max\{a_1 a_2,\, a_1 b_2, \,b_1 a_2, \,b_1 b_2\}],$
\item $ \mbox{For }ab>0,\quad [a,b]^{-1}=[\min\{\frac{1}{a},\frac{1}{b}\},\max\{\frac{1}{a},\frac{1}{b}\}],$
\item $ \mbox{For } a\geq 0, b>0, n\geq 0 \mbox{ or }a< 0, b\geq 0, n\geq 0,\quad
[a,b]^n = [\min\{a^n, b^n\},\max\{a^n,b^n\}],$
\item $\mbox{For } a>0,b>0,n>0,\quad [a,b]^{-n} = ([a,b]^n)^{-1}.$
\end{enumerate}
The result is an exact, but obviously non-unique, pair of bounds for $R$.

Since multiple application of these rules can lead to prohibitively large denominators 
for the resulting rational bounds, it is convenient to round off after every elementary 
operation.  
The lower and upper bounds are rounded down (resp. up) to the nearest rational number 
with a power of 10 in the denominator and a specified number of significant digits 
in the numerator.  Any common factors in numerator and denominator can be divided out.  
For example, with the 4-digit rational rounding used in this article, the interval 
$[\pi/2,\pi]$ would be rounded to $[157/100,3142/10000]$.

The Mathematica interval functions used to obtain rigorous bounds throughout this article are listed in sec. E.1 of \cite{ESupplement}.

\subsection*{A.2 \quad Proof of $P_k(m)$ estimates in lemma \ref{lemma:PQexpand}}

The initial terms in the stated formulae for $P_0(m)$ and $P_1(m)$ result from
the Taylor expansion of (\ref{eq:P0}) and (\ref{eq:P1}) up to second order in
the $x$-component and third order in the $y$-component (in the latter case,
the zeroth order term is trivially equal to unity).  Our task is to establish the stated numerical bounds
on the remainder coefficients $r_i,\, i=0x, 1x, 1y, 01x, 10y,11y$.

We begin with 
\beq\label{eq:r0xdef}
r_{0x}=(2m+1) P_0(m)_x-1-m-(2m+1)\lambda+\frac{m}{6}(m+1)\lambda^2.
\eeq
Defining auxiliary variables $u$ and $v$ by the expansion formulae
$$
x\cot(x)= 1-\frac{x^2}{3}-(1+u)\frac{x^4}{45},
$$
$$
\theta=\sin^{-1}\left(\frac{\lambda}{2}\right)=\frac{\lambda}{2}+(1+v)\frac{\lambda^3}{48}.
$$
Inserting these expansions reduces $P_0(m)_x$, hence $r_{0x}$, to ratios of 
polynomials in $\lambda, m, u$, and $v$.  The numbers of terms is in these 
polynomials is enormous, but the expressions are easily manipulated with 
computer assistance.  Moreover, all but a handful of terms have sufficiently 
high powers of $\lambda$ that they contribute insignificantly to the final 
estimates (recall our cut-off condition $\lambda<\lambda_0=10^{-4}$). 
The detailed calculations have been relegated to the electronic supplement \cite{ESupplement} (E.7.1). 

To obtain a bounding interval for $r_{0x}$, we use the following bounds \cite{ESupplement} (E.7.1) for our variables:
$$
\lambda\in[0,10^{-4}],\qquad m\in [2, \lambda^{-1/4}],\qquad u\in[-3/2500,3/1000], \qquad  v\in [0,2\times 10^{-8}].
$$
Each monomial of the form $\pm\lambda^a m^b u^d v^f$ 
is assigned a single bounding interval using the rules of interval arithmetic (see appendix A.1) 
with fixed precision rational rounding down (resp.~up) of lower (resp.~upper) bounds.  
The monomial bounds are then be rigorously combined, again using interval arithmetic, 
to give the final bounds for $r_{0x}$  stated in the lemma.  The same technique was 
used to establish \cite{ESupplement} (E.7.1-2) the stated bounds on the other remainder factors in the lemma.

We have implemented the process described above with the aid of Mathematica functions involving 
only elementary logical and arithmetic manipulations of symbols and integers. The specific 
procedures are listed in the Electronic Supplement \cite{ESupplement}, sec.E.1. Note that while our 
estimates are rigorous, they are not optimal, and can easily be changed by permuting the 
order in which monomial bounds are combined.  

A potential difficulty of the method comes from the fact that the upper bound for $m$ grows as 
$\lambda^{-1/4}$, which is numerically unbounded for $\lambda$ in the chosen range.  The problem 
is avoided by a simple trick based on the fact that all polynomial expressions which we wish to 
bound have only monomials of the form 
$\pm\lambda^a m^b \cdots$ where $b$ is either negative or no greater than $4 a$. In the first case 
we use $m^b\in [0,2^b]$, and in the second case we use $\lambda^a m^b\in[0,\lambda_0^{a-b/4}]$.
If the original rational function is not in this form, we try, if possible, to make it 
acceptable by dividing both numerator and denominator by a common power of $m$.

\subsection*{A.3\quad Proof of $Q_k(n)$ estimates in lemma \ref{lemma:PQexpand}}
Taylor expansion of (\ref{eq:Q0}) and (\ref{eq:Q1}) has been carried out with respect to $\lambda$ 
with $\tau$ held fixed.  The remainders, expressed as rational functions of 
$\lambda$ and $\tau$ are
\beq\label{eq:Qremainders}
r_{2x}=\frac{{\cal N}_{2x}}{{\cal D}_{2x}},\quad
r_{3x}=\frac{{\cal N}_{3x}}{{\cal D}_{3}},\quad r_{3y}=\frac{{\cal N}_{3y}}{{\cal D}_{3}},
\eeq

For $r_{2x}$, we have
$$
\begin{array}{ll}
{\cal N}_{2x}=&-64 + 64 \tau^2 + 32 \lambda + 32 \tau \lambda - 32 \tau^2 \lambda - 32 \tau^3 \lambda + 32 \lambda^2 - 
 16 \tau \lambda^2 - 40 \tau^2 \lambda^2 + 
 \\ 
 & 16 \tau^3 \lambda^2 +8 \tau^4 \lambda^2 - 24 \tau \lambda^3 + 
 16 \tau^2 \lambda^3 + 8 \tau^3 \lambda^3 - 4 \tau \lambda^4 + 10 \tau^2 \lambda^4 - 4 \tau^3 \lambda^4 - 
 2 \tau^4 \lambda^4 + 
 \\
 &\sqrt{4-\lambda^2} (32 - 32 \tau^2 - 16 \lambda - 16 \tau \lambda + 16 \tau^2 \lambda + 16 \tau^3 \lambda - 12 \lambda^2 + 
    8 \tau \lambda^2 +
\\    
&  16 \tau^2 \lambda^2 - 8 \tau^3 \lambda^2 - 4 \tau^4 \lambda^2 - 2 \lambda^3 + 5 \tau \lambda^3 - 
    6 \tau^3 \lambda^3 + 2 \tau^4 \lambda^3 + \tau^5 \lambda^3),
\\
\\
{\cal D}_{2x}=& 32\lambda^2\left(-8\tau+2\lambda^2\tau 
  +\sqrt{4-\lambda^2}(\lambda-\lambda\tau^2) \right).
\end{array}
$$

Application of interval bounds to these expressions runs into a serious difficulty: 
the denominator has a zero (barely!) within the allowed range $0\leq\tau\leq1$, namely at 
$$
\tau_1=\frac{2\sqrt{4-\lambda^2}-4-\lambda^2}{\lambda\sqrt{4-\lambda^2}}=\frac{\lambda}{4}+O(\lambda^3).
$$
Fortunately, the numerator shares the same zero, and we can divide both 
${\cal N}_{2x}$ and ${\cal D}_{2x}$ by the factor $(\tau-\tau_1)$. 
The result takes the form
\beq\label{eq:newr2x}
r_{2x}=\frac{ {\cal N}'_{2x} }{ {\cal D}'_{2x} }
 \left(2\sqrt{4-\lambda^2}-4+\lambda^2\right)^{-5},
\eeq
where $ {\cal D}'_{2x} $ and (especially) ${\cal N}'_{2x}$ are lengthy polynomials 
in $\lambda$, $\tau$, and $\sqrt{4-\lambda^2}$.
 
To apply the interval methods used successfully in appendix A.2, it turns out that 
we need to use a high-order estimate for the square root, namely
$$
\sqrt{4-\lambda^2}=2-\frac{\lambda^2}{4}-\frac{\lambda^4}{64}-\frac{\lambda^6}{512}-\frac{5 \lambda^8}{16384}-\frac{7\lambda^{10}}{131072}-\frac{21 v \lambda^{12}}{2097152}.
$$
Inserting the interval estimates
$$
v\in [1,\frac{1117}{1000}],\qquad
\tau\in[0,1],\qquad
\lambda\in[0,\lambda_0],
$$
in each of the 3 polynomial components of $r_{2x}$ in (\ref{eq:newr2x}) and applying the 
rules of appendix A.1 to combine them, we arrive, finally at the stated bounds.
Applying the same techniques to $r_{3x},\,r_{3y},\,r_{4x},\,r_{5x}$, we get the 
remaining estimates in lemma \ref{lemma:PQexpand}.  
Details of the calculations may be found in \cite{ESupplement}, sec. E.7.3.

\subsection*{A.4\quad $\tau$ estimates}
In addition to the contents of lemma \ref{lemma:PQexpand} concerning the cartesian 
coordinates of the $\Lambda_m^{\rm out}$ vertices, we will also need information 
relating to their $\tau$-coordinates (see (\ref{eq:etatau})).  
Estimates for the latter may be obtained by combining the formulae of lemma 
\ref{lemma:PQexpand}. In the following lemma, the remainder terms are 
rigorously bounded \cite{ESupplement} using the same techniques discussed in appendices A.2 and A.3.
\begin{lemma}\label{lemma:tauP}
For $\lambda\leq \lambda_0= 10^{-4}$ and $2\leq m \leq \lambda^{-1/4}$, we have the expansions
\begin{eqnarray}
\tau(P_0(m)) &=& \mu^{-1}+\frac{\lambda}{4}\left(5+2\mu^{-1}+\mu^{-2}\right)\nonumber\\
&& +\frac{\lambda^2}{48}\left(24+(57-16 m(m+1))\mu^{-1}+12\mu^{-2}+3\mu^{-3}\right) + r_{0\tau}\lambda^3,
\label{eq:tauP0m}
\end{eqnarray}
where 
$
\mu=2m+1,\; r_{0\tau}\in \left[-\frac{4541}{1000},\frac{3087}{1000}\right],
$
and
\begin{eqnarray}
\tau(P_1(m))&=& \mu'^{-1} +\frac{\lambda}{4}\left(1+2\mu'^{-1}+5\mu'^{-2}\right)\nonumber\\
&&+\frac{\lambda^2}{48}\left(24+(129-16(m+2)(m+3))\mu'^{-1}+60\mu'^{-2} +75\mu'^{-3}\right)+ r_{1\tau}\lambda^3,
\label{eq:tauP1m}
\end{eqnarray}
where 
$
 \mu'=2m-1,\; r_{1\tau}\in \left[-\frac{633}{20},\frac{731}{50}\right].
$
Moreover,
$$
\tau(P_0(1))=\frac{1}{3}+\frac{13}{9}\lambda+\frac{19}{27}\lambda^2+r_{01\tau}\lambda^3,  
$$
with $r_{01\tau}\in [\frac{383}{2500},\frac{2819}{1000}]$.
\end{lemma}
As a corollary, we get bounds on $\tau(Q_0(n))$ restricted to a single value of $m$.
\begin{corollary}\label{corollary:mtauexpand}
For $ \lambda\in [0, 10^{-4}], m\in [3,\lambda^{-1/4}]$,
$$  
m\tau(Q_0(n)) = r_\tau \in \left[\frac{2}{5},\frac{1001}{1000}\right].
$$ 
\end{corollary}

\proof
For those atoms $\Lambda^{\rm in}_n$ which intersect $\Lambda^{\rm out}_m$, we have 
the inequalities
\beq\label{eq:tauInequalities}
\tau(P_1(m-1))\geq \tau(Q_0(n))=\tan(\gamma(Q_0(n-1))-2\theta)\geq
\frac{\tau(P_0(m))-\tan(2\theta)}{1+\tan(2\theta)}.
\eeq
where, for $\lambda<10^{-4}$,
$$
\tan (2\theta)=\frac{\lambda\sqrt{4-\lambda^2}}{2-\lambda^2}<(1+10^{-8})\lambda.
$$
From lemma \ref{lemma:tauP} it is not difficult to derive\cite{ESupplement}
$$
\tau(P_0(m))=\frac{1}{2m+1}+\left(\frac{5}{4}+r_0\right)\lambda,\qquad r_0\geq -\frac{1467}{2\times 10^6},
$$
$$
\tau(P_1(m))=\frac{1}{2m-1}+\left(\frac{1}{4}+r_2\right)\lambda,\qquad r_2<\frac{153}{500}.
$$
With the inequalities (valid in the cut-off regime)
$$
\frac{m}{2m+1}\geq\frac{3}{7},\qquad
\frac{m}{2 m-3}\leq 1,\qquad 
0\leq\lambda m\leq \frac{1}{1000}
$$
we obtain
$$
\frac{2}{5} < m\tau(Q_0(n))< \frac{1001}{1000}.
$$
\endproof

\subsection*{A.5\quad Proof of lemma \ref{lemma:nP}}
The differences in the $n$-coordinate are related to those in $\tau$ by 
\beq\label{eq:ndiff}
n_2-n_1=\frac{1}{2\theta}(\tan^{-1}(\tau_1)-\tan^{-1}(\tau_2))
  =\frac{1}{2\sin^{-1}(\lambda/2)}\tan^{-1}\left(\frac{\tau_1-\tau_2}{1+\tau_1\tau_2}\right).
\eeq
From lemma \ref{lemma:tauP} we have
$$
\begin{array}{ll}
\displaystyle\tau(P_1(m))-\tau(P_0(m))=\frac{2}{4 m^2-1} + O(\lambda),\\ \\
\displaystyle\tau(P_1(m-1))-\tau(P_0(m-1))=\frac{2}{4 m^2 -8 m+3} +O(\lambda)\\ \\
\displaystyle\frac{3}{5}\lambda < \tau(P_0(m-1))-\tau(P_1(m)) <\frac{17}{10}\lambda.
\end{array}
$$
Inserting these estimates in (\ref{eq:ndiff}) and considering the small-argument 
behaviour of the inverse trigonometric functions, we obtain the desired relations. 

\subsection*{A.6\quad  Proof of lemma \ref{lemma:Zmnexpand}}
The method of proof is a familiar one.  The remainder factors $r_{Zx}$ and $r_{Zy}$ are 
constructed by combining (\ref{eq:Zmnformula}) and (\ref{eq:Zmnexpand}) with the 
insertion of $P_0(m),P_1(m-1),Q_0(n-1)$, and $Q_1(n)$ from lemma \ref{lemma:Zmnexpand}.  
Here care must be taken to express $Q_0(n-1)$ in terms of the variable $\tau=\tau(Q_0(n))$:
$$
\tau(Q_0(n-1))=\tan\left(\frac{\pi}{2}-2 n \theta + 2\theta\right)
  =\frac{\tau + \tan (2\theta)}{1-\tau\tan(2\theta)}.
$$
The resulting expressions (see \cite{ESupplement} for details) for $r_{Zx}$ and $r_{Zy}$ are 
enormous polynomials in $\lambda, m, \tau$, and the various remainder factors appearing in the 
expansion formulae for the atom vertices.  Because the $Z(m,n)$ lies within a single atom 
$\Lambda_m^{\rm out}$, we improve our estimates by replacing $\tau$ by $r/m$, where, from 
appendix A.5, $r\in [2/5, 1001/1000]$. Application of the rules of interval arithmetic 
(appendix A.1) then yields the uniform bounds stated in the lemma.

\subsection*{A.7\quad Crossover calculation}
We now return to the crossover phenomenon for a more precise treatment. 
First of all, we observe that the angle between the line segments $[P_0(m),P_1(m)]$ 
and $[Q_0(n(P_0(m))$, $Q_1(n(P_0(m))]$ changes sign as one proceeds from right to left in $\Lambda$.  
For some $m=m^*, n=n^*=n(P_0(m^*))$, the segments coincide and the atoms 
$\Lambda_{m^*}^{\rm out}$ and $\Lambda_{n^*}^{\rm in}$ have their vertical edges nearly parallel. 
We adopt $m=m^*$ as our precise definition of the crossover.  
By the construction of section \ref{section:ReturnMap}, the orientation of the line 
through $[P_0(m),P_1(m)]$ differs from the vertical by a generalized rotation 
$C^{-4m}$, while $[Q_0(n), Q_1(n)]$ is rotated by $C^{4n+3}$.  
Thus crossover corresponds to
$$
C^{4m^* +4 n^* +3} =\pm {\bf 1}.
$$
Obviously there are infinitely many solutions, but the relevant one, corresponding to a 
single crossover as $m$ varies from 1 to $M$, has a negative sign, with $m^*$ satisfying
$$
(4m^* +4 n(P_0(m^*)) +3 )\sin^{-1}\left(\frac{\lambda}{2}\right)=\frac{\pi}{2},
$$ 
where $n(P_0(m*))$ is calculated by solving $Q_0(n)=X$ for $n$, then setting $X=P_0(m*)$.  
The equation for $m*$ can be solved perturbatively to obtain
\begin{eqnarray}
m^*(\lambda)&=& \frac{1}{\sqrt{2\lambda}}\left(1+\frac{\lambda}{3} +\frac{\lambda^2}{120} 
+\frac{31\lambda^3}{1008}+O(\lambda^4)\right), \label{eq:mstar}
\\
n^*(\lambda)&=&\frac{\pi}{4\lambda}-\frac{1}{\sqrt{2\lambda}}-\frac{3}{4}
-\frac{\sqrt{\lambda}}{3\sqrt{2}}-\frac{\pi\lambda}{96}
-\frac{\lambda^{3/2}}{120\sqrt{2}} +O(\lambda^2)\nonumber
\end{eqnarray}
In particular, $m^*(\lambda)$ diverges to infinity as $\lambda\rightarrow 0$.

\subsection*{A.8\quad Proof of $\Lambda_m^{\rm out}$ area estimate}
We are seeking uniform bounds, for  $3\leq m \leq M$ and $0<\theta\leq10^{-4}$ on 
$$
r(\theta,m)=\frac{\theta^{-2} m^{-2} c((2m-1)\theta)}{s((2m-3)\theta)s((2m-1)\theta)s((2m+1)\theta)}-1, 
$$
where
$$
c(x)=\cos(x),\qquad s(x)=\frac{\sin x}{x}.
$$
Note that in calculation of bounds we will not be able to benefit from the 
previous cut-off on $m$. 

Following our usual practice, we define auxiliary variables $u$ and $v$ via 
$$
c(x)=1-\frac{x^2}{2}+\frac{(1-u)x^4}{24},\qquad
s(x)=1-\frac{x^2}{6}+\frac{(1-v)x^4}{120}.
$$
Bounds $0\leq u, v \leq 1002/1000$ are straightforward to establish \cite{ESupplement} (sec. E.11.1).
Insertion of these expressions reduces $r(\theta,m)$ to a ratio of polynomials 
in the variables $\theta, m, u$, and $v$.  Although each polynomial can be 
uniformly bounded above and below using the same methods used throughout this 
article, there is a serious problem: the denominator range includes zero.  
On the other hand, numerical evaluation of the denominator shows that it is 
almost certainly negative definite.  

To proceed, we isolate the leading contribution consisting of terms 
proportional to $m^k \theta^k$  (for all others, the power of $m$ is less 
than that of $\theta$). There is no $u$ dependence of these terms. 
With respect to $v$, we make the worst-case assumption, setting it equal 
to zero for all negative terms and to its maximum value of $1002/1000$.  
The leading part of the denominator then reduces to
$$
{\cal D}_1(\theta,m)= h(m^2 \theta^2),
$$
\begin{eqnarray*}
h(x)&=&-1728000 + 3456000 x - \frac{11513088}{5} x^2 + 1433600 x^3\\
&& +\frac{2328576}{25}x^4 + \frac{384768768}{3125}x^5 + \frac{24096096064}{1953125}x^6.
\end{eqnarray*}
We need to show that this polynomial is negative definite over the 
range $0<x<\frac{\pi^2}{16}< x_1=\frac{6171}{10000}$.  
Direct application of interval arithmetic fails, as expected, producing a range 
of values containing zero.  On the other hand, the same method, applied to the 
derivative $h'(x)$ produces a pair of positive bounds.Thus $h(x)$ is a monotone increasing function.  
Since it is negative at both $0$ and $x_1$ it must be negative everywhere on 
the interval. Specifically,
$$
{\cal D}_1(\theta,m)\in [-1728000,-110000]
$$
Meanwhile,the numerator ${\cal N}(\theta,m,u,v)$ and the non-leading part of the 
denominator, ${\cal D}_2(\theta,m,u,v)$ can be bounded by our original method, 
obtaining \cite{ESupplement} (sec. E.11.1)
$$
{\cal N}(\theta,m,u,v) \in [-24730000, 27090000],   \qquad {\cal D}_2((\theta,m,u,v)\in [-31800000, 14920].
$$

Combining the bounding intervals according to the rules of appendix A.1, we get, finally
$$
r(\theta,m)\in\left[-285,\frac{1301}{5}\right].
$$
\endproof

\end{document}